\RequirePackage{fix-cm}
\documentclass[smallextended,envcountsame,numbook]{svjour3-ArXiv}      % onecolumn (second format)
\smartqed  % flush right qed marks, e.g. at end of proof
\usepackage{tikz-cd}
\usetikzlibrary{arrows.meta, positioning, calc}
\usetikzlibrary{shapes.misc}
\tikzset{
  nodebox/.style={draw, rounded corners=2pt, inner sep=3pt, minimum size=6mm}
}

\usepackage{appendix}

\usepackage{graphicx}
\usepackage{slashed}
\usepackage{arydshln}

 \usepackage{hyperref}
 \hypersetup{
     colorlinks=true,
     linkcolor=blue,
     filecolor=magenta,      
     urlcolor=cyan,
     citecolor=red,
}
\usepackage{amsfonts, amssymb, amsbsy, amsmath}
 \usepackage{blkarray}
 \usepackage{graphicx}
 \usepackage{float}
 \usepackage{epstopdf}
 \usepackage{algorithmic}
 \ifpdf
   \DeclareGraphicsExtensions{.eps,.pdf,.png,.jpg}
 \else
   \DeclareGraphicsExtensions{.eps}
 \fi
 \usepackage{mathrsfs}
 \usepackage{amscd,empheq} %,yhmath}
 \usepackage{yfonts,bm}

 \usepackage[pdf]{graphviz}

 \usepackage[all,2cell]{xy}
 \usepackage{pb-diagram}
 \usepackage{pb-xy}

 \usepackage{tikz}
 \usepackage{tikz-cd}
 \usetikzlibrary{patterns}
 \usetikzlibrary{positioning}
 \usetikzlibrary{shapes,decorations.markings,arrows.meta}
 \usetikzlibrary{matrix,arrows}
 \tikzset{hfit/.style={rounded rectangle, inner xsep=0pt},
           vfit/.style={rounded corners}}

 \usepackage{morewrites}
 \usepackage{scalerel}
 \usepackage{mathabx}
 \usepackage{palatino}
 
\usepackage{enumitem}
\setlist[enumerate]{leftmargin=.5in}
\setlist[itemize]{leftmargin=.5in}

\usepackage{amsopn}

\newcommand{\rmGamma}{\mathrm{\Gamma}}

\newcommand{\dff}{{\mathrm d}}

\newcommand{\sfS}{{\mathsf \Sigma}}
\newcommand{\sfL}{{\mathsf \Lambda}}

\newcommand{\bt}{{\bm 0}}
\newcommand{\smbt}{{\scalebox{0.6}{$\bm{0}$}}}
\newcommand{\tp}{{\bm 1}}
\newcommand{\smtp}{{\scalebox{0.6}{$\bm{1}$}}}

\newcommand{\spi}{\mathsf{p}}

\newcommand{\IIonei}{\bfQ_{\mathbf 1}}
\newcommand{\IIi}{\bfQ_{\mathbf 2}}

\newcommand{\IIthreei}{\bfQ_{\mathbf 3}}

\newcommand{\IIfouri}{\bfQ_{\mathbf 4}}
\newcommand{\IIfivei}{\bfQ_{\mathbf 5}}
\newcommand{\IIsixi}{\bfQ_{\mathbf 6}}
\newcommand{\IIeighti}{\bfQ_{\mathbf 8}}
\newcommand{\smin}{\setminus}

\newcommand{\pred}{\blacktriangleleft}
\renewcommand{\subset}{\subseteq}

\newcommand{\rmTheta}{\mathrm{\Theta}}

\definecolor{Blue}{RGB}{0,162,255}
\definecolor{Orange}{RGB}{243,144,25}
\definecolor{Red}{RGB}{236,93,87}

\newcommand{\bR}{{R}}

\newcommand{\K}{{\mathbb{K}}}
\newcommand{\N}{{\mathbb{N}}}

\newcommand{\Z}{{\mathbb{Z}}}

\newcommand{\sC}{{\mathsf C}}

\newcommand{\sE}{{\mathsf E}}

\newcommand{\sO}{{\mathsf O}}

\newcommand{\sM}{{\mathsf M}}

\newcommand{\sP}{{\mathsf P}}
\newcommand{\sQ}{{\mathsf Q}}

\newcommand{\sS}{{\mathsf S}}

\newcommand{\sU}{{\mathsf U}}

\newcommand{\sk}{{\mathsf k}}

\renewcommand{\imath}{{\mathsf i}}
\renewcommand{\jmath}{{\mathsf j}}

\newcommand{\Hom}{\mathrm{Hom}}

\newcommand{\id}{\mathrm{id}}
\newcommand{\coker}{\mathrm{coker}}

\DeclareMathOperator{\Inv}{\mathrm{Inv}}

\newcommand{\image}{\mathrm{im~}}

\newcommand{\cone}{\mathrm{Cone}}

\newcommand{\Sub}{\mathsf{Sub}}

\newcommand{\sRmod}{\mathtt{Mod}}
\newcommand{\sRmodR}{\mathtt{Mod}_{R}}
\newcommand{\sRmodRQ}{\mathtt{Mod}_{R,\shQ}}
\newcommand{\sRmodRQd}{\mathtt{Mod}_{R,\shQ^{\rm diff}}}
\newcommand{\CdR}{\mathtt{Cd}(\bR)}
\newcommand{\KdR}{\mathtt{Kd}(\bR)}
\newcommand{\DdR}{\mathtt{Dd}(\bR)}
\newcommand{\DQ}{\mathscr{D}_{\shQ^{\rm diff}}(\bR)}
\newcommand{\QS}{\mathtt{Qis}}

\newcommand{\shQ}{\mathtt{Q}}
\newcommand{\bfE}{\text{\bf E}}

\newcommand{\bfC}{\text{\bf C}}
\newcommand{\bfS}{\text{\bf S}}
\newcommand{\bfQ}{\text{\bf Q}}

\newcommand{\Cd}{\mathtt{Cd}}

\newcommand{\Eexact}{\mathscr E}

{\begin{snugshade}\begin{quote}}
{\hfill\end{quote}\end{snugshade}}

\definecolor{shadecolor}{rgb}{0.8,0.8,0.8}

\begin{document}

\title{Representations of spectral systems over finite distributive lattices}
%  %\thanks{Grants or other notes
% %about the article that should go on the front page should be
% %placed here. General acknowledgments should be placed at the end of the article.}
% }

\titlerunning{Representations of spectral systems}
%\titlerunning{Spectral representations of Cartan-Eilenberg systems}
%\titlerunning{Tessellated Morse decompositions and connection matrices revisited}        
% if too long for running head

\author{Robert Vandervorst %etc.
}


\institute{
R.C.A.M. Vandervorst \at
              %Department of Mathematics\\
              VU University Amsterdam\\
              Amsterdam, The Netherlands\\
              \email{r.c.a.m.vander.vorst@vu.nl}
}

\date{Version date: \today}

\maketitle

\begin{abstract}
A representation theorem is proved for spectral systems indexed by the lattice of down-sets of a finite poset $\sP$. Over an arbitrary associative ring, every excisive, homotopy-exact chain system of differential modules   is 
isomorphic to the chain system of a projective $\sP$-graded differential module
 with a single upper-triangular differential. Consequently, every Cartan--Eilenberg system that is excisive and chain generated, together with its spectral system, admits such a differential model. An example over $\Z/4$ shows that chain generation cannot be dropped. 
 A chain system over $\Z/8$ with short exact rows, whose maps commute only up to homotopy, shows that homotopy-exactness cannot be weakened to distinguished triangles in the derived category.
 Under finiteness conditions the model is finitely generated. The proof combines a derived horseshoe lemma with induction over the lattice, and yields a short proof and extension of Franzosa's connection matrix theory.

\end{abstract}

\tableofcontents

\section{Prelude}
\label{prelude}
Spectral sequences are presented page by page, where each page is the homology of the previous one. Already for filtrations indexed by the integers the bookkeeping is intricate. For filtrations indexed by a finite distributive lattice, represented as the lattice $\sO(\sP)$ of down-sets\footnote{By Birkhoff's representation theorem every finite distributive lattice is of this form, with $\sP$ isomorphic to the poset of its join-irreducible elements, cf.\ \cite[Ch.~5]{davey:priestley}.} of a finite poset $\sP$, the pages are replaced by a \emph{spectral system} \cite{Matschke}. The latter is defined as a collection of modules $S^{\gamma\delta}_{\alpha\beta}$ over an \emph{associative unital ring $\bR$}, indexed by ordered quadruples $\alpha\subseteq\beta\subseteq\gamma\subseteq\delta$ in $\sO(\sP)$, connected by differentials and exact sequences, cf.\ Sect.\ \ref{specsystems}. This paper shows that, under natural and sharp hypotheses, this information is carried by a single   projective $\sP$-graded $\bR$-module $A$ with an upper-triangular differential $\dff_A$. The pair $(A,\dff_A)$ is called a \emph{$\sP$-graded differential module}, cf.\ Sects.\ \ref{diffmodules} and \ref{chainsystems}.

Let $A=\bigoplus_{p\in\sP}G_pA$ be a $\sP$-graded $\bR$-module with differential $\dff_A$, $\dff_A^2=0$, whose matrix entries satisfy
\[
\dff_A(p,q)\colon G_qA\to G_pA,\quad \dff_A(p,q)\neq 0\implies p\le q,\quad\text{e.g.}~~
\dff_A=\begin{pmatrix}
\dff_{11} & \dff_{12} & \dff_{13}\\ 0 & \dff_{22} & \dff_{23}\\ 0 & 0 & \dff_{33}
\end{pmatrix},
\]
cf.\ \eqref{diffasmatrix}, where the example is for the poset $\sP=\{1<2<3\}$. For down-sets $\alpha\subseteq\beta$ the block of $\dff_A$ with rows and columns in the convex set $\beta\smin\alpha$ is a differential on $G_{\beta\smin\alpha}A=\bigoplus_{p\in\beta\smin\alpha}G_pA$. For $\alpha\subseteq\beta\subseteq\gamma$ the short exact sequences $0\to G_{\beta\smin\alpha}A\to G_{\gamma\smin\alpha}A\to G_{\gamma\smin\beta}A\to 0$ yield exact triangles in homology, and for $(\alpha,\beta,\gamma)\le(\alpha',\beta',\gamma')$ commutative squares
\[
\begin{tikzcd}[column sep=small]
E^\beta_\alpha \arrow[rr,"i"] & & E^\gamma_\alpha \arrow[dl,"j"]\\
& E^\gamma_\beta \arrow[ul,"k"] &
\end{tikzcd}
\qquad\qquad
\begin{tikzcd}
E^\gamma_\beta \arrow[r,"k"] \arrow[d,"\ell"'] & E^\beta_\alpha \arrow[d,"\ell"]\\
E^{\gamma'}_{\beta'} \arrow[r,"k'"] & E^{\beta'}_{\alpha'}
\end{tikzcd}
\]
where $E^\beta_\alpha=H(G_{\beta\smin\alpha}A)$, and likewise for $E^\gamma_\alpha$ and $E^\gamma_\beta$. These data form a \emph{Cartan--Eilenberg system} $\bfE(A,\dff_A)$, cf.\ Defn.\ \ref{CEsys} and \cite[Ch.~XV, \S7]{CE}, \cite{HelleRognes,SpV4}, which in turn determines a spectral system $\bfS(A,\dff_A)$ of consecutive homology pages, cf.\ \eqref{Sterms}. The diagonal blocks $\dff_A(p,p)$ compute the \emph{basic $E$-terms} $E^\beta_\alpha=H\bigl(G_pA,\dff_A(p,p)\bigr)$, $\beta\smin\alpha=\{p\}$, which depend only on $p$.
The off-diagonal blocks record the connecting homomorphisms, the differentials on all pages of the spectral system, and the extensions between them. For an abstract Cartan--Eilenberg system $\bfE$, a (projective) $\sP$-graded differential module $(A,\dff_A)$ with $\bfE(A,\dff_A)\cong\bfE$ is called a \emph{(projective) differential model} for $\bfE$, cf.\ Sect.~\ref{specsystems}. Such a model compresses the entire system into one matrix.

Systems that come from chains are described by \emph{chain systems} $\bfC$, cf.\ Defn.\ \ref{Csys} and \cite{SpV4}, which consist of differential modules $C^\beta_\alpha$, for $\alpha\subseteq\beta$ in $\sO(\sP)$, with distinguished triangles $C^\beta_\alpha\to C^\gamma_\alpha\to C^\gamma_\beta\to\sfS C^\beta_\alpha$ in the derived category $\DdR$ of differential modules, cf.\ Sects.\ \ref{diffmodules} and \ref{triang}. In this paper we consider chain systems that are \emph{homotopy-exact}, cf.\ Defn.\ \ref{hexrows}. For such systems the maps between the modules $C^\beta_\alpha$ are chain maps that compose as chain maps, i.e.\ the system is a functor to $\CdR$, and not only to $\DdR$, and every row $C^\beta_\alpha\xrightarrow{\,i\,}C^\gamma_\alpha\xrightarrow{\,j\,}C^\gamma_\beta$ is homotopy-exact, which entails that the map $\rmTheta\colon\cone(i)\to C^\gamma_\beta$, $\rmTheta(y,x)=j(y)$, is a quasi-isomorphism, cf.\ Defn.\ \ref{homotopyexact}. Homotopy-exact sequences yield distinguished triangles in $\DdR$. A chain system is \emph{excisive} if the maps $\ell\colon C^\alpha_{\alpha\cap\beta}\to C^{\alpha\cup\beta}_\beta$ are isomorphisms in $\DdR$ for all $\alpha,\beta\in\sO(\sP)$. Chain systems with short exact rows, or with the weakly exact rows of \cite{fran}, are homotopy-exact. Lattice-filtered differential modules, in particular $\sP$-graded ones, yield excisive, homotopy-exact chain systems, cf.\ \eqref{sheactno2}-\eqref{dirsumpgm}, and so do the relative chain complexes of filtered spaces under the hypotheses of the excision theorem, e.g.\ for index filtrations \cite{fran86} and filtrations by subcomplexes, cf.\ Sect.\ \ref{chainsystems}. A Cartan--Eilenberg system $\bfE$ is \emph{chain generated} if $\bfE\cong\bfE(\bfC)$, the homology of a homotopy-exact chain system $\bfC$, cf.\ Defn.\ \ref{def:chaingenerated}.
\begin{theorem}[cf.\ Thms.\ \ref{mainres987655} and \ref{thm:chaingeneratedmodel}]
\label{preludethm}
Let $\bR$ be an associative unital ring and $\sP$ a finite poset. Every excisive, homotopy-exact chain system $\bfC$ over $\sO(\sP)$ is isomorphic to the chain system $\bfC(A,\dff_A)$ of a projective $\sP$-graded differential module $(A,\dff_A)$, via $K$-projective resolutions $G_{\beta\smin\alpha}A\to C^\beta_\alpha$, cf.\ Defn.~\ref{cofibrantrepl}, that commute with the chain maps $i$ and $j$ in $\CdR$. Consequently, every excisive, chain generated Cartan--Eilenberg system $\bfE$, and its spectral system $\bfS$, have a projective differential model: $\bfE\cong\bfE(A,\dff_A)$ and $\bfS\cong\bfS(A,\dff_A)$.
\end{theorem}

Conversely, $\bfC(A,\dff_A)$ is an excisive chain model for $\bfE(A,\dff_A)$. Its rows are short exact and its excision maps are identities, since $\alpha\smin(\alpha\cap\beta)=(\alpha\cup\beta)\smin\beta$. The hypotheses are therefore \emph{sharp}. A Cartan--Eilenberg system has a differential model if and only if it is excisive and chain generated. In particular, it allows a projective differential model. Not every Cartan--Eilenberg system is chain generated. Over $\bR=\Z/4$ a single exact triangle with $i$, $j$ and $k$, given by multiplication by $2$, admits no differential model, cf.\ Sect.\ \ref{counterex12}. Whether the Cartan--Eilenberg system of every excisive chain system in $\DdR$ is chain generated is open, cf.\ Rmk.\ \ref{notcgsys}.

The result is about \emph{differential modules} over an arbitrary ring $\bR$. We assume no homological grading  and chain complexes are a special case.
In general $G_pA$ is a $K$-projective resolution of the differential modules $C^\beta_\alpha$, $\beta\smin\alpha=\{p\}$. Moreover, $\dff_A(p,p)$ need not vanish. For an excisive Cartan--Eilenberg system, the basic $E$-terms $E^\beta_\alpha$, $\beta\smin\alpha=\{p\}$, are isomorphic via excision maps. If the basic $E$-terms are projective, there is a \emph{strict} model, i.e.\ $\dff_A(p,p)=0$ and $G_pA\cong E^\beta_\alpha$, $\beta\smin\alpha=\{p\}$ for all $p\in\sP$, cf.\ \eqref{strictdiff12}, Thm.~\ref{mainres987655proj} and Cor.~\ref{mainres987655-22}.

For most computational purposes differential models should be finitely generated. Section~\ref{specialcases} gives finitely generated models under conditions on the basic $E$-terms. For example, if $\bR$ is hereditary and Noetherian, e.g.\ a principal ideal domain, and the basic $E$-terms are finitely generated, then $G_pA=P_0\oplus P_1$ is the differential module of a projective resolution $0\to P_1\to P_0\to E^\beta_\alpha\to 0$, $\beta\smin\alpha=\{p\}$, by finitely generated modules, cf.\ Prop.~\ref{smallresol12}, Thm.~\ref{mainres987655chain33} and Cor.~\ref{mainres987655-33}. For arbitrary associative unital rings $\bR$, if the basic $E$-terms admit finite resolutions by finitely generated projective modules, then there is a finitely generated projective differential model, cf.\ Thm.~\ref{mainres987655chain44} and Cor.~\ref{mainres987655-44}.

The main motivation for Theorem~\ref{preludethm} comes from Conley index theory \cite{conley,fran}. For a Morse representation of a flow, or a map, indexed by a finite poset $\sP$, cf.\ \cite{lsa,robbin:salamon2}, the relative chain complexes of an index filtration \cite{fran86} form an excisive, homotopy-exact chain system with short exact rows, and the basic $E$-terms are the homology Conley indices of the Morse sets. A differential model is then a \emph{connection matrix}, or \emph{$C$-connection matrix}, in the sense of \cite{fran}. The off-diagonal entries $\dff_A(p,q)$, $p<q$, carry information about connecting orbits between the Morse sets labeled by $q$ and $p$. Franzosa \cite{fran} proves existence when $\bR$ is a principal ideal domain. If the basic $E$-terms (Conley indices) are free, a strict differential (connection matrix) exists, and if the basic $E$-terms are not free, resolutions of the Conley indices yield a $C$-connection matrix, cf.\ Rmk.~\ref{franzosaresult}. For a computational approach, cf.\ \cite{hms}. Theorem~\ref{preludethm} provides a short proof in terms of homological algebra, and extends the theory to arbitrary coefficient rings $\bR$ and to differential modules.

The extension to arbitrary rings is not only of algebraic interest. For a dynamical system given by a map $\bm f$,   the homological Conley index of an isolated invariant set is
the shift equivalence class of the map induced by $\bm f$ on the homology of an index pair,
cf.\ \cite{Mrozek90,FranksRich}. Let $\sU=\{F_\alpha U\}$, $\alpha\in\sO(\sP)$, be
a lattice of
trapping regions for a Morse representation of $\bm f$, cf.\ \cite{lsa,FranksRich}.
For every index pair $(F_\beta U,F_\alpha U)$, $\alpha\subset\beta$, the relative chain complexes $C(F_\beta U,F_\alpha U)$ over the polynomial ring $\Z[z]$, with $z$ acting by $\bm f$ on chains,
define the Conley indices.
This yields
a chain system over $\bR=\Z[z]$, whose basic $E$-terms are the homology Conley indices of the Morse sets $M_{\beta\setminus\alpha}=\Inv(F_\beta U,F_\alpha U)$\footnote{The notation $\Inv(F_\beta U,F_\alpha U)$ denotes the maximal invariant sets in $F_\beta U\setminus F_\alpha U$ with respect to $\bm f$.} together with their index maps.
A differential model over $\Z[z]$ is then a $C$-connection matrix that also records the index maps. %, cf.\
With coefficients in a field $\K$ the ring $\K[z]$ is a principal ideal domain, cf.\ \cite{Rich2}. However, with integer coefficients, which detect torsion in the Conley indices, the ring $\Z[z]$ is neither a principal ideal domain, nor hereditary.
Chain complexes over $\Z[z]$ need not be formal, and the approaches of \cite{fran,Rich2} do not apply. In this  case Theorem~\ref{preludethm} is applicable and extends the connection matrix theory to arbitrary rings.
In Section \ref{dynexample} we provide a detailed example; a more in-depth analysis
will be the subject of future study.

In Franzosa's approach the restrictions on the ring $\bR$ stem from the method of proof. The analogue of the horseshoe lemma in \cite[Prop.~4.7]{fran} uses a basis of the cycles of a free complex and a complementary subcomplex. Both exist over principal ideal domains; over hereditary rings the cycles are projective, not necessarily free, and a complementary subcomplex still exists. Theorem~\ref{SESinDR} replaces the element-wise construction by the lifting property of $K$-projective differential modules up to homotopy, cf.\ Prop.~\ref{prop:kprojequiv}, which holds over every ring, as does the existence of $K$-projective resolutions, cf.\ App.~\ref{resolfordiffmod}. The idea in \cite{fran} is to replace each chain complex $C^\beta_\alpha$, $\beta\smin\alpha=\{p\}$, by a free resolution of length one of its homology, with the differential of the resolution as diagonal entry $\dff_A(p,p)$. This presupposes that $C^\beta_\alpha$ is determined by its homology up to quasi-isomorphism, i.e.\ $C^\beta_\alpha\sim\bigl(H(C^\beta_\alpha),0\bigr)$ in $\DdR$ (formality), which holds over hereditary rings, cf.\ Lem.~\ref{heredinderived}, but not in general. In our approach, Step~(ii) of the proof of Theorem~\ref{mainres987655} chooses $K$-projective resolutions of the differential modules $C^\beta_\alpha$, $\beta\smin\alpha=\{p\}$, as opposed to resolutions of their homology. Resolutions built from homology are only used in special cases, cf.\ Sect.~\ref{specialcases}. Another difference is that in \cite{fran} the chain complexes are indexed by convex sets of $\sP$, so that excision holds by construction, whereas here a chain system is indexed by pairs of down-sets and the excision maps are only required to be quasi-isomorphisms.

Sections~\ref{diffmodules} and~\ref{triang} are preliminaries about differential modules, homotopy categories and derived categories, cf.\ \cite{ABI2007,gelfand,HolmJorg1,MacLane,Spalt,weibel}. Since the literature mostly treats chain complexes, the standard facts are proved for differential modules where needed, cf.\ Apps.~\ref{rottrian} and~\ref{resolfordiffmod}. Chain systems are defined in Section~\ref{chainsystems}, Cartan--Eilenberg systems and spectral systems in Section~\ref{specsystems}. The main result is also sharp in the sense that homotopy-exactness cannot be weakened to the requirement of Definition~\ref{Csys} that the rows are distinguished triangles in $\DdR$. In Section~\ref{csysnothe} we give a counterexample with $\sP=\{1<2<3\}$, over commutative rings $\bR$ with an element $z$ such that $z^3=0$ and $z^2\neq 0$, e.g.\ $\Z/8$, $\Z/p^3$ and $\K[z]/(z^3)$. The counterexample starts with an appropriately chosen $\sP$-graded differential module $(A,\dff_A)$. Transporting the row of one triple of $\bfC(A,\dff_A)$ along an automorphism of its middle term yields a chain system $\bfC'$ in $\DdR$ with the same terms,  base\footnote{The terms $C^\beta_\smbt$ over the empty down-set $\bt$, with the maps $i$ between them.} and bottom triangles,\footnote{The rows of the triples $(\bt,\alpha,\beta)$} cf.\ Sect.~\ref{csysnothe}, whose rows are still short exact, but whose maps compose only up to homotopy. It is then shown that $\bfC'$ is \emph{not} isomorphic to $\bfC(B,\dff_{B})$ for any $\sP$-graded differential module $(B,\dff_{B})$, and hence not to any homotopy-exact chain system. Nevertheless, $\bfE(\bfC')\cong\bfE(A,\dff_A)$. Sections~\ref{chainsystems}--\ref{specialcases} can be read with Section~\ref{diffmodules}, Definition~\ref{homotopyexact} and the definitions of Section~\ref{chainsystems}. The proof of Theorem~\ref{mainres987655} uses Sections~\ref{triang} and~\ref{hslfordiffmod}, and Appendices~\ref{rottrian} and~\ref{resolfordiffmod}.

\section{Differential modules}
\label{diffmodules}
A \emph{differential module} $(A,\dff_A)$ consists of a left $\bR$-module $A$ over an associative unital ring $\bR$ and an endomorphism (differential) $\dff_A\colon A\to A$ with $\dff_A^2=0$, cf.\ \cite[II.1]{MacLane}.
The \emph{homology} of $(A,\dff_A)$ is defined as $H(A) := \ker\dff_A/\image \dff_A$.
A homomorphism $f\colon (A,\dff_A) \to (B,\dff_B)$ of differential modules is a homomorphism of $R$-modules that
satisfies $\dff_B f = f \dff_A$ and is called a
\emph{chain homomorphism}, or \emph{chain map}.
Differential modules and chain maps form the category of differential modules $\CdR$.
Given a chain map $f\colon (A,\dff_A)\to (B,\dff_B)$, the homology homomorphism $H(f)\colon H(A)\to H(B)$ is defined by
$H(f)[a] = [f(a)]$, for all $a\in \ker \dff_A$. Homology defines a functor $H\colon \CdR \to \sRmodR$ to the category $\sRmodR$ of left $\bR$-modules.

Two chain maps $f,g\colon A\to B$ are \emph{chain-homotopic} if there exists a homomorphism $h\colon A\to B$ such that $f-g= h\dff_A+\dff_B h$. Notation $f\pmb{\sim} g$. Homotopic chain maps induce the same map in homology.
Two chain homomorphisms $f\colon A\to B$ and $g\colon B\to A$ with $gf\pmb{\sim} \id_A$ and $fg\pmb{\sim} \id_B$ are called \emph{chain-equivalences}, and $A$ and $B$ are \emph{chain-equivalent}, notation $A\simeq B$.
The \emph{homotopy category} $\KdR$ has the same objects as $\CdR$, and the morphisms $A\to B$ in $\KdR$ are the homotopy classes of chain maps; homotopy is compatible with composition and sums, and chain-equivalences become isomorphisms.
A chain map $f\colon A\to B$ for which $H(f)\colon H(A)\to H(B)$ is an isomorphism is called a \emph{quasi-isomorphism}, or \emph{qis}.
The \emph{derived category} $\DdR$ is the localization $\DdR = \KdR[\QS^{-1}]$, where $\QS$ is the class of quasi-isomorphisms. Both $\KdR$ and $\DdR$ are additive categories.
Since $\QS$ is a localizing class in $\KdR$, cf.\ \cite[Ch.~III]{gelfand} and Rmk.~\ref{Qshapedanalogon}, morphisms in $\DdR$ are represented
as \emph{roof diagrams} given by:
\begin{equation}
    \label{roof12}
\begin{tikzcd}
& W \arrow[dl, "s"'] \arrow[dr, "w"] & \\
A & & B
\end{tikzcd}
\end{equation}
where $s$ is a quasi-isomorphism, and $w$ is a morphism in the homotopy category $\KdR$, i.e.\ a homotopy class of chain maps, cf.\ \cite[Lem.\ 8, p.\ 148]{gelfand}. We also write a roof as a pair $(w,s)$, or as the \emph{fraction} $w s^{-1} = w/s$.
The \emph{localization functor} $\KdR\to \DdR$, defined by $f\mapsto f/\id$, makes quasi-isomorphisms invertible and is an additive functor.
Two roofs
\[
\begin{tikzcd}
& W \arrow[dl, "s"'] \arrow[dr, "w"] & \\
A & & B
\end{tikzcd}
\quad\quad \approx \quad\quad
\begin{tikzcd}
& W' \arrow[dl, "s'"'] \arrow[dr, "w'"] & \\
A & & B
\end{tikzcd}
\]
are equivalent, $(w,s)\approx(w',s')$, if there exists a commutative diagram in $\KdR$
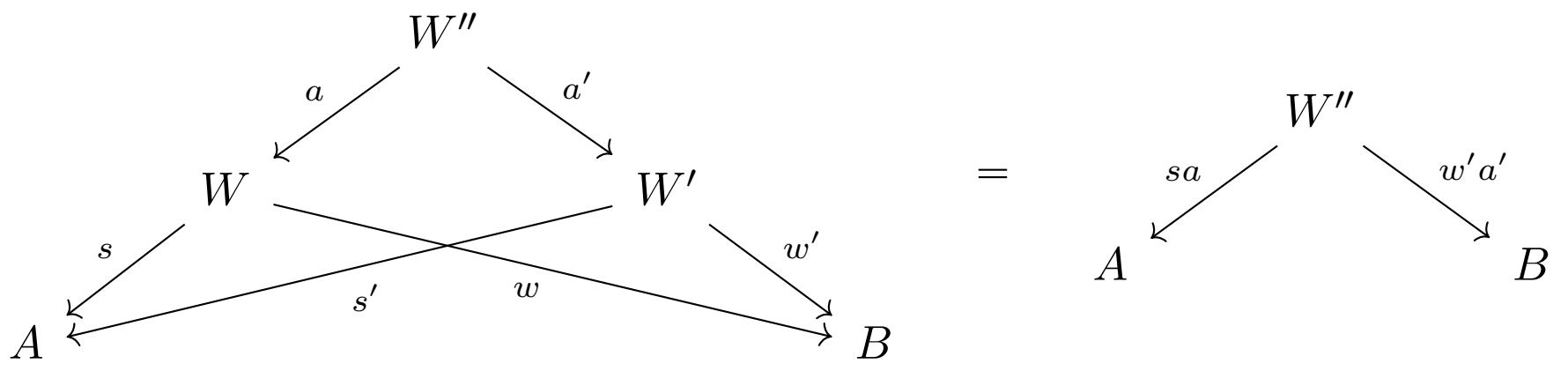
\[
\begin{tikzcd}
  &                                      & W'' \arrow[ld, "a"'] \arrow[rd, "a'"] &                                         &   \\
  & W \arrow[ld, "s"'] \arrow[rrrd, "w"'] &                                          & W' \arrow[rd, "w'"] \arrow[llld, "s'"] &   \\
A &                                      &                                          &                                         & B
\end{tikzcd}
\quad = \quad
\begin{tikzcd}
& W'' \arrow[dl, "sa"'] \arrow[dr, "w'a'"] & \\
A & & B
\end{tikzcd}
\]
with $a$ a quasi-isomorphism, i.e.\ $sa=s'a'$ and $wa=w'a'$, cf.\ \cite[Lem.\ 8, p.\ 148]{gelfand}.
For two roofs $f=w/s\colon A\to B$ and $f'=w'/s'\colon B\to C$ there exists a diagram
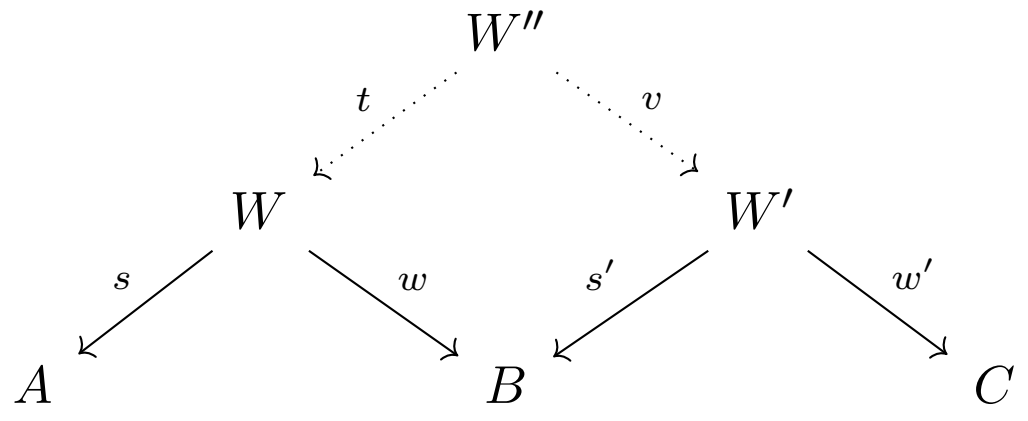
\[
\begin{tikzcd}
  &                                    & W'' \arrow[ld, "t"', dotted] \arrow[rd, "v", dotted] &                                       &   \\
  & W \arrow[ld, "s"'] \arrow[rd, "w"] &                                                      & W' \arrow[rd, "w'"] \arrow[ld, "s'"'] &   \\
A &                                    & B                                                    &                                       & C
\end{tikzcd}
\]
with $t$ a quasi-isomorphism and $wt=s'v$, and the composition is given by  $f'f=(w'v)/(st)\colon A\to C$.
Homology defines a functor $H\colon\DdR\to\sRmodR$ by $H(w/s)=H(w)H(s)^{-1}$.
An isomorphism in $\CdR$ is denoted by $\cong$, an isomorphism in $\KdR$, i.e.\ a chain-equivalence, by $\simeq$, and an isomorphism in the derived category $\DdR$ by $\sim$.
\begin{lemma}
    \label{isoisqis}
    A morphism $f=w/s\colon A\to B$ in the derived category $\DdR$ is an isomorphism if and only if $w\colon W\to B$ is a quasi-isomorphism, if and only if $H(f)$ is an isomorphism.
\end{lemma}
\begin{proof}
If $f$ is an isomorphism in $\DdR$, then $H(f)$ is an isomorphism. Since $H(f) = H(w)H(s)^{-1}$ and $H(s)$ is an isomorphism, $H(f)$ is an isomorphism if and only if $H(w)$ is, i.e.\ if and only if $w$ is a quasi-isomorphism. Conversely, if $w$ is a quasi-isomorphism, then $(s,w)$, or $f^{-1} = s/w$, is the inverse in $\DdR$.
    \qed
\end{proof}
\begin{remark}
    In particular, a chain map $g\colon A\to B$ is an isomorphism $g/\id$ in $\DdR$ if and only if $g$ is a quasi-isomorphism. It is important to point out that an isomorphism $A\sim B$ in the derived category $\DdR$ is \emph{not} necessarily represented by a chain map $A\to B$. For example, for $\bR=\Z$, $A=(\Z/2,0)$ and $B=(\Z^2,\dff_B)$ with $\dff_B(a,b)=(2b,0)$, the map $B\to A$, $(a,b)\mapsto a \bmod 2$, is a quasi-isomorphism, while every chain map $A\to B$ is zero. If $A$ is $K$-projective, every morphism $A\to B$ in $\DdR$ is represented by a chain map, cf.\ Defn.~\ref{Kproj} and Prop.~\ref{prop:kprojequiv}(iii).
\end{remark}

\begin{remark}
    \label{Qshapedanalogon}
    The theory of differential modules and associated derived categories may be regarded as a particular case of the $\shQ$-shaped derived categories of Holm and Jørgensen, cf.\ \cite{HolmJorg1}. For a small $\Z$-preadditive category $\shQ$, the category of \emph{$\shQ$-shaped $R$-modules}, denoted $\sRmodRQ$, consists of the $\Z$-linear functors from $\shQ$ to $\sRmodR$. For the quiver $\rmGamma^{\rm diff}$ with one vertex and one loop $\dff$, the path category is given by the polynomial ring $\Z[\dff]$, and the relation $\dff^2$ defines the quotient category $\shQ^{\rm diff} = \Z[\dff]/\langle\dff^2\rangle$, the algebra of dual numbers over $\Z$. It satisfies the standing assumptions of \cite{HolmJorg1} over $\Z$: its Hom-module $\Z\oplus\Z\dff$ is free of rank two, the identity is a Serre functor, and the radical $\Z\dff$ has square zero. The associated category of $\shQ^{\rm diff}$-shaped $R$-modules is $\sRmodRQd$ and, following \cite{HolmJorg1},
\[
\sRmodRQd \cong \sRmod_{R[\epsilon]/\langle\epsilon^2\rangle} \cong \CdR,\quad\text{and}\quad
\DQ = \sRmodRQd[W^{-1}] \cong \DdR,
\]
where $R[\epsilon]/\langle\epsilon^2\rangle$ is the ring of dual numbers and the weak equivalences $W$ in $\sRmodRQd$ are the quasi-isomorphisms; inverting them in $\CdR$ or in $\KdR$ yields the same category, since homotopic chain maps are equal when inverting quasi-isomorphisms.
In \cite{HolmJorg1} various properties of associated derived categories are given which are applicable for the category of differential modules. For example the derived category $\DQ\cong \DdR$ is a triangulated category, cf.\ Rmk.\ \ref{Qshapedanalogon-22}.
\end{remark}

\section{Triangulated structures}
\label{triang}
The homotopy category $\KdR$ and the derived category $\DdR$ are \emph{triangulated categories}, satisfying the standard axioms (TR1) -- (TR4) of a triangulated category, cf.\ {\cite[Defn.~10.2.1]{weibel}}.
We now explain the construction of distinguished triangles in both categories.

\subsection{The mapping cone and distinguished triangles}
\label{mappingconesect}
In the category  $\CdR$ we define the \emph{shift functor} $\sfS\colon\CdR\to \CdR$ via:
\[
\sfS(A,\dff_A) = (A,-\dff_A),\quad \text{and}\quad\sfS f=f,~~\forall f\in \hom(A,B),
\]
{so that $\sfS^2=\id$.} An important construction is
  the \emph{mapping cone} which is defined as follows. Let $f\colon B\to A$
be a chain map, then
\begin{equation}
    \label{dagmapcone}
\cone(f) = \Bigl( A\oplus B,\dff_{f}\Bigr),\quad
\dff_{f} = \begin{pmatrix}
    \dff_A~& f~\\0~& -\dff_B~
\end{pmatrix},
\end{equation}
which yields the short exact sequence
\begin{equation}
    \label{SESforcone}
\begin{tikzcd}
0 \arrow[r] & A \arrow[r, "i", tail] & \cone(f) \arrow[r, "j", two heads] & \sfS B \arrow[r] & 0,
\end{tikzcd}
\end{equation}
where $i$ and $j$ are the standard inclusion and projection maps respectively.
{Every short exact sequence of differential modules yields an exact triangle in homology, with the connecting homomorphism as third map, cf.\ \cite[Ch.~IV, \S1]{CE}. For \eqref{SESforcone}, using $H(\sfS B)=H(B)$, the} long exact sequence on homology becomes:
\begin{equation}
    \label{LESinHom}
\begin{tikzcd}
{} \arrow[r] & H(A) \arrow[r, "H(i)"] & H\bigl(\cone(f)\bigr) \arrow[r, "H(j)"] & H(B) \arrow[r, "H(f)"] & H(A) \arrow[r] & {}
\end{tikzcd}
\end{equation}
where $H(f)$ is the connecting homomorphism and which can   also be represented as {an} exact triangle:
\begin{equation}
    \label{SESforcone2}
\begin{tikzcd}[row  sep=large]
H(A) \arrow[rr, "H(i)"] &                         & H\bigl(\cone(f) \bigr) \arrow[ld, "H(j)"] \\
                        & H(B) \arrow[lu, "H(f)"] &
\end{tikzcd}
\end{equation}
\begin{lemma}
    \label{qisviacone}
    The chain map $f$ is a quasi-isomorphism if and only if $H\bigl(\cone(f)\bigr)=0$.
\end{lemma}
\begin{proof}
If $H\bigl(\cone(f)\bigr)=0$, the exactness of $0\to H(B)\xrightarrow[]{H(f)} H(A)\to 0$ implies that $H(f)$ is an isomorphism and thus $f$ is a quasi-isomorphism. Conversely, if $f$ is a quasi-isomorphism, then $H(f)$ is an isomorphism and therefore {the} exact triangle in \eqref{SESforcone2} yields  $\image H(j) = \ker H(f)=0$, $\ker H(i) = \image H(f) = H(A)$ and therefore $\image H(i) =0$. Finally, $\ker H(j) = \image H(i) =0$ and therefore $H\bigl( \cone(f)\bigr)=0$.
    \qed
\end{proof}
A \emph{triangle} $T$ in $\KdR$ (or $\DdR)$ is a diagram
    \begin{equation}
    \label{disttrian12}
\begin{tikzcd}
X \arrow[r,"r"] & Y \arrow[r,"s"] & Z \arrow[r,"t"] & \sfS X
\end{tikzcd}
\end{equation}
consisting of objects and morphisms in $\KdR$ (or $\DdR$) and is denoted as a {sextuple} $T=(X,Y,Z,r,s,t)$. {An \emph{isomorphism} of triangles $T\to T'$, with $T'=(X',Y',Z',r',s',t')$, is a triple of isomorphisms $u\colon X\to X'$, $v\colon Y\to Y'$ and $w\colon Z\to Z'$ with $vr=r'u$, $ws=s'v$ and $(\sfS u)t=t'w$.}
The sequence of chain maps
\begin{equation}
    \label{shexseq2}
\begin{tikzcd}
B \arrow[r, "f"] & A \arrow[r, "i"] & \cone(f) \arrow[r, "j"] & \sfS B
\end{tikzcd}
\end{equation}
is called a \emph{standard triangle}, which is a triangle in $\CdR$, $\KdR$ and $\DdR$, and
is denoted by  $\bigl(B,A,\cone(f),f,i,j \bigr)$.

\begin{definition}
    \label{defndisttri}
  A triangle $T$   is a   \emph{distinguished triangle} in $\KdR$ if it is isomorphic in $\KdR$ to a standard triangle. A triangle $T$
is a {distinguished triangle} in $\DdR$ if it is isomorphic in $\DdR$  to a standard triangle.
\end{definition}
\begin{theorem}
    \label{KDdRtriangle}
    With the above choices of distinguished triangles the categories $\KdR$ and $\DdR$ are triangulated categories.
\end{theorem}
\begin{proof}
{For $\KdR$ the proof of \cite[Prop.~10.2.4]{weibel} for chain complexes carries over verbatim, since it only uses mapping cones and chain maps between them. Homology $H\colon\KdR\to\sRmodR$ takes distinguished triangles to exact triangles, cf.\ \eqref{SESforcone2}, and $\QS$ is the class of morphisms $s$ for which $H(s)$ is an isomorphism. By \cite[Prop.~10.4.1]{weibel}, $\QS$ is a multiplicative system in $\KdR$, and the derived category $\DdR=\KdR[\QS^{-1}]$ is a triangulated category whose distinguished triangles are the triangles isomorphic to images of distinguished triangles in $\KdR$, i.e.\ to standard triangles, cf.\  Rmk.~\ref{Qshapedanalogon}.}
    \qed
\end{proof}
In the category $\KdR$ short exact sequences
do not necessarily yield {distinguished} triangles in $\KdR$, cf.\ \cite[Rem.\ 10.1.7]{weibel}.
We now explain how this is possible in the derived category $\DdR$.
For a short exact sequence
$0\to X\xrightarrow[]{r} Y \xrightarrow[]{s} Z \to 0$, $r$ is injective, $s$ is surjective and $\image r=\ker s$. This yields the isomorphism
$r^\dagger\colon X\xrightarrow[]{\cong} \image r=\ker s$, given by $x\mapsto r(x)\in \image r\subset Y$, and the short exact sequence
\begin{equation}
    \label{shexseq3}
    \begin{tikzcd}
0 \arrow[r] & \cone(r^\dagger) \arrow[r, "\xi", "\subset"'] \arrow[r] & \cone(r) \arrow[r, "\rmTheta"] & Z \arrow[r] & 0
\end{tikzcd}
\end{equation}
where $\rmTheta$ is defined by $\rmTheta(y,x) = s(y)$
and {$\image \xi = \image r\oplus X = \ker s\oplus X = \ker \rmTheta$}.
The map $r^\dagger\colon (X,\dff_X)\to (\image r,\dff_Y)$ is a chain isomorphism, {hence a quasi-isomorphism, and thus
$H\bigl(\cone(r^\dagger)\bigr)=0$ by Lemma \ref{qisviacone}.}
If we apply the homology functor to the short exact sequence in \eqref{shexseq3} we obtain:
\[
\begin{tikzcd}
0 \arrow[r] & H\bigl(\cone(r)\bigr) \arrow[r, "H(\rmTheta)"] & H(Z) \arrow[r] & 0
\end{tikzcd},
\]
which shows that $\rmTheta\colon \cone(r)\to Z$ is a quasi-isomorphism.
In the derived category $\DdR$ we can invert the quasi-isomorphisms and consider the diagram
\begin{equation}
    \label{constructwithr}
\begin{tikzcd}[column sep=large, row  sep=large]
X \arrow[r, "r"] \arrow[equal]{d} & Y \arrow[r, "s"] \arrow[equal]{d} & Z \arrow[r, "t"] \arrow[d, "\rmTheta^{-1}"] & \sfS X \arrow[equal]{d} \\
X \arrow[r, "r"]                    & Y \arrow[r, "i"]                    & \cone(r) \arrow[r, "j"]    & \sfS X
\end{tikzcd}
\end{equation}
{Its middle square commutes since $\rmTheta i=s$.} The latter proves that the short exact sequence $0\to X\xrightarrow[]{r} Y \xrightarrow[]{s} Z \to 0$ yields a distinguished triangle
$(X,Y,Z,r,s,j\rmTheta^{-1})$ in $\DdR$  given by \eqref{disttrian12} with $t=j\rmTheta^{-1}$, cf.\ \cite[Ex.~10.4.9]{weibel}, which we call the triangle \emph{induced} by the short exact sequence. On homology $H(t)=-\delta$, where $\delta\colon H(Z)\to H(X)$ is the connecting homomorphism of the short exact sequence; for \eqref{SESforcone} this gives $H(t)=-H(f)$, cf.\ \eqref{firstrot12}.
\begin{remark}
    Let us emphasize   that $\rmTheta$ is not necessarily a {chain-equivalence} and therefore
$\rmTheta$ may not be invertible in $\KdR$. In fact, short exact sequences do not necessarily yield distinguished triangles in $\KdR$, {for any choice of $t$: for $\bR=\Z$ the short exact sequence $0\to(\Z/2,0)\xrightarrow[]{2}(\Z/4,0)\to(\Z/2,0)\to0$ does not, cf.\ \cite[Rem.~10.1.7, Exerc.~10.1.2]{weibel};} but {they} do in $\DdR$.
\end{remark}
\begin{remark}
    \label{sometimesinKdR}
    If the short exact sequence
$0\to X\xrightarrow[]{r} Y \xrightarrow[]{s} Z \to 0$ is split { as a sequence of $\bR$-modules}, then $\rmTheta$ is a {chain-equivalence} and  is invertible in $\KdR$, with $t=j\rmTheta^{-1}$. {Indeed, an $\bR$-linear splitting identifies $Y$ with $X\oplus Z$ and $\dff_Y(x,z)=(\dff_Xx+\tau z,\dff_Zz)$, where $\tau\colon\sfS Z\to X$ is a chain map since $\dff_Y^2=0$. Hence $Y=\cone(\tau)$, the short exact sequence is \eqref{SESforcone} for $f=\tau$, and the statement is Proposition~\ref{firstshift}.}
\end{remark}

Distinguished triangles always come from short exact sequences in both categories as we show now. Let $T=(X,Y,Z,r,s,t)$ be  a distinguished triangle in   $\KdR$  displayed in \eqref{disttrian12}. We want to show that there is a short exact sequence generating it.
Consider a chain map $f\colon B\to A$ with associated short exact sequence given by \eqref{SESforcone}. The triangle
\begin{equation}
    \label{firstrot12}
    \begin{tikzcd}
 A \arrow[r, "i"] & \cone(f) \arrow[r, "j"] & \sfS B  \arrow[r, "-\sfS f"]  & \sfS A,
\end{tikzcd}
\end{equation}
 is distinguished, since it is the rotation via (TR2) of the standard triangle, cf.\ \eqref{shexseq2}. Stating it like this is using the fact that   $\KdR$ is triangulated.\footnote{The same applies to distinguished triangles in $\DdR$.} The form in \eqref{firstrot12} can also be obtained via a direct argument. By Proposition~\ref{firstshift}, \eqref{firstrot12} is the triangle induced by the short exact sequence \eqref{SESforcone}, which is split as a sequence of $\bR$-modules. By the rotation axiom (TR2), and since $\sfS^{-1}=\sfS$, the triangle $T'=(\sfS Z,X,Y,-\sfS t,r,s)$ is distinguished and therefore isomorphic in $\KdR$ to a standard triangle:
\[
\begin{tikzcd}[column sep=large, row  sep=large]
\sfS Z \arrow[r, "-\sfS t"] \arrow[d] & X \arrow[r, "r"] \arrow[d] & Y \arrow[r, "s"] \arrow[d] & Z \arrow[d] \\
B \arrow[r, "f"]                      & A \arrow[r, "i"]           & \cone(f) \arrow[r, "j"]    & \sfS B
\end{tikzcd}
\]
Rotating forward shows that $T$ is isomorphic to {\eqref{firstrot12}, i.e.\ to the distinguished triangle induced by the split short exact sequence \eqref{SESforcone}}. This argument works equally well in the derived category $\DdR$. We obtain the following equivalences for the derived category:
\begin{theorem}\label{cor:reduction}
For a triangle $T=(X,Y,Z,r,s,t)$ in $\DdR$ the following statements are equivalent:
\begin{enumerate}[label=\textup{(\roman*)},leftmargin=2.2em,itemsep=0.1em,topsep=0pt]
  \item[{\rm (i)}]$T$ is isomorphic to the triangle induced by a short exact
        sequence of differential modules;
  \item[{\rm (ii)}] $T$ is isomorphic to the triangle induced by a short exact
        sequence of differential modules { which is split as a sequence of $\bR$-modules};
  \item[{\rm (iii)}] $T$ is a distinguished triangle.
\end{enumerate}
\end{theorem}
\begin{proof}
The first argument above shows that (i) implies (iii) and by the second argument we conclude that (iii) implies (ii). Finally, (ii) implies (i) is obvious.
    \qed
\end{proof}
\begin{remark}
    In the homotopy category $\KdR$ distinguished triangles come from short exact {sequences}, split as sequences of $\bR$-modules, but not necessarily vice versa, cf.\ the remark after \eqref{constructwithr}. The latter is true in the derived category $\DdR$ since quasi-isomorphisms are invertible.
\end{remark}

\subsection{$K$-projective differential modules and $K$-projective resolutions}
The notion of $K$-projectivity plays a crucial role for understanding morphisms (roofs) in the derived category and is a concept that is also used in the theory of chain complexes, cf.\ \cite{Spalt}.
\begin{definition}
    \label{Kproj}
    A differential module $(P,\dff_P)$ is \emph{$K$-projective} if for every acyclic differential module {$(N,\dff_N)$}, i.e. {$H(N)=0$}, it holds that every chain map $p\colon P\to N$ satisfies $p\pmb{\sim} 0$, i.e. $p$ represents the zero homomorphism in $\KdR$.% (null-homotopic).
\end{definition}

The notion of $K$-projectivity behaves well with respect to homomorphisms in the derived category. For this we use properties of hom functors  in triangulated categories, cf.\ \cite{Hanratty}{ and \cite[Ex.~10.2.8]{weibel}}.

\begin{proposition}\label{prop:kprojequiv}
For a differential module $(P,\dff_P)$
the following statements are equivalent:
\begin{enumerate}[label=\textup{(\roman*)},leftmargin=2.2em,itemsep=0.1em,topsep=0pt]
\item [{\rm (i)}] $P$ is $K$-projective;
\item [{\rm (ii)}] $\Hom_{\KdR}(P,-)$ carries quasi-isomorphisms to bijections;
\item [{\rm (iii)}] $\Hom_{\KdR}(P,A)\to\Hom_{\DdR}(P,A)$ is bijective for every $A=(A,\dff_A)$.
\end{enumerate}
\end{proposition}
\begin{proof}
(i) implies (ii):
Let $s\colon W\to A$ be a quasi-isomorphism. Consider the associated distinguished triangle
\[
\begin{tikzcd}
W \arrow[r, "s"] & A \arrow[r, "i"] & \cone(s) \arrow[r, "j"] & \sfS W
\end{tikzcd}
\]
By Lemma \ref{qisviacone}, the cone $\cone(s)$ is  acyclic. If we apply the functor $\Hom_{\KdR}(P,-)$, we obtain the long exact sequence:
\[
\begin{aligned}
  \cdots\longrightarrow&\Hom_{\KdR}\big(P,\sfS\cone(s)\big)\longrightarrow\Hom_{\KdR}(P,W)
  \xrightarrow{\ s_*\ }\Hom_{\KdR}(P,A)\longrightarrow\\
  \longrightarrow&\Hom_{\KdR}\big(P,\cone(s)\big)\longrightarrow\cdots,
  \end{aligned}
\]
where $s_*$ {is} the induced map under $\Hom_{\KdR}(P,-)$.
Since, by assumption, $P$ is $K$-projective, we have $\Hom_{\KdR}\big(P,\cone(s)\big)=\Hom_{\KdR}\big(P,\sfS\cone(s)\big)=0$, and thus the long exact sequence reduces to
\[
0\longrightarrow\Hom_{\KdR}(P,W)
  \xrightarrow{\ s_*\ }\Hom_{\KdR}(P,A)\longrightarrow 0,
\]
which {proves} that $s_*$ is bijective.

(ii) implies (i): For $N$ acyclic, the map $0\to N$ is a quasi-isomorphism and $\Hom_{\KdR}(P,0)=0$. Since $s_*$ is a bijection, it follows that $\Hom_{\KdR}(P,N)=0$.

(ii) implies (iii): Fix $A$. We show the localization map
$\sfL\colon\Hom_{\KdR}(P,A)\to\Hom_{\DdR}(P,A)$, defined by $[f]\mapsto
f/\id_P$, is bijective.
Let $w/s$ represent a morphism in $\Hom_{\DdR}(P,A)$, with $w\colon W\to A$ a chain map  and $s\colon W\to P$ a quasi-isomorphism. Then, by assumption,
\[
  s_*\colon\ \Hom_{\KdR}(P,W)\longrightarrow\Hom_{\KdR}(P,P),
  \qquad
  s_*\big([g]\big)=[s g],
\]
is invertible.
By choosing
$[\id_P]\in\Hom_{\KdR}(P,P)$ we obtain a homotopy
class $[g]$ with $s_*[g]=[\id_P]$ and therefore $[sg]=[\id_P]$, i.e.
$sg\pmb{\sim} \id_P$ with a chain map $g\colon P\to W$, which is a quasi-isomorphism since $H(s)H(g)=\id$.
The roof  $w/s$ is equivalent to {$f/\id_P$, with} the chain map $f=wg\colon P\to A$, via the following common refinement:
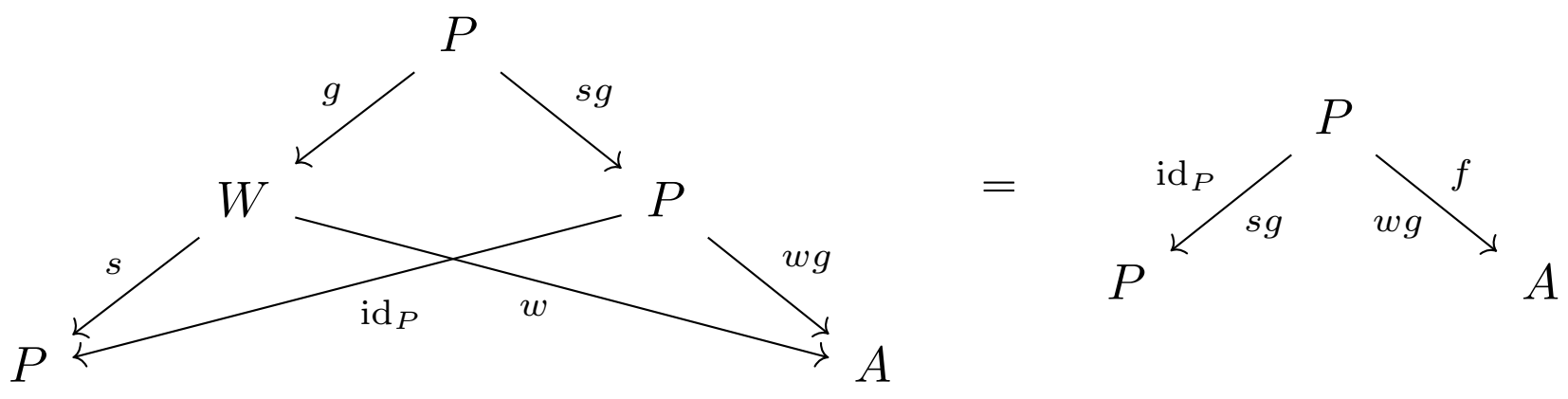
\[
\begin{tikzcd}
  &                                      & P \arrow[ld, "g"'] \arrow[rd, "sg"] &                                         &   \\
  & W \arrow[ld, "s"'] \arrow[rrrd, "w"'] &                                          & P \arrow[rd, "wg"] \arrow[llld, "\id_P"] &   \\
P &                                      &                                          &                                         & A
\end{tikzcd}
\quad = \quad
\begin{tikzcd}
& P \arrow[dl, "sg","\id_P"'] \arrow[dr, "wg"',"f"] & \\
P & & A
\end{tikzcd}
\]
which proves surjectivity.

Suppose {$f\colon P\to A$} is a chain map with
$\sfL[f]=0$. In order to prove injectivity we show that $[f]=0$. Two roofs {$f/\id_P$} and {$0/\id_P$} are equivalent if and only if there is a
quasi-isomorphism $a\colon W\to P$ with $f a \pmb{\sim} 0$, i.e.
\[
\begin{tikzcd}
  &                                      & W \arrow[ld, "a"'] \arrow[rd, "a"] &                                         &   \\
  & P \arrow[ld, "\id_P"'] \arrow[rrrd, "0"'] &                                          & P \arrow[rd, "f"] \arrow[llld, "\id_P"] &   \\
P &                                      &                                          &                                         & A
\end{tikzcd}
\quad = \quad
\begin{tikzcd}
& P \arrow[dl, "a"'] \arrow[dr, "0"',"f a"] & \\
P & & A
\end{tikzcd}
\]
If we apply (ii) as in the previous step, we obtain
 a chain map
$b\colon P\to W$ whose defining property, as a preimage of $[\id_P]$ under
$a_*$, is $a b\pmb{\sim}\id_P$. We conclude that
$f\pmb{\sim} f (ab) = (fa) b\pmb{\sim} 0$, and thus
$[f]=0$ in {$\Hom_{\KdR}(P,A)$}.

(iii) implies (i): An acyclic differential module $N$ is a zero object in $\DdR$ (not necessarily in $\KdR$). The map $i\colon 0\to N$ represents an {isomorphism} in $\DdR$. For any $f\colon P\to N$ in $\DdR$ we have that $[i^{-1} f]=0$ in $\Hom_{\DdR}(P,0)$. Therefore $[f] = [i (i^{-1} f)] = [i][i^{-1} f] = [i]0=0$ in $\DdR$ and thus $\Hom_{\DdR}(P,N)=0$. By (iii) the localization map is bijective and therefore $\Hom_{\KdR}(P,N)=0$.
    \qed
\end{proof}

A consequence of Proposition \ref{prop:kprojequiv} is that in the derived category mappings  out of a $K$-projective module may be represented by  genuine chain maps.
\begin{corollary}
    \label{genchainmaps}
    Let $(P,\dff_P)$ be $K$-projective. Then, a homomorphism in $\Hom_{\DdR}(P,A)$ is represented by a chain map $f\colon P\to A$ for every differential module $(A,\dff_A)$.
\end{corollary}
\begin{proof}
The surjectivity part in the proof of Proposition \ref{prop:kprojequiv} shows that a morphism in $\Hom_{\DdR}(P,A)$,
represented by $w/s$, can be represented by a genuine chain map $f=wg$.
    \qed
\end{proof}
\begin{proposition}\label{prop:kproj}
The ${K}$-projective differential modules are closed under $\sfS^{\pm1}$, finite direct sums,
direct summands, and mapping cones; i.e.\ they form a thick triangulated
subcategory.
\end{proposition}

The properties in Proposition \ref{prop:kproj} are {proved} in Appendix \ref{rottrian}, cf.\ Prop.~\ref{Kprojcones}. {Projectivity of the module $P$ and $K$-projectivity of $(P,\dff_P)$ are independent properties: for $\bR=\Z/4$, the differential module $(\Z/4,\dff)$ with $\dff a=2a$ is free and acyclic, but $\id$ is not null-homotopic, so it is not $K$-projective; for $\bR=\Z$, $\cone(\id_{\mathbb{Q}})$ is contractible, hence $K$-projective, while $\mathbb{Q}$ is not a projective $\Z$-module.}
\begin{definition}
    \label{cofibrantrepl}
    A \emph{$K$-projective resolution} for a differential module $(A,\dff_A)$ is a quasi-isomorphism $\pi\colon (P,\dff_P)\to (A,\dff_A)$, for which $P$ is a projective module and $(P,\dff_P)$ is $K$-projective.
\end{definition}

In \cite{HolmJorg1} a {surjective} $K$-projective resolution is referred to as a \emph{cofibrant replacement} of $A$.
The existence of a $K$-projective resolution is an immediate consequence of the theory of {$\shQ$-shaped categories, cf.\ Rmk.~\ref{Qshapedanalogon-22}}. In Appendix \ref{resolfordiffmod} we give a self-contained proof for completeness.

\begin{proposition}
    \label{existcofrepl}
    Every differential module  admits a $K$-projective resolution, and the $K$-projective resolutions are unique up to homotopy: if $\pi\colon P\to A$ and $\pi'\colon P'\to A$ are $K$-projective resolutions, then there is a chain-equivalence $g\colon P\to P'$, unique up to homotopy, with $\pi' g\pmb{\sim}\pi$.
\end{proposition}
\begin{proof}
By Proposition \ref{prop:semifreeexist} every differential module allows a projective flag resolution, cf.\ App.\ \ref{resolfordiffmod}, which provides a $K$-projective resolution in particular, cf.\ Prop.~\ref{Kprojprop}. By Proposition \ref{prop:kprojequiv}(ii), applied to $P$ and the quasi-isomorphism $\pi'$, there is a unique homotopy class $[g]$ with $\pi'g\pmb{\sim}\pi$, and likewise a unique $[g']$ with $\pi g'\pmb{\sim}\pi'$. Then $\pi(g'g)\pmb{\sim}\pi$, and injectivity of $\pi_*$ on $\Hom_{\KdR}(P,P)$ yields $g'g\pmb{\sim}\id_P$; in the same way $gg'\pmb{\sim}\id_{P'}$.
    \qed
\end{proof}
With $K$-projective resolutions we can prove the following result that allows us to regard morphisms in the derived category $\DdR$ as actual chain maps in the homotopy category.
\begin{proposition}
    \label{morphviaKproj}
    For a $K$-projective resolution $\pi\colon (P,\dff_P)\to (A,\dff_A)$, the map
    \[
    \Hom_{\KdR}(P,B) \xrightarrow[]{\cong} \Hom_{\DdR}(A,B),\quad [f]\mapsto f/\pi,
    \]
    is bijective for every differential module $(B,\dff_B)$. In particular, every morphism $A\to B$ in $\DdR$ may be represented  by a chain map $P\to B$ in $\KdR$, unique up to homotopy of the numerator.
\end{proposition}
\begin{proof}
{The map is the composite of the localization map
\[
\Hom_{\KdR}(P,B)\to\Hom_{\DdR}(P,B),\quad [f]\mapsto f/\id_P,
\]
which is bijective by Proposition \ref{prop:kprojequiv}(iii), and the bijection $\Hom_{\DdR}(P,B)\to\Hom_{\DdR}(A,B)$, $h\mapsto h\pi^{-1}$, since $\pi$ is an isomorphism in $\DdR$ and $(f/\id_P)\pi^{-1}=f/\pi$.}
    \qed
\end{proof}

\begin{remark}\label{Qshapedanalogon-22}
The $\shQ$-shaped framework of Remark~\ref{Qshapedanalogon} has additional properties.
The category $\CdR$ is equipped with two
Quillen-equivalent abelian model structures, cf.\ \cite{Hovey},
whose hypotheses hold
{for $\shQ^{\rm diff}$ over the ground ring $\Z$ and every unital ring $\bR$, \cite[Thm.~6.1, Prop.~6.3]{HolmJorg3}}. The resulting dictionary below, which complements Remark~\ref{Qshapedanalogon},
is proved in {\cite{HolmJorg3,HolmJorg2}}, cf.\ \cite{HolmJorg1} for an
introduction:
\begin{enumerate}[label=\textup{(\roman*)},leftmargin=2.2em,itemsep=0.1em,topsep=0pt]
\item the \emph{exact} objects $\Eexact$
      are exactly the
      \emph{acyclic} differential modules, $(A,\dff_A)\in\Eexact\iff H(A)=0$;
\item the \emph{cofibrant} objects (projective structure) are the
      $K$-projective differential modules with projective
      underlying $R$-module; dually the \emph{fibrant} objects
      (injective structure) are the $K$-injective ones with
      injective underlying module;
\item every $(A,\dff_A)$ admits a \emph{cofibrant replacement}, a surjective
      quasi-isomorphism $\pi\colon P\to A$ with $P$ $K$-projective and
      projective as an $R$-module,
      unique up to homotopy.
\end{enumerate}
\end{remark}

\section{The horseshoe lemma}
\label{hslfordiffmod}
The derived category $\DdR$ allows a direct version of the {horseshoe} lemma, cf.\ \cite[Lem.~2.2.8]{weibel}. We formulate it for sequences of chain maps which are exact up to homotopy.

\begin{definition}
    \label{homotopyexact}
    A sequence of chain maps $X\xrightarrow[]{r}Y\xrightarrow[]{s}Z$ is called
    \emph{homotopy-exact} if the map
\begin{equation}
    \label{thetamap12}
\rmTheta\colon\cone(r)\xrightarrow[]{\sim} Z,\quad \rmTheta(y,x)=s(y),
\end{equation}
is a quasi-isomorphism.
\end{definition}
The fact that $\rmTheta$ is a quasi-isomorphism implies that $\rmTheta$ is a chain map and therefore
\[
\begin{aligned}
\rmTheta\dff_r(y,x)&=s\bigl(\dff_Y y+r(x)\bigr)=s\bigl(\dff_Y y\bigr)+s\big(r(x)\big)\\
&=\dff_Z s(y) +s\big(r(x)\big) =\dff_Z s(y),
\end{aligned}
\]
which implies that $sr=0$, and where $\dff_r$ is the differential on $\cone(r)$, cf.\ \eqref{dagmapcone}.
Homotopy-exact sequences do not require $r$ to be injective, nor $s$ to be surjective.
The \emph{connecting homomorphism} of a homotopy-exact sequence is the morphism $t=j\rmTheta^{-1}\colon Z\to\sfS X$ in $\DdR$, where $i\colon Y\to\cone(r)$ and $j\colon\cone(r)\to\sfS X$ are the standard inclusion and projection. Since $\rmTheta i=s$, Diagram \eqref{constructwithr} commutes in $\DdR$. Hence $(X,Y,Z,r,s,t)$ is isomorphic to the standard triangle of $r$, and is therefore a distinguished triangle in $\DdR$. In particular, it induces an exact triangle in homology. By \eqref{shexseq3} every short exact sequence $0\to X\xrightarrow[]{r}Y\xrightarrow[]{s}Z\to 0$ is homotopy-exact, and $t$ is the {third map of the triangle induced by it, cf.\ \eqref{constructwithr}; on homology $H(t)=-\delta$, where $\delta$ is the connecting homomorphism in homology}.
\begin{remark}
    \label{HEterm}
  The {term} homotopy-exact sequence is motivated by
the construction in  \eqref{shexseq3} and \eqref{constructwithr}, which  replaces $Z$ by $\cone(r)$ and still captures the same homology information and yields a distinguished triangle in $\DdR$.
The word homotopy reflects that $\cone(r)$ is obtained from $Y$ not by dividing out $r(X)$, but by gluing in the contractible module $\cone(\id_X)$, with contracting homotopy $(x,x')\mapsto(0,x)$, along $r$. The cone $\cone(r)$ is the pushout of $Y\xleftarrow{r}X\to\cone(\id_X)$, $x\mapsto(x,0)$,
similar to the topological mapping cone construction, cf.\ \cite[Chap.\ 0]{Hatcher}.
\end{remark}
\begin{lemma}
    \label{natforhomexact}
Let $X\xrightarrow[]{r}Y\xrightarrow[]{s}Z$ and $X'\xrightarrow[]{r'}Y'\xrightarrow[]{s'}Z'$ be homotopy-exact sequences with connecting homomorphisms $t$ and $t'$ respectively, and let $h'\colon X\to X'$, $h\colon Y\to Y'$ and $h''\colon Z\to Z'$ be chain maps with  $hr=r'h'$ and $h''s=s'h$ in $\CdR$. Then, $\sfS h' t = t'h''$ in $\DdR$, which is summarized  by
 the commutative diagram in $\DdR$:
\begin{equation}
     \label{inducedcomm12}
\begin{tikzcd}[column sep=large, row sep=large]
X \arrow[r, "r"] \arrow[d, "h'"'] & Y \arrow[r, "s"] \arrow[d, "h"'] & Z \arrow[r, "t"] \arrow[d, "h''"] & \sfS X \arrow[d, "\sfS h'"] \\
X' \arrow[r, "r'"]                & Y' \arrow[r, "s'"]               & Z' \arrow[r, "t'"]                & \sfS X'
\end{tikzcd}
\end{equation}
of distinguished triangles.
\end{lemma}
\begin{proof}
    The map $h\oplus h'\colon\cone(r)\to\cone(r')$, $(y,x)\mapsto\bigl(h(y),h'(x)\bigr)$, is a chain map since $hr=r'h'$. It satisfies $j(h\oplus h')=\sfS h'\,j$ and $\rmTheta'(h\oplus h')=h''\rmTheta$, because $\rmTheta'\bigl(h(y),h'(x)\bigr)=s'h(y)=h''s(y)$. The latter gives $\rmTheta'^{-1}h''=(h\oplus h')\rmTheta^{-1}$ in $\DdR$, and therefore $t'h'' = j\rmTheta'^{-1}h'' = j(h\oplus h')\rmTheta^{-1}=\sfS h'\, j\rmTheta^{-1}=\sfS h'\, t$.
\qed
\end{proof}

For a chain map $\lambda\colon\sfS Q\to P$, the cone $\cone(\lambda)=P\oplus Q$ has the differential $\dff_\lambda(x,y)=\bigl(\dff_P x+\lambda(y),\dff_Q y\bigr)$, cf.\ \eqref{dagmapcone}. The split short exact sequence $0\to P\xrightarrow[]{i'}\cone(\lambda)\xrightarrow[]{j'}Q\to 0$ of the standard inclusion and projection, cf.\ \eqref{SESforcone}, is homotopy-exact, and its connecting homomorphism is the chain map $-\sfS\lambda$, cf.\ Prop.~\ref{firstshift}. Indeed, the chain map $Q\to\cone(i')$, $y\mapsto\bigl((0,y),-\lambda(y)\bigr)$, is a right inverse of $\rmTheta$ and therefore represents $\rmTheta^{-1}$, and its composition with the standard projection $\cone(i')\to\sfS P$ is $-\lambda$. This is the triangle in \eqref{firstrot12}.

\begin{theorem}[Derived horseshoe lemma]
    \label{SESinDR}
    Let $A\xrightarrow[]{i}B\xrightarrow[]{j}C$ be a homotopy-exact sequence with induced connecting homomorphism denoted by $k\colon C\to\sfS A$. Then:
    \begin{enumerate}[label=\textup{(\roman*)},leftmargin=2.2em,itemsep=0.1em,topsep=0pt]
        \item [{\rm (i)}]
        Let $P$ and $Q$ be differential modules, and let $\pi'\colon P\to A$ and $\lambda\colon\sfS Q\to P$ be chain maps and $\phi\colon Q\to B$   a homomorphism of $R$-modules. Then, $\pi\colon\cone(\lambda)\to B$, given by $\pi(x,y)=i\pi'(x)+\phi(y)$,
        is a chain map if and only if
        \begin{equation}
            \label{chainmapcond}
            \dff_B\phi-\phi\,\dff_Q=i\pi'\lambda,
        \end{equation}
        in which case $j\phi\colon Q\to C$ is a chain map, and the diagram
        \begin{equation}
            \label{distinDR12}
            \begin{tikzcd}[column sep=large, row sep=large]
            P \arrow[d, "\pi'"'] \arrow[r, "i'"] & \cone(\lambda) \arrow[d, "\pi"', dotted] \arrow[r, "j'"] & Q \arrow[dl,"\phi"']\arrow[d, "j\phi",dotted] \arrow[r,"-\sfS \lambda"] &  \sfS P \arrow[d,"\sfS \pi'"]\\
            A \arrow[r, "i"]                     & B \arrow[r, "j"]                 & C \arrow[r,"k"]                    & \sfS A
            \end{tikzcd}
        \end{equation}
        commutes, the first two squares   in $\CdR$ and the third square   in $\DdR$; the diagonal arrow $\phi$ is only a homomorphism of $\bR$-modules. If two of the three chain maps $\pi'$, $\pi$ and $j\phi$ are quasi-isomorphisms, then so is the third. Every chain map $\pi\colon\cone(\lambda)\to B$, for which $\pi i'=i\pi'$ in $\CdR$, is as given above and $\phi(y)=\pi(0,y)$.
        \item [{\rm (ii)}] Let $\pi'\colon P\to A$ and $\pi''\colon Q\to C$ be $K$-projective resolutions. Then, there exist a chain map $\lambda\colon\sfS Q\to P$ and a homomorphism $\phi\colon Q\to B$ which satisfy \eqref{chainmapcond}, such that $j\phi\pmb{\sim}\pi''$  and  the diagram
        \begin{equation}
            \label{distinDR1234}
            \begin{tikzcd}[column sep=large, row sep=large]
            P \arrow[d, "\pi'"'] \arrow[r, "i'"] & \cone(\lambda) \arrow[d, "\pi"',dotted] \arrow[r, "j'"] & Q \arrow[dl,"\phi"',dotted]\arrow[d, "\pi''"] \arrow[r,"-\sfS \lambda",dotted] &  \sfS P \arrow[d,"\sfS \pi'"]\\
            A \arrow[r, "i"]                     & B \arrow[r, "j"]                 & C \arrow[r,"k"]                    & \sfS A
            \end{tikzcd}
        \end{equation}
        commutes, the first square   in $\CdR$ and the  {last} two squares   in $\DdR$.
         If $j$ is surjective, then there exists a map $\phi$   such that $j\phi=\pi''$. The chain map $\pi\colon\cone(\lambda)\to B$ defined in (i) is a $K$-projective resolution.
    \end{enumerate}
\end{theorem}
In (ii) the second square of \eqref{distinDR1234}  commutes in $\DdR$; however,
$j\pi=(j\phi)j'$ is an identity of chain maps, and the identity $j\phi=\pi''$ holds in $\DdR$.

\begin{proof}[of Theorem \ref{SESinDR}]
(i): For $(x,y)\in\cone(\lambda)$ we have
\[
\pi\dff_\lambda(x,y)=i\pi'(\dff_P x)+i\pi'(\lambda(y))+\phi(\dff_Q y),\quad\text{and}\quad \dff_B\pi(x,y)=\dff_B i\pi'(x)+\dff_B\phi(y).
\]
Since $i\pi'$ is a chain map, $\pi$ is a chain map if and only if \eqref{chainmapcond} holds. The property, $ji=0$, combined with \eqref{chainmapcond}, yields
\[
 \dff_Cj\phi= j\dff_B \phi = ji\pi'\lambda + j\phi\dff_Q=j\phi\dff_Q,
\]
which shows that $j\phi$ is a chain map. Similarly, $\pi i'(x)=\pi(x,0)=i\pi'(x)$ and $j\pi(x,y)=ji\pi'(x)+j\phi(y)=(j\phi)\big(j'(x,y)\big)$, using that  $ji=0$. Both rows of \eqref{distinDR12} are homotopy-exact, with connecting homomorphisms $-\sfS\lambda$ and $k$, and the naturality of connecting homomorphisms, with $h'=\pi'$, $h=\pi$ and $h''=j\phi$, shows that the third square commutes in $\DdR$. Both rows induce exact triangles in homology, and the commutative squares induce a morphism between them. If two of the three maps $H(\pi')$, $H(\pi)$ and $H(j\phi)$ are isomorphisms, then so is the third by the five lemma. Finally, if $\pi\colon\cone(\lambda)\to B$ is a chain map with $\pi i'=i\pi'$, then $\pi(x,y)=\pi(x,0)+\pi(0,y)=i\pi'(x)+\phi(y)$ with $\phi(y)=\pi(0,y)$.

(ii): The map $\id\oplus\pi'\colon\cone(i\pi')\to\cone(i)$, $(b,x)\mapsto\bigl(b,\pi'(x)\bigr)$, is a chain map. The chain maps $\pi'$, $\id$, $\id\oplus\pi'$ and $\sfS\pi'$ form a morphism from the standard triangle of $i\pi'$ to the standard triangle of $i$, cf.\ \eqref{shexseq2},
\begin{equation}
    \label{morphofstand12}
\begin{tikzcd}[column sep=large, row sep=large]
P \arrow[r, "i\pi'"] \arrow[d, "\pi'"'] & B \arrow[r, "i^\dagger"] \arrow[d, "\id"'] & \cone(i\pi') \arrow[r, "j^\dagger"] \arrow[d, "\id\oplus\pi'"] & \sfS P \arrow[d, "\sfS \pi'"] \\
A \arrow[r, "i"]                        & B \arrow[r, "i^\ddagger"]                   & \cone(i) \arrow[r, "j^\ddagger"]                                & \sfS A
\end{tikzcd}
\end{equation}
which commutes in $\CdR$, and where $i^\dagger,i^\ddagger,j^\dagger$, and $j^\ddagger$ are the standard inclusion and projection maps. By the long exact sequence \eqref{LESinHom} and the five lemma, $\id\oplus\pi'$ is a quasi-isomorphism. Consequently, the composition $\rmTheta_{\pi'}=\rmTheta(\id\oplus\pi')$, with $\rmTheta$ as defined in \eqref{thetamap12}, given by
\[
\begin{tikzcd}[column sep=large, row sep=large]
\cone(i\pi') \arrow[r, "\id\oplus\pi'"] & \cone(i) \arrow[r, "\rmTheta"] & C,
\end{tikzcd}
\quad (b,x)\mapsto j(b),
\]
is a quasi-isomorphism, because $\rmTheta$ is by assumption.

We now construct a chain map $\psi$ such that the diagram
\[
\begin{tikzcd}[row sep=large]
                                                   & \cone(i\pi') \arrow[d, "\rmTheta_{\pi'}"] \\
Q \arrow[ru, "\psi", dashed] \arrow[r, "\pi''"'] & C
\end{tikzcd}
\]
commutes in $\KdR$, i.e.\ $\rmTheta_{\pi'}\psi{\pmb{\sim}}\pi''$.
By Proposition \ref{prop:kprojequiv}(ii), applied to the $K$-projective {differential} module $Q$ and the quasi-isomorphism $\rmTheta_{\pi'}$, the map
\[
(\rmTheta_{\pi'})_*\colon\Hom_{\KdR}\bigl(Q,\cone(i\pi')\bigr)\to\Hom_{\KdR}(Q,C),\quad [g]\mapsto[\rmTheta_{\pi'}g],
\]
is bijective. Surjectivity yields a chain map $\psi$ with $[\rmTheta_{\pi'}\psi]=[\pi'']$, i.e.\ $\rmTheta_{\pi'}\psi{\pmb{\sim}}\pi''$, and injectivity shows that $\psi$ is unique up to homotopy. We make the construction explicit, since it also provides the homotopy, which is used below.
Consider $\cone(\rmTheta_{\pi'})=C\oplus\cone(i\pi')$, with differential $\dff_{\rmTheta_{\pi'}}$, given by
$
(c,w)\mapsto\bigl(\dff c+\rmTheta_{\pi'}(w),-\dff_{i\pi'}w\bigr),
$
cf.\ \eqref{dagmapcone}, and $\cone(\rmTheta_{\pi'})$ is acyclic by Lemma \ref{qisviacone}. Let $i''\colon C\to\cone(\rmTheta_{\pi'})$, $c\mapsto(c,0)$, be the standard inclusion. Then $i''\pi''\colon Q\to\cone(\rmTheta_{\pi'})$ is a chain map into an acyclic differential module, and since $Q$ is $K$-projective, $i''\pi''\pmb{\sim} 0$, cf.\ Defn.\ \ref{Kproj}. Hence, there exists a homomorphism $h\colon Q\to\cone(\rmTheta_{\pi'})$, such that $i''\pi''=\dff h+h\dff$. Write {$h=(s,\psi)$}, with homomorphisms $s\colon Q\to C$ and $\psi\colon Q\to\cone(i\pi')$. For $y\in Q$ this gives
\[
\begin{aligned}
i''\pi''(y) &= \bigl(\pi''(y),0\bigr)=\dff_{\rmTheta_{\pi'}} h(y)+h(\dff_Q y)\\
&=\Bigl(\dff_C s(y)+s(\dff_Q y)+\rmTheta_{\pi'}\psi(y),\ -\dff_{i\pi'}\psi(y)+\psi(\dff_Q y)\Bigr),
\end{aligned}
\]
and therefore $\psi$ is a chain map, and $s$ satisfies
\[
\pi''-\rmTheta_{\pi'}\psi=\dff_C s+s\dff_Q,
\]
i.e.\ $\rmTheta_{\pi'}\psi\pmb{\sim}\pi''$ with homotopy $s$. Conversely, a chain map $\psi$ together with a homotopy $s$ between $\pi''$ and $\rmTheta_{\pi'}\psi$ defines a null-homotopy {$h=(s,\psi)$} of $i''\pi''$.
Since $\psi$ takes values in $\cone(i\pi')$ we write {$\psi=(\phi,-\lambda)$},
with homomorphisms $\phi\colon Q\to B$ and $\lambda\colon Q\to P$, i.e.\ $-\lambda= j^\dagger\psi$. For $y\in Q$,
\[
\dff_{i\pi'}\psi(y)=\bigl(\dff_B\phi(y)-i\pi'\lambda(y),\ \dff_P\lambda(y)\bigr),\quad\text{and}\quad \psi(\dff_Q y)=\bigl(\phi(\dff_Q y),\ -\lambda(\dff_Q y)\bigr).
\]
The map $\psi$ is a chain map if and only if
\[
\dff_B\phi-\phi\,\dff_Q=i\pi'\lambda\quad\text{and}\quad \dff_P\lambda=-\lambda\,\dff_Q.
\]
The first identity is given by \eqref{chainmapcond}, and the second identity yields the chain map $\lambda\colon\sfS Q\to P$ ($\dff_{\sfS Q} =-\dff_Q$). The sign in $\psi=(\phi,-\lambda)$ is chosen such that $\lambda$, and not $-\lambda$, appears in the differential $\dff_\lambda(x,y)=\bigl(\dff_P x+\lambda(y),\dff_Q y\bigr)$ on $\cone(\lambda)$. Furthermore,
\[
\rmTheta_{\pi'}\psi(y)=\rmTheta\bigl(\phi(y),-\pi'\lambda(y)\bigr)=j\phi(y),
\]
{Therefore,} $j\phi=\rmTheta_{\pi'}\psi$ is a chain map and  $\pi''-j\phi=\dff_C s+s\dff_Q$,
i.e.\ $j\phi\pmb{\sim}\pi''$ with homotopy $s$. {By (i), diagram \eqref{distinDR12} commutes for these $\lambda$ and $\phi$, and since $j\phi\pmb{\sim}\pi''$, so does \eqref{distinDR1234}, the first square in $\CdR$ and the last two squares in $\DdR$.}

Suppose $j$ is surjective. Let $s\colon Q\to C$ be a homotopy with $\pi''-j\phi=\dff_C s+s\dff_Q$. Since $Q$ is a projective module, $s=j\tilde s$ for some homomorphism $\tilde s\colon Q\to B$. Replace $\phi$ by $\phi+\dff_B\tilde s+\tilde s\dff_Q$. Then \eqref{chainmapcond} still holds, since $\dff_B(\dff_B\tilde s+\tilde s\dff_Q)-(\dff_B\tilde s+\tilde s\dff_Q)\dff_Q=0$, and $j\phi=\pi''$.

Finally, $\pi'$ and $j\phi$ are quasi-isomorphisms, and by (i) so is $\pi$. Since $P$ and $Q$ are projective and $K$-projective, so is $\cone(\lambda)$, cf.\ Prop.\ \ref{prop:kproj}, and $\pi$ is a $K$-projective resolution.
    \qed
\end{proof}

\begin{remark}
    \label{comparisonHSL}
    For a short exact sequence $0\to A\xrightarrow[]{i}B\xrightarrow[]{j}C\to 0$ in $\CdR$, the chain map $\psi$ in the proof of Theorem \ref{SESinDR}(ii) can also be constructed in two separate lifting steps, which is the classical construction in the proof of the horseshoe lemma, cf.\ \cite[Lem.\ 2.2.8]{weibel}. Since $Q$ is a projective module and $j$ is surjective, there exists a homomorphism $\phi_0\colon Q\to B$ with $j\phi_0=\pi''$. Then $j(\dff_B\phi_0-\phi_0\dff_Q)=\dff_C\pi''-\pi''\dff_Q=0$, and since $\ker j=\image i$ and $i$ is injective, $\dff_B\phi_0-\phi_0\dff_Q=i\mu$ for a unique homomorphism $\mu\colon Q\to A$. From
\[
\dff_B(\dff_B\phi_0-\phi_0\dff_Q)=-\dff_B\phi_0\dff_Q=-(\dff_B\phi_0-\phi_0\dff_Q)\dff_Q
\]
we obtain
\[
0=\dff_Bi\mu + i\mu \dff_Q = i\big( \dff_A\mu +\mu\dff_Q\big).
\]
The injectivity of $i$ then implies that
$\dff_A\mu=-\mu\dff_Q$, i.e.\ $\mu\colon\sfS Q\to A$ is a chain map. Since $\sfS Q$ is $K$-projective, cf.\ Prop.\ \ref{prop:kproj}, and $\pi'$ is a quasi-isomorphism, Proposition \ref{prop:kprojequiv}(ii) yields a chain map $\lambda\colon\sfS Q\to P$ and a homomorphism $g\colon Q\to A$ with $\pi'\lambda-\mu=\dff_A g-g\dff_Q$. Then, $\phi=\phi_0+ig$ satisfies
\[
\begin{aligned}
\dff_B\phi-\phi\dff_Q &= \dff_B\phi_0-\phi_0\dff_Q + \dff_B ig-ig\dff_Q\\
&=i\mu+i(\dff_A g-g\dff_Q)= i\mu + i(\pi'\lambda-\mu) = i\pi'\lambda,
\end{aligned}
\]
and $j\phi=j\phi_0=\pi''$, using the fact that $ji=0$. We conclude that
 \eqref{chainmapcond} holds, $\psi=(\phi,-\lambda)$ is a chain map, and $\rmTheta(\id\oplus\pi')\psi=\pi''$. In this construction $\phi_0$ is found first and $\lambda$ second, and the first step uses that $j$ is surjective and $\ker j=\image i$. For a general homotopy-exact sequence neither property is available. The proof of Theorem \ref{SESinDR}(ii) replaces the quotient $B/i(A)$ by the cone $\cone(i)$, cf.\ Rmk.\ \ref{HEterm}, and performs both lifting steps at once: $\phi$ and $\lambda$ are the two components of a single lift $\psi$ of $\pi''$ along the quasi-isomorphism $\rmTheta_{\pi'} =\rmTheta(\id\oplus\pi')$.
\end{remark}

A derived version of the horseshoe lemma also holds for distinguished triangles in $\DdR$.

\begin{corollary}[Derived horseshoe lemma for distinguished triangles]
    \label{SESinDRcor}
    Let $T = (A,B,C,i,j,k)$ be a distinguished triangle in $\DdR$ and let $\pi'\colon P\to A$ and $\pi''\colon Q\to C$ be $K$-projective resolutions for $(A,\dff_A)$ and $(C,\dff_C)$ respectively.
    Then, there exist a chain map $\lambda\colon \sfS Q\to P$ and a $K$-projective resolution
    $\pi\colon \cone(\lambda) \to B$ such that the diagram
    \begin{equation}
        \label{distinDR13}
    \begin{tikzcd}[column sep=large, row sep=large]
     P \arrow[d, "\pi'"',"\sim"] \arrow[r, "i'"] & \cone(\lambda) \arrow[d, "\pi"',"\sim",dotted] \arrow[r, "j'"] & Q \arrow[d, "\pi''","\sim"'] \arrow[r,"-\sfS \lambda"] &  \sfS P \arrow[d,"\sfS \pi'","\sim"']\\
     A \arrow[r, "i"]                     & B \arrow[r, "j"]                 & C \arrow[r,"k"]                    & \sfS A
    \end{tikzcd}
    \end{equation}
    commutes in $\DdR$, where the maps $i'$ and $j'$ are the standard inclusion and projection.
\end{corollary}
\begin{proof}
By Theorem \ref{cor:reduction} there exist a short exact sequence $0\to X\xrightarrow[]{r}Y\xrightarrow[]{s}Z\to 0$ and an isomorphism $(\varpi',\varpi,\varpi'')$ in $\DdR$ from the triangle $(X,Y,Z,r,s,t)$ in \eqref{constructwithr} to $(A,B,C,i,j,k)$. Since $P$ and $Q$ are $K$-projective, the isomorphisms $\varpi'^{-1}\pi'$ and $\varpi''^{-1}\pi''$ in $\DdR$ are represented by chain maps $\rho'\colon P\to X$ and $\rho''\colon Q\to Z$, cf.\ Cor.\ \ref{genchainmaps}, which are $K$-projective resolutions, by Lemma~\ref{isoisqis}. The short exact sequence is homotopy-exact with connecting homomorphism $t$, and $s$ is surjective. Theorem \ref{SESinDR}(ii) yields a chain map $\lambda\colon\sfS Q\to P$, a homomorphism $\phi\colon Q\to Y$ with $s\phi=\rho''$, and a $K$-projective resolution $\rho\colon\cone(\lambda)\to Y$ such that \eqref{distinDR12} commutes, with vertical maps $\rho'$, $\rho$, $\rho''$ and $\sfS\rho'$. The morphism $\varpi\rho$ in $\DdR$ is represented by a chain map $\pi\colon\cone(\lambda)\to B$, which is a $K$-projective resolution. In $\DdR$,
\[
\pi i'=\varpi\rho i'=\varpi r\rho'=i\varpi'\rho'=i\pi',\quad j\pi=j\varpi\rho=\varpi''s\rho=\varpi''\rho''j'=\pi''j',
\]
and $k\pi''=k\varpi''\rho''=\sfS\varpi'\,t\rho''=\sfS\varpi'\,\sfS\rho'(-\sfS\lambda)=\sfS\pi'(-\sfS\lambda)$.
    \qed
\end{proof}

\section{Chain systems}
\label{chainsystems}
Consider a  poset $\sQ$ and let $\bm{\sQ_1}$ be the associated posetal category. By $\bm{\sQ_2}$ we denote the arrow category of $\bm{\sQ_1}$, which consists of ordered pairs $(\alpha,\beta)$, $\alpha\le\beta$, with $(\alpha,\beta)\le(\alpha',\beta')$ if $\alpha\le\alpha'$ and $\beta\le\beta'$.
The (functor) categories $\bm{\sQ_n}$ consist of concatenations of $n-1$ arrows and are given by $n$-tuples
$\alpha_1\le\cdots\le\alpha_{n}$ which may be regarded as a path of length {$n-1$} in the directed graph induced by the poset $\sQ$, ordered componentwise. In particular, $\IIthreei$ is the category of triangles in $\IIonei$.
Following \cite{HelleRognes} we consider the covariant functors $\spi_0$, $\spi_1$ and $\spi_2$ acting from $\IIthreei$ to $\IIi$,
given by
\[
(\alpha,\beta,\gamma) \xmapsto{\spi_0} (\alpha,\beta),\quad (\alpha,\beta,\gamma) \xmapsto{\spi_1} (\alpha,\gamma),\text{~~and~~} (\alpha,\beta,\gamma) \xmapsto{\spi_2} (\beta,\gamma),
\]
respectively, and  the natural transformations
$\imath\colon \spi_0\Rightarrow\spi_1$ and
$\jmath\colon \spi_1\Rightarrow\spi_2$
whose components
are given by
\[
(\alpha,\beta) \xmapsto{\imath_{(\alpha,\beta,\gamma)}} (\alpha,\gamma) \text{~~and~~}
(\alpha,\gamma) \xmapsto{\jmath_{(\alpha,\beta,\gamma)}} (\beta,\gamma),
\]
respectively.
\begin{definition}[cf.\ \cite{SpV4}]
    \label{Csys}
    A \emph{chain system}  over a poset $\sQ$ in the category $\DdR$ consists of a covariant functor $\sC\colon \bfQ_2 \to \DdR$
and a natural transformation $\sk\colon \sC \spi_2 \Rightarrow \sfS\sC\spi_0$
such that
\begin{equation}
    \label{funchain1}
\begin{tikzcd}[column sep=large]
\sC \spi_0 \arrow[r, "\sC\imath", Rightarrow] & \sC\spi_1 \arrow[r, "\sC\jmath", Rightarrow] & \sC\spi_2 \arrow[r, "\sk", Rightarrow] & \sfS\sC\spi_0
\end{tikzcd}
\end{equation}
is a distinguished triangle in $\DdR$, where the natural transformations $\sC\imath$ and $\sC\jmath$ are the right whiskerings of $\sC$ and $\imath$, and $\sC$ and $\jmath$ respectively.
A chain system over $\sQ$ is denoted by $\bfC = \bigl(\IIi,\sC,\sk\bigr)$.
\end{definition}

 If we unpack the above definition we set $\sC(\alpha,\beta) = C^\beta_\alpha$,  a differential module with differential $\dff_{C^\beta_\alpha}$ for every $(\alpha,\beta)\in \IIi$, and
\begin{enumerate}[label=\textup{(\roman*)},leftmargin=2.2em,itemsep=0.1em,topsep=0pt]
    \item [(i)]  {morphisms $\ell\colon C^\beta_\alpha \to C^{\beta'}_{\alpha'}$ in $\DdR$ for all $(\alpha,\beta)\le (\alpha',\beta')$, such that $\ell'\ell=\ell''$ for $(\alpha,\beta)\le (\alpha',\beta')\le (\alpha'',\beta'')$, where $\ell'\colon C^{\beta'}_{\alpha'} \to C^{\beta''}_{\alpha''}$ and $\ell''\colon C^\beta_\alpha \to C^{\beta''}_{\alpha''}$, and $\ell=\id$ for $(\alpha,\beta) = (\alpha',\beta')$};
 \item[(ii)] distinguished triangles in $\DdR$
 \begin{equation}
     \label{sheactno1}
\begin{tikzcd}[column sep=large, row sep=large]
 C^\beta_\alpha \arrow[r, "i"] & C^\gamma_\alpha \arrow[r, "j"] & C^\gamma_\beta \arrow[r,"k"] & \sfS C_\alpha^\beta,
\end{tikzcd}
 \end{equation}
for all $(\alpha,\beta,\gamma)\in \IIthreei$, where $i$, $j$ and $k$  {are} morphisms in $\DdR$.
\item[(iii)] morphisms between distinguished triangles:
\[
\begin{tikzcd}[column sep=large, row sep=large]
C^\beta_\alpha \arrow[r, "i"] \arrow[d] & C^\gamma_\alpha \arrow[r, "j"] \arrow[d] & C^\gamma_\beta \arrow[r, "k"] \arrow[d] & \sfS C_\alpha^\beta \arrow[d] \\
C^{\beta'}_{\alpha'} \arrow[r, "i'"]    & C^{\gamma'}_{\alpha'} \arrow[r, "j'"]    & C^{\gamma'}_{\beta'} \arrow[r, "k'"]    & \sfS C_{\alpha'}^{\beta'}
\end{tikzcd}
\]
for all pairs of ordered triples $(\alpha,\beta,\gamma)\le (\alpha',\beta',\gamma')$.\footnote{In the special cases $(\alpha,\beta)\le (\alpha,\gamma)$ and $(\alpha,\gamma)\le (\beta,\gamma)$ the morphisms $\ell$ are denoted by $i$ and $j$ respectively.}
\item[(iv)] {for $\sQ=\sO(\sP)$,} a chain system is \emph{excisive} if $\ell\colon C^\alpha_{\alpha\cap\beta}\xrightarrow[]{\sim} C^{\alpha\cup \beta}_{\beta}$ in $\DdR$ for all $\alpha,\beta\in \sO(\sP)$.
\end{enumerate}
\begin{remark}
    {A chain system is in general not assumed to be excisive. For the homotopy-exact chain systems of Definition~\ref{hexrows} the excision maps $\ell$ are chain maps, and they need not be isomorphisms in $\CdR$, even if they are isomorphisms in $\DdR$. %, cf.\ the end of this section.
    In this paper excision is always considered in $\DdR$.}
\end{remark}
\begin{definition}
    \label{hexrows}
   A \emph{homotopy-exact chain system}  over a poset $\sQ$  consists of a covariant functor $\sC\colon \bfQ_2 \to \CdR$
such that
\begin{equation}
    \label{funchain2}
\begin{tikzcd}[column sep=large]
\sC \spi_0 \arrow[r, "\sC\imath", Rightarrow] & \sC\spi_1 \arrow[r, "\sC\jmath", Rightarrow] & \sC\spi_2
\end{tikzcd}
\end{equation}
is a homotopy-exact sequence in $\CdR$. Via the induced natural transformation
$\sk\colon \sC \spi_2 \Rightarrow \sfS\sC\spi_0$, the triple
 $\bfC=\bigl(\IIi,\sC,\sk\bigr)$ is a chain system  in $\DdR$.

\end{definition}
 If we unpack the above definition we set $\sC(\alpha,\beta) = C^\beta_\alpha$ as before, and obtain {the properties (i)--(iv) listed after Definition~\ref{Csys}, now with chain maps $\ell$, in particular $i$ and $j$, and with $\ell'\ell=\ell''$ in $\CdR$. The rows}
 \begin{equation}
     \label{sheactno12}
\begin{tikzcd}[column sep=large, row sep=large]
 C^\beta_\alpha \arrow[r, "i"] & C^\gamma_\alpha \arrow[r, "j"] & C^\gamma_\beta
\end{tikzcd}
 \end{equation}
{are homotopy-exact sequences in $\CdR$} for all $(\alpha,\beta,\gamma)\in \IIthreei$. The chain map $\rmTheta\colon \cone(i)\to C^\gamma_\beta$, $\rmTheta(y,x) = j(y)$, is a quasi-isomorphism, cf.\ Defn.~\ref{homotopyexact}, and $k=j^\dagger\rmTheta^{-1}$ in $\DdR$ is the connecting homomorphism of the row, where $j^\dagger\colon\cone(i)\to\sfS C^\beta_\alpha$ is the standard projection, cf.\ \eqref{SESforcone}. In the morphisms between distinguished triangles in (iii) the first two squares commute in $\CdR$ and the third square in $\DdR$, cf.\ \eqref{inducedcomm12}. For $\sQ=\sO(\sP)$, a homotopy-exact chain system is excisive if the associated chain system in $\DdR$ is excisive, i.e.\ if the maps $\ell\colon C^\alpha_{\alpha\cap\beta}\to C^{\alpha\cup \beta}_{\beta}$ are isomorphisms in $\DdR$.
\begin{definition}
\label{chainsysmorph}
Let $\bfC=\bigl(\IIi,\sC,\sk\bigr)$ and $\bfC'=\bigl(\IIi,\sC',\sk'\bigr)$ be chain systems over a poset $\sQ$ in the category $\DdR$, with maps $\ell$, $i$, $j$, $k$ and $\ell'$, $i'$, $j'$, $k'$ respectively.
\begin{enumerate}[label=\textup{(\roman*)},leftmargin=2.2em,itemsep=0.1em,topsep=0pt]
\item[(i)] A {morphism of chain systems} $f\colon\bfC\to\bfC'$ is a natural transformation $f\colon\sC\Rightarrow\sC'$ that is compatible with the connecting homomorphisms, i.e.\ $\sk'\circ f\spi_2=\sfS f\spi_0\circ\sk$.
The morphism $f$ is an \emph{isomorphism} of chain systems if all $f^\beta_\alpha$ are isomorphisms in $\DdR$.
\item[(ii)] If $\bfC$ and $\bfC'$ are homotopy-exact chain systems, with functors $\sC,\sC'\colon\IIi\to\CdR$, then morphisms and isomorphisms between $\bfC$ and $\bfC'$ are those defined in (i) between the associated chain systems in $\DdR$, cf.\ Lem.\ \ref{natforhomexact}. A \emph{strict morphism} $f\colon\bfC\to\bfC'$ is a natural transformation $f\colon\sC\Rightarrow\sC'$ of functors $\IIi\to\CdR$, i.e.\ a family of chain maps $f^\beta_\alpha\colon C^\beta_\alpha\to{C'}^\beta_\alpha$ with $\ell'f^\beta_\alpha=f^{\beta'}_{\alpha'}\ell$ in $\CdR$ for all $(\alpha,\beta)\le(\alpha',\beta')$.
\end{enumerate}
In particular, a chain system $\bfC$ in $\DdR$ is {isomorphic to a homotopy-exact chain system} $\bfC'$ if there is an isomorphism as given in (i) between $\bfC$ and the induced chain system in $\DdR$ associated with $\bfC'$.
\end{definition}

If we unpack the above definition in (i), $f$ consists of morphisms $f^\beta_\alpha\colon C^\beta_\alpha\to {C'}^\beta_\alpha$ in $\DdR$, for all $(\alpha,\beta)\in\IIi$, such that $\ell' f^\beta_\alpha=f^{\beta'}_{\alpha'}\ell$ for all $(\alpha,\beta)\le(\alpha',\beta')$, and such that for all $(\alpha,\beta,\gamma)\in\IIthreei$ the diagram
\begin{equation}
\label{chainsysmorphdiag}
\begin{tikzcd}[column sep=large, row sep=large]
C^\beta_\alpha \arrow[r, "i"] \arrow[d, "f^\beta_\alpha"'] & C^\gamma_\alpha \arrow[r, "j"] \arrow[d, "f^\gamma_\alpha"] & C^\gamma_\beta \arrow[r, "k"] \arrow[d, "f^\gamma_\beta"] & \sfS C^\beta_\alpha \arrow[d, "\sfS f^\beta_\alpha"] \\
{C'}^\beta_\alpha \arrow[r, "i'"] & {C'}^\gamma_\alpha \arrow[r, "j'"] & {C'}^\gamma_\beta \arrow[r, "k'"] & \sfS {C'}^\beta_\alpha
\end{tikzcd}
\end{equation}
commutes in $\DdR$. For strict morphisms in (ii) we have a family of chain maps $f^\beta_\alpha\colon C^\beta_\alpha\to{C'}^\beta_\alpha$ with $\ell'f^\beta_\alpha=f^{\beta'}_{\alpha'}\ell$ in $\CdR$ for all $(\alpha,\beta)\le(\alpha',\beta')$, i.e. the first two squares in \eqref{chainsysmorphdiag} commute in $\CdR$, and by Lemma~\ref{chainsysmorphlemma}(ii) the third square {commutes} in $\DdR$.

\begin{lemma}
\label{chainsysmorphlemma}
Let $\bfC$ and $\bfC'$ be chain systems over $\sQ$ in $\DdR$, and let
 {$f$ be a family of morphisms}
$f^\beta_\alpha\colon C^\beta_\alpha\to{C'}^\beta_\alpha$ {in $\DdR$}, $(\alpha,\beta)\in\IIi$.
Then,
\begin{enumerate}[label=\textup{(\roman*)},leftmargin=2.2em,itemsep=0.1em,topsep=0pt]
\item[{\rm (i)}] the maps $f^\beta_\alpha$ commute with all maps $\ell$ in $\DdR$ if and only if they commute in $\DdR$ with all maps $i$ and $j$. If $\bfC$ and $\bfC'$ are homotopy-exact and   $f^\beta_\alpha$ are chain maps, then the same holds in $\CdR$ for the chain maps $\ell$, $i$ and $j$ of the functors $\sC,\sC'\colon\IIi\to\CdR$;
\item[{\rm (ii)}] if $\bfC$ and $\bfC'$ are homotopy-exact, then every strict morphism $f\colon\bfC\to\bfC'$ is a morphism of chain systems, and it is an isomorphism of chain systems if and only if all chain maps $f^\beta_\alpha$ are quasi-isomorphisms;\footnote{In this case the inverse $f^{-1}$ is a morphism of chain systems which in general is not strict. A chain map $g$, regarded as the morphism $g/\id$ in $\DdR$, is an isomorphism if and only if $g$ is a quasi-isomorphism, by Lemma \ref{isoisqis}. Its inverse is the roof $\id/g$, which in general is not represented by a chain map.
}
\item[{\rm (iii)}]
if $\bfC\cong\bfC'$, then $\bfC$ is excisive if and only if $\bfC'$ is excisive.
\end{enumerate}
\end{lemma}
\begin{proof}
(i): For $(\alpha,\beta)\le(\alpha',\beta')$ we have $(\alpha,\beta)\le(\alpha,\beta')\le(\alpha',\beta')$ in $\IIi$. By functoriality, $\ell\colon C^\beta_\alpha\to C^{\beta'}_{\alpha'}$ is the composition of  $i\colon C^\beta_\alpha\to C^{\beta'}_\alpha$   and   $j\colon C^{\beta'}_\alpha\to C^{\beta'}_{\alpha'}$, and likewise for $\bfC'$. Conversely, the maps $i$ and $j$ are maps $\ell$. Due to functoriality of $\sC$ and $\sC'$ the same holds in $\CdR$ for functors $\sC,\sC'\colon\IIi\to\CdR$ and chain maps $f^\beta_\alpha$.

(ii): The first two squares of \eqref{chainsysmorphdiag} commute in $\CdR$, and the third square commutes in $\DdR$ by the naturality of the connecting homomorphisms of the homotopy-exact sequences, Lemma \ref{natforhomexact}.

(iii):
For the excision maps of $\bfC$ and $\bfC'$ we have $f^{\alpha\cup\beta}_\beta\,\ell=\ell'\,f^\alpha_{\alpha\cap\beta}$. If the components of $f$ are isomorphisms, then $\ell$ is an isomorphism if and only if $\ell'$ is.
\qed
\end{proof}

Homotopy-exact chain systems  are specific cases of chain systems in $\DdR$. Many special cases give rise {to} homotopy-exact chain systems and therefore chain systems in $\DdR$.
 For example, short exact sequences:
\begin{equation}
    \label{sheactno2}
\begin{tikzcd}[column sep=large, row sep=large]
0 \arrow[r] & C^\beta_\alpha \arrow[r, "i"] & C^\gamma_\alpha \arrow[r, "j"] & C^\gamma_\beta \arrow[r] & 0,
\end{tikzcd}
\end{equation}
for all $(\alpha,\beta,\gamma)\in \IIthreei$
are homotopy-exact, cf.\ Sect.\ \ref{hslfordiffmod}.
The same applies to the \emph{weakly exact sequences} as introduced in \cite{fran}, i.e. sequences:
\begin{equation}
    \label{sheactno3}
\begin{tikzcd}[column sep=large, row sep=large]
C^\beta_\alpha \arrow[r, "i"] & C^\gamma_\alpha \arrow[r, "j"] & C^\gamma_\beta,
\end{tikzcd}
\end{equation}
that satisfy $i$ injective, $ji=0$, and $\varpi\colon C^\gamma_\alpha/i (C_\alpha^\beta) \xrightarrow[]{\sim} C^\gamma_\beta$, where $\varpi$ is defined as: $\varpi[x] := j(x)$. Indeed, if $[x]=[x']$, then $x-x'=i(y)$. Moreover, $j(x)-j(x')=j(x-x')=j(i(y))=0$, which shows that $\varpi$ is well-defined and $\varpi \pi = j$, with $\pi$ defined by $x\mapsto [x]$. {Since $i$ is injective, the chain map $\cone(i)\to C^\gamma_\alpha/i(C^\beta_\alpha)$, $(y,x)\mapsto[y]$, is a quasi-isomorphism by \eqref{shexseq3}, and $\rmTheta$ is its composition with $\varpi$. Hence} the sequence \eqref{sheactno3} is weakly exact {if and only if} $\rmTheta\colon\cone(i)\to C^\gamma_\beta$ {is} a quasi-isomorphism, i.e.\ if and only if it is homotopy-exact.

A chain system is a more general notion than graded or filtered  differential modules. Let $\sQ=\sO(\sP)$, with $\sP$ a finite poset. A \emph{lattice-filtered module} is a lattice homomorphism $\sO(\sP)\to \Sub A$, into the lattice of submodules of $A$, given by $\alpha \mapsto F_\alpha A$, cf.\ \cite{robbin:salamon2,SpV4}. For a lattice-filtered differential module we add the requirement that
 the differential $\dff_A$ is filtered, i.e.  $\dff_A F_\alpha A\subset F_\alpha A$, for all $\alpha\in \sO(\sP)$. Using quotients we have the short exact sequences:
\begin{equation}
    \label{firstexactseq1233}
\begin{tikzcd}
0 \arrow[r] & F_\alpha A \arrow[r, "i","\subset"'] & F_\beta A \arrow[r, "j"] & \displaystyle{\frac{F_\beta A}{F_\alpha A}}
\arrow[r]  & 0
\end{tikzcd},\quad (\alpha,\beta)\in \bm{\sO(\sP)_2}.
\end{equation}
The same can be written for triples $(\alpha,\beta,\gamma)\in \bm{\sO(\sP)_3}$. If we set $C^\beta_\alpha=F_\beta A/F_\alpha A$ we obtain a chain system $\bfC(A,\dff_A)$ induced by a lattice-filtered differential module $(A,\dff_A)$. The latter is actually a homotopy-exact chain system. Since {$F_{\alpha\cup\beta}A/F_\beta A =(F_\alpha A+F_\beta A)/F_\beta A
\cong F_\alpha A/(F_\alpha A \cap F_\beta A) = F_\alpha A/F_{\alpha\cap \beta}A$}, the chain system $\bfC(A,\dff_A)$ of a lattice-filtered differential module is excisive in $\CdR$. The homomorphisms $\ell\colon C^\beta_\alpha\to C^{\beta'}_{\alpha'}$ are defined by {inclusion followed by projection: $F_\beta A/F_\alpha A \xrightarrow[]{\subset} F_{\beta'} A/F_{\alpha} A \xrightarrow[]{\pi} F_{\beta'} A/F_{\alpha'} A$}.

A special case of a lattice-filtered  module is a \emph{poset-graded module}, or $\sP$-graded module, i.e.\ a module with a direct sum decomposition {into submodules $G_pA$}:
\begin{equation}
    \label{dirsumpgm}
    A=\bigoplus_{p\in \sP} G_pA.
\end{equation}
A poset-graded differential module is a special case of a lattice-filtered module via the formula
$F_\alpha A = \bigoplus_{p\in \alpha} G_pA$. For a poset-graded differential module we add the requirement of a filtered differential, i.e. $\dff_A F_\alpha A\subset F_\alpha A$, for all $\alpha\in \sO(\sP)$. In the case of a poset-graded differential module the differential  $\dff_A$ may be represented as a matrix $\dff_A = \bigl[\dff_A(p,q)\bigr]_{p,q}$ with entries:
\begin{equation}
    \label{diffasmatrix}
    \dff_A(p,q)\colon G_qA \to G_pA,\quad \dff_A(p,q)\neq 0~~\implies~~p\le q.
\end{equation}
A differential is \emph{strict} if the above inequality is strict, i.e.
\begin{equation}
    \label{strictdiff12}
    \dff_A(p,p) =0,\quad\forall~p\in \sP.
\end{equation}
In this case we refer to a \emph{strict} differential module.
For more details on lattice-filtered and poset-graded modules see \cite{robbin:salamon2} and \cite{SpV4}.
A poset-graded module is projective if and only if $G_pA$ is projective for all $p\in \sP$.

We can also consider $\alpha\mapsto F_\alpha A$ as a poset-filtered differential module. In this case we have that $\alpha\subset\beta$ implies $F_\alpha A\subset F_\beta A$ and the differential is filtered: $\dff_A F_\alpha A\subset F_\alpha A$, for all $\alpha\in \sO(\sP)$. In this case we again obtain a homotopy-exact chain system $\bfC(A,\dff_A)$. However, the latter is not excisive in $\CdR$ in general.
    A chain complex for an $\sO(\sP)$-filtered topological space, i.e.\ a lattice homomorphism $\alpha\mapsto X_\alpha$ into the subspaces of $X$, is \emph{not} $\sO(\sP)$-filtered, but is poset-filtered. For singular chains, $C(X_\alpha)+C(X_\beta)\subsetneq C(X_{\alpha\cup\beta})$ in general. The latter does yield an {excisive} chain system in $\DdR$  if the inclusions $C(X_\alpha)+C(X_\beta)\to C(X_{\alpha\cup\beta})$ are quasi-isomorphisms, e.g.\ if the interiors of $X_\alpha$ and $X_\beta$ cover $X_{\alpha\cup\beta}$, cf.\ \cite[Prop.~2.21]{Hatcher}, or for subcomplexes of a CW complex and for index filtrations \cite{fran86}.

\section{Main result}
\label{mainresult}
A $\sP$-graded differential module $A=\bigoplus_{p\in \sP} G_pA$ with an $\sO(\sP)$-filtered differential $\dff_A$ yields a chain system $\bfC(A,\dff_A)$ via the exact sequences
\[
\begin{tikzcd}
0 \arrow[r] & G_{\beta\smin\alpha}A \arrow[r, "i"] & G_{\gamma\smin\alpha}A \arrow[r, "j"] & G_{\gamma\smin\beta}A \arrow[r] & 0
\end{tikzcd}
\]
where $G_{\beta\smin\alpha}A=\bigoplus_{p\in \beta\smin\alpha} G_pA$, for all {$\alpha\subset\beta\subset\gamma$}. An $\sO(\sP)$-filtered differential {induces} a differential on the modules $G_{\beta\smin\alpha}A\cong F_\beta A/F_\alpha A$, and the chain system $\bfC(A,\dff_A)$ is said to be \emph{$\sP$-split}. Its rows are short exact sequences, and thus $\bfC(A,\dff_A)$ is a homotopy-exact chain system.

\begin{theorem}
    \label{mainres987655}
    Let $R$ be an   associative unital ring and $\sP$ a finite poset.
    Every excisive, homotopy-exact chain system $\bfC$
    {over $\sO(\sP)$} is isomorphic to a $\sP$-split chain system $\bfC(A,\dff_A)$, for some projective $\sP$-graded differential module $(A,\dff_A)$.
    The isomorphism is given by $K$-projective resolutions $f^\beta_\alpha\colon G_{\beta\smin\alpha}A\to C^\beta_\alpha$, that commute with the chain maps $i$ and $j$ in $\CdR$, i.e.\ by a strict morphism $\bfC(A,\dff_A)\to\bfC$, cf.\ Defn.\ \ref{chainsysmorph}(ii).
\end{theorem}
We follow the inductive procedure as used in \cite{SpV4}, with the  derived horseshoe lemma, Theorem \ref{SESinDR}, as the main tool. The only choices are the $K$-projective resolutions in Step (ii), one for each $q\in\sP$, and the maps provided by Theorem \ref{SESinDR}(ii). All other maps are given by explicit formulas.

\begin{proof}[of Theorem \ref{mainres987655}]
The induction hypothesis for a down-set $\alpha\in\sO(\sP)$ is:
\begin{enumerate}[label=\textup{(\roman*)},leftmargin=2.2em,itemsep=0.1em,topsep=0pt]
\item[(a)] $F_\alpha A=\bigoplus_{p\in \alpha}G_pA$ is a projective $\bR$-module with an upper-triangular differential $\dff_{F_\alpha A}$, i.e.\ $\dff_{F_\alpha A}(p,p')\neq 0$ implies $p\le p'$, for all $p,p'\in \alpha$, cf.\ \eqref{diffasmatrix}, and the differential modules $(G_pA,\dff_{G_pA})$ are $K$-projective;
\item[(b)] there are $K$-projective resolutions $f^{\alpha'}_\smbt\colon F_{\alpha'}A\to C^{\alpha'}_\smbt$ for all down-sets $\alpha'\subset\alpha$;
\item[(c)] the resolutions commute with the inclusions in $\CdR$:
\begin{equation}
\label{natinCdR}
f^{\beta'}_\smbt i' = i f^{\alpha'}_\smbt,\quad\text{for all}~~\alpha'\subset\beta'\subset\alpha.
\end{equation}
\end{enumerate}
Note that \eqref{natinCdR} is an identity of chain maps and not only an identity in $\DdR$.

{\it Step (i):} Put $F_\smbt A=0$ and $f^\smbt_\smbt=0$. Since $C^\smbt_\smbt=0$, because the row of the triple $(\bt,\bt,\bt)$ is homotopy-exact with $i=j=\id$, so that $\id=ji=0$ on $C^\smbt_\smbt$, the induction hypothesis holds for $\alpha=\bt$.

{\it Step (ii):} Suppose the induction hypothesis holds for $\alpha\neq\sP$.
Let $q\in \sP\smin \alpha$ be  a minimal element, i.e. $\alpha^\dagger = \alpha\cup \{q\}\in \sO(\sP)$.
Define $\beta=\big\downarrow q$ and $\beta^\dagger\prec \beta$, where $\prec$ is the covering relation, and thus $\beta = \beta^\dagger\cup \{q\}$. Since $q$ is minimal in $\sP\smin\alpha$ we have $\beta^\dagger\subset\alpha$, and
$f'=f^{\beta^\dagger}_\smbt\colon F_{\beta^\dagger}A\to C^{\beta^\dagger}_\smbt$ is a $K$-projective resolution by the induction hypothesis. Let $f''\colon G_qA\to C^\beta_{\beta^\dagger}$ be a $K$-projective resolution, cf.\ Prop.\ \ref{existcofrepl}. Theorem \ref{SESinDR}(ii), applied to the homotopy-exact sequence $C^{\beta^\dagger}_\smbt\xrightarrow[]{i}C^\beta_\smbt\xrightarrow[]{j}C^\beta_{\beta^\dagger}$ of the triple $(\bt,\beta^\dagger,\beta)$ and the resolutions $f'$ and $f''$, yields a chain map $\lambda_q\colon\sfS G_qA\to F_{\beta^\dagger}A$ and a homomorphism $\phi_q\colon G_qA\to C^\beta_\smbt$ such that
\begin{equation}
\label{chainmapcondq}
\dff_{C^\beta_\smbt}\phi_q-\phi_q\dff_{G_qA} = if'\lambda_q,\quad\text{and}\quad j\phi_q {\pmb{\sim}} f''.
\end{equation}
Define  $F_\beta A=\cone(\lambda_q)$, i.e.\ $F_\beta A=F_{\beta^\dagger}A\oplus G_qA$ with differential $\dff_{F_\beta A}(x,y) = \dff_{\lambda_q}(x,y) = \bigl(\dff_{F_{\beta^\dagger}A} x+\lambda_q(y),\dff_{G_qA } {y}\bigr)$. By Theorem \ref{SESinDR},
\begin{equation}
\label{newfbeta}
f^\beta_\smbt\colon F_\beta A\to C^\beta_\smbt,\quad f^\beta_\smbt(x,y)=if'(x)+\phi_q(y),
\end{equation}
is a $K$-projective resolution with $f^\beta_\smbt i'=if'$ and $jf^\beta_\smbt=(j\phi_q)j'$ in $\CdR$.\footnote{If the rows of $\bfC$ are short exact sequences, then $\phi_q$ can be chosen such that $j\phi_q=f''$.}  {The construction is summarized in the following diagram:}
\begin{equation}
    \label{firstdiaginprf61}
\begin{tikzcd}[column sep=large, row sep=large]
F_{\beta^\dagger}A \arrow[d, "f'"',"f_\smbt^{\beta^\dagger}"] \arrow[r, "i'"] & F_{\beta}A \arrow[r, "j'"] \arrow[d, "f"',"f_\smbt^\beta", dotted] & G_qA \arrow[r, "-\sfS \lambda_q", dotted] \arrow[d, "f''"] \arrow[ld, "\phi_q"', dotted] & \sfS F_{\beta^\dagger}A \arrow[d, "\sfS f'"] \\
C^{\beta^\dagger}_\smbt \arrow[r, "i"]              & C^{\beta}_\smbt \arrow[r, "j"]                     & C_{\beta^\dagger}^\beta \arrow[r, "k"]                                                   & \sfS C^{\beta^\dagger}_\smbt
\end{tikzcd}
\end{equation}
where the first square commutes in $\CdR$ and the  {last two squares} in $\DdR$.  {For $\beta^\dagger=\bt$, i.e.\ $q$ minimal in $\sP$, we have $\lambda_q=0$ and $F_\beta A=G_qA$.}

Finally, let $i'\colon F_{\beta^\dagger}A\to F_\alpha A$ be the inclusion for the pair $\beta^\dagger\subset\alpha$.\footnote{This is not the inclusion $F_{\beta^\dagger}A\to F_\beta A$ used above, since $F_\beta A \not\subseteq F_\alpha A$.} By (a), $F_{\beta^\dagger}A$ is a differential submodule of $F_\alpha A$, so that $i'$ is a chain map, and $i'\lambda_q\colon\sfS G_qA\to F_\alpha A$ is the chain map $\lambda_q$ regarded as a map into $F_\alpha A$. Define
\[
F_{\alpha^\dagger}A=\cone(i'\lambda_q)=F_\alpha A\oplus G_qA,\quad \dff_{F_{\alpha^\dagger}A}(x,y) = \bigl(\dff_{F_\alpha A} x+i'\lambda_q(y),\dff_{G_qA} y\bigr).
\]
Since $\lambda_q$ takes values in $F_{\beta^\dagger}A=\bigoplus_{p<q}G_pA$, the new column of the differential has entries $\dff_{F_{\alpha^\dagger}A}(p,q)\neq 0$ only for $p\le q$, and $\dff_{F_{\alpha^\dagger}A}(q,q)=\dff_{G_qA}$. Hence $\dff_{F_{\alpha^\dagger}A}$ is upper-triangular and (a) holds for $\alpha^\dagger$. In particular,
$F_\beta A=\cone(\lambda_q)$, as defined above, is a differential submodule of $F_{\alpha^\dagger}A$.

{\it Step (iii):} The next step is to construct the $K$-projective resolutions for all down-sets containing $q$ (not just $\beta=\bigl\downarrow q$). Let $\delta\subset\alpha^\dagger$ be a down-set with $q\in \delta$. Then, $q$ is maximal in $\delta$
and consequently, $\gamma=\delta\smin\{q\}\subset\alpha$ is a down-set. Since $\beta\subset\delta$,
\begin{equation}
\label{irrjoin}
\delta=\beta\cup\gamma,\quad\text{and}\quad \beta\cap\gamma=\beta^\dagger.
\end{equation}
Then, $F_{\beta\cup\gamma}A=F_\gamma A\oplus G_qA=\cone(i'\lambda_q)$, with $i'\colon F_{\beta^\dagger}A\to F_\gamma A$. Define
\begin{equation}
\label{fbetagamma}
f^{\beta\cup\gamma}_\smbt\colon F_{\beta\cup\gamma}A\to C^{\beta\cup\gamma}_\smbt,\quad
f^{\beta\cup\gamma}_\smbt(x,y) = i_\gamma f^\gamma_\smbt(x) + i_\beta\phi_q(y),
\end{equation}
with $i_\gamma\colon C^\gamma_\smbt\to C^{\beta\cup\gamma}_\smbt$ and $i_\beta\colon C^\beta_\smbt\to C^{\beta\cup\gamma}_\smbt$. For $\gamma=\beta^\dagger$ we have $i_\beta=\id$, and this is \eqref{newfbeta}. All maps $f^{\beta\cup\gamma}_\smbt$ use the same chain map $\lambda_q$ and the same homomorphism $\phi_q$.

Apply Theorem \ref{SESinDR}(i) to the homotopy-exact sequence $C^\gamma_\smbt\xrightarrow[]{i_\gamma}C^{\beta\cup\gamma}_\smbt\xrightarrow[]{j_\gamma}C^{\beta\cup\gamma}_\gamma$ of the triple $(\bt,\gamma,\beta\cup\gamma)$, with $\pi'=f^\gamma_\smbt$, with $\lambda=i'\lambda_q$, and with $\phi=i_\beta\phi_q$. Let $i^\beta_{\beta^\dagger}\colon C^{\beta^\dagger}_\smbt\to C^\beta_\smbt$ and $i^\gamma_{\beta^\dagger}\colon C^{\beta^\dagger}_\smbt\to C^\gamma_\smbt$ be the maps $i$ of the associated homotopy-exact sequences of the triples $(\bt,\beta^\dagger,\beta)$ and $(\bt,\beta^\dagger,\gamma)$, so that $i=i^\beta_{\beta^\dagger}$ in \eqref{chainmapcondq}. By \eqref{natinCdR}, $f^\gamma_\smbt i'=i^\gamma_{\beta^\dagger}f' = i^\gamma_{\beta^\dagger} f^{\beta^\dagger}_\smbt$, and by functoriality of $\sC$, $i_\beta i^\beta_{\beta^\dagger}=i_\gamma i^\gamma_{\beta^\dagger}$. Hence, by \eqref{chainmapcondq},
and the fact that $i_\beta$ is a chain map:
\[
i_\beta i^\beta_{\beta^\dagger} f^{\beta^\dagger}_\smbt \lambda_q=i_\beta \dff_{C^\beta_\smbt} \phi_q-i_\beta\phi_q\dff_{G_qA} =
\dff_{C^{\beta\cup \gamma}_\smbt}(i_\beta\phi_q)-(i_\beta\phi_q)\dff_{G_{q}A}.
\]
Using the compositions for the $i$-maps we obtain:
\[
\dff_{C^{\beta\cup \gamma}_\smbt}(i_\beta\phi_q)-(i_\beta\phi_q)\dff_{G_{q}A} = i_\beta i^\beta_{\beta^\dagger} f^{\beta^\dagger}_\smbt \lambda_q = i_\gamma i^\gamma_{\beta^\dagger}f^{\beta^\dagger}_\smbt\lambda_q = i_\gamma f^\gamma_\smbt\,i'\lambda_q,
\]
which verifies \eqref{chainmapcond}. By Theorem \ref{SESinDR}(i), $f^{\beta\cup\gamma}_\smbt$ is a chain map with $f^{\beta\cup\gamma}_\smbt i'_\gamma=i_\gamma f^\gamma_\smbt$ and $j_\gamma f^{\beta\cup\gamma}_\smbt=(j_\gamma i_\beta\phi_q)j'_\gamma$ in $\CdR$.
Summarizing, we obtain the diagram:
\begin{equation}
    \label{seconddiaginprf61}
    \begin{tikzcd}[column sep=large, row sep=large]
F_{\gamma}A \arrow[d, "f^\gamma_\smbt"'] \arrow[r, "i_\gamma'"] & F_{\beta\cup\gamma}A \arrow[r, "j_\gamma'"] \arrow[d, "f","f^{\beta\cup\gamma}_\smbt"', dotted] & G_qA \arrow[d, "j_\gamma i_\beta \phi_q"] \arrow[ld, "i_\beta\phi_q"'] \arrow[r, "-\sfS i'\lambda_q"] & \sfS F_{\gamma}A \arrow[d, "\sfS f^\gamma_\smbt"] \\
C^{\gamma}_\smbt \arrow[r, "i_\gamma"]       & C^{\beta\cup\gamma}_\smbt \arrow[r, "j_\gamma"]               & C^{\beta\cup\gamma}_\gamma \arrow[r, "k_\gamma"]                                                             & \sfS C^{\gamma}_\smbt
\end{tikzcd}
\end{equation}
By functoriality, the homomorphism $j_\gamma i_\beta=\ell j$, where $j\colon C^\beta_\smbt\to C^\beta_{\beta^\dagger}$ is the map $j$ in \eqref{chainmapcondq} and $\ell\colon C^\beta_{\beta^\dagger}\to C^{\beta\cup\gamma}_\gamma$. By \eqref{irrjoin}, $\ell\colon C^\beta_{\beta\cap\gamma}\to C^{\beta\cup\gamma}_\gamma$ is an excision map and therefore a quasi-isomorphism, since $\bfC$ is excisive. Moreover, \eqref{chainmapcondq} implies that $j\phi_q\pmb{\sim} f''$, and thus the chain map $j_\gamma i_\beta\phi_q=\ell j\phi_q$ is a quasi-isomorphism. By Theorem \ref{SESinDR}(i) then also $f^{\beta\cup\gamma}_\smbt$ is a quasi-isomorphism. By construction $F_{\beta\cup\gamma}A=\cone(i'\lambda_q)$ is projective and $K$-projective, cf.\ Prop.\ \ref{prop:kproj}, and therefore $f^{\beta\cup\gamma}_\smbt$ is a $K$-projective resolution.

Finally, we check that \eqref{natinCdR} is satisfied for all
$\alpha'\subset\beta'\subset\alpha^\dagger$. If $q\notin\beta'$, then \eqref{natinCdR} holds by the induction hypothesis. If $q\notin\alpha'$ and $q\in\beta'$, then $\alpha'\subset \gamma''=\beta'\smin\{q\}$, and thus $\alpha'\subset \gamma''\subset\alpha$.
For $x\in F_{\alpha'}A$,
Equation \eqref{fbetagamma}, with $\beta'=\beta\cup\gamma''$, implies that $f^{\beta'}_\smbt(x,0) = i^{\beta'}_{\gamma''}f^{\gamma''}_\smbt(x)$.
By \eqref{natinCdR}, using {$\alpha'\subset\gamma''\subset \alpha$},
\[
f^{\beta'}_\smbt(x,0) = i^{\beta'}_{\gamma''}f^{\gamma''}_\smbt(x) = i^{\beta'}_{\gamma''} i^{\gamma''}_{\alpha'} f^{\alpha'}_\smbt(x) = i^{\beta'}_{\alpha'}f^{\alpha'}_\smbt(x).
\]
If $q\in\alpha'$, define $\gamma'=\alpha'\smin\{q\}$. Then, $\gamma'\subset\gamma''\subset \alpha$. For $x\in F_{\alpha'\smin\{q\}}A$ and $y\in G_qA$, using the functoriality of $\sC$, \eqref{natinCdR} and \eqref{fbetagamma},\footnote{For down-sets $\delta'\subset\delta$ we write $i^\delta_{\delta'}\colon C^{\delta'}_\smbt\to C^\delta_\smbt$ and $j^\delta_{\delta'}\colon C^\delta_\smbt\to C^\delta_{\delta'}$ for the maps $i$ and $j$ of the triple $(\bt,\delta',\delta)$. By functoriality of $\sC$, $i^\delta_{\delta'}i^{\delta'}_{\delta''}=i^\delta_{\delta''}$ for $\delta''\subset\delta'\subset\delta$.}
\[
\begin{aligned}
f^{\beta'}_\smbt(x,y) &= i^{\beta'}_{\gamma''}f^{\gamma''}_\smbt(x)+i^{\beta'}_\beta\phi_q(y) = i_{\gamma''}^{\beta'}i_{\gamma'}^{\gamma''}f^{\gamma'}_\smbt(x)+i_{\alpha'}^{\beta'}i^{\alpha'}_\beta\phi_q(y)\\
&= i_{\alpha'}^{\beta'}i_{\gamma'}^{\alpha'}f^{\gamma'}_\smbt(x)+i_{\alpha'}^{\beta'}i^{\alpha'}_\beta\phi_q(y) = i_{\alpha'}^{\beta'} \Bigl( i_{\gamma'}^{\alpha'}f^{\gamma'}_\smbt(x)+i^{\alpha'}_\beta\phi_q(y)\Bigr) \\
&= i f^{\alpha'}_\smbt(x,y).
\end{aligned}
\]
Hence \eqref{natinCdR} holds for all $\alpha'\subset\beta'\subset\alpha^\dagger$, which completes the induction step, and  (a)--(c) hold for $\alpha^\dagger\subset \sP$.

{\it Step (iv):} By Steps (i)--(iii) and induction over the finite poset $\sP$ we conclude that (a)--(c)  {hold} for $\alpha=\tp$. This yields the projective $\sP$-graded differential module $A=\bigoplus_{p\in\sP}G_pA$, with an $\sO(\sP)$-filtered differential $\dff_A$, and $K$-projective resolutions $f^\alpha_\smbt\colon F_\alpha A\to C^\alpha_\smbt$ which satisfy \eqref{natinCdR} for all $\alpha\subset\beta$ in $\sO(\sP)$.

 We now construct $f^\beta_\alpha$ for all $\alpha\subset\beta$. Since $\dff_A$ is upper-triangular, $F_{\beta}A=F_{\alpha}A\oplus G_{\beta\smin\alpha}A=\cone(\lambda_{\beta\smin\alpha})$, where $\lambda_{\beta\smin\alpha}\colon\sfS G_{\beta\smin\alpha}A\to F_{\alpha}A$
is the chain map
 given by the component of $\dff_A$ from $G_{\beta\smin\alpha}A$ to $F_{\alpha}A$, and $i'$ and $j'$ are the standard inclusion and projection. By \eqref{natinCdR}, $f^{\beta}_\smbt i'=if^{\alpha}_\smbt$, and Theorem \ref{SESinDR}(i), applied to the row of the triple $(\bt,\alpha,\beta)$ with $\pi'=f^{\alpha}_\smbt$ and $\pi=f^{\beta}_\smbt = if^\alpha_\smbt + \phi_{\beta\smin\alpha}$, where $\phi_{\beta\smin\alpha}(y)=f^\beta_\smbt(0,y)$, yields the diagram
\begin{equation}
    \label{thirddiaginprf61}
    \begin{tikzcd}[column sep=large, row sep=large]
F_{\alpha}A \arrow[d, "f^\alpha_\smbt"'] \arrow[r, "i'"] & \cone(\lambda_{\beta\smin\alpha}) \arrow[r, "j'"] \arrow[d, "f^\beta_\smbt"', dotted] & G_{\beta\smin\alpha}A \arrow[d, "j \phi_{\beta\smin\alpha}","f^\beta_\alpha"',dotted] \arrow[ld, "\phi_{\beta\smin\alpha}"'] \arrow[r, "-\sfS \lambda_{\beta\smin\alpha}"] & \sfS F_{\alpha}A \arrow[d, "\sfS f^\alpha_\smbt"] \\
C^{\alpha}_\smbt \arrow[r, "i"]                          & C^{\beta}_\smbt \arrow[r, "j"]                                                        & C^{\beta}_\alpha \arrow[r, "k"]                                                                                                                   & \sfS C^{\alpha}_\smbt
\end{tikzcd}
\end{equation}
where the chain map $f^{\beta}_{\alpha}\colon G_{\beta\smin\alpha}A\to C^{\beta}_{\alpha}$ is defined by
\begin{equation}
\label{deffbetaalpha}
f^{\beta}_{\alpha} = j\phi_{\beta\smin\alpha}\quad \implies \quad f^{\beta}_{\alpha}(y)=jf^{\beta}_\smbt(0,y),\quad\text{i.e.}\quad f^{\beta}_{\alpha}j' = jf^{\beta}_\smbt.
\end{equation}
The chain maps $f^{\alpha}_\smbt$ and $f^{\beta}_\smbt$ are quasi-isomorphisms, which implies that $f^{\beta}_{\alpha}$ {is a} quasi-isomorphism by the five lemma, cf.\ Thm.\ \ref{SESinDR}(i). As a finite iterated cone of the modules $G_pA$, $p\in\beta\smin\alpha$,  $G_{\beta\smin\alpha}A$ is projective and $K$-projective, cf.\ Prop.\ \ref{prop:kproj}, and therefore $f^{\beta}_{\alpha}$ is a $K$-projective resolution. For $\alpha=\bt$ this is $f^{\beta}_\smbt$, and $f^\beta_{\beta^\dagger}=j\phi_q$ for $\beta=\big\downarrow q$.

Let $\alpha\subset\beta\subset\gamma$. Without indices, $i$, $j$, $i'$ and $j'$ denote the maps of the triple $(\alpha,\beta,\gamma)$  in Diagram \eqref{diagforall}.
By \eqref{deffbetaalpha}, \eqref{natinCdR}, functoriality of $\sC$, and the identities ${j'}^\gamma_\alpha{i'}^\gamma_\beta=i'{j'}^\beta_\alpha$ and $j'{j'}^\gamma_\alpha={j'}^\gamma_\beta$,\footnote{For down-sets $\delta'\subset\delta$ we write ${i'}^\delta_{\delta'}\colon F_{\delta'}A\to F_\delta A$ and ${j'}^\delta_{\delta'}\colon F_\delta A\to G_{\delta\smin\delta'}A$ for the maps of the triple $(\bt,\delta',\delta)$, cf.\ Step (iii). } and ${j}^\gamma_\alpha{i}^\gamma_\beta=i{j}^\beta_\alpha$ and $j{j}^\gamma_\alpha={j}^\gamma_\beta$,
\[
\begin{aligned}
f^{\gamma}_{\alpha}i'{j'}^\beta_\alpha &= f^{\gamma}_{\alpha}{j'}^\gamma_\alpha{i'}^\gamma_\beta  \overset{\eqref{deffbetaalpha}}{=} j^\gamma_\alpha f^{\gamma}_\smbt{i'}^\gamma_\beta \overset{\eqref{natinCdR}}{=} j^\gamma_\alpha i^\gamma_\beta f^{\beta}_\smbt
= i\,j^\beta_\alpha f^{\beta}_\smbt \overset{\eqref{deffbetaalpha}}{=} if^{\beta}_{\alpha}{j'}^\beta_\alpha,\\
f^{\gamma}_{\beta}j'{j'}^\gamma_\alpha &= f^{\gamma}_{\beta}{j'}^\gamma_\beta \overset{\eqref{deffbetaalpha}}{=} j^\gamma_\beta f^{\gamma}_\smbt = j\,j^\gamma_\alpha f^{\gamma}_\smbt \overset{\eqref{deffbetaalpha}}{=} jf^{\gamma}_{\alpha}{j'}^\gamma_\alpha.
\end{aligned}
\]
For the triples $(\bt,\alpha,\beta)$ and $(\bt,\alpha,\gamma)$ the maps ${j'}^\beta_\alpha$ and ${j'}^\gamma_\alpha$ are surjective, and thus
${j'}^\beta_\alpha F_\beta A = G_{\beta\smin\alpha}A$ and ${j'}^\gamma_\alpha F_\gamma A = G_{\gamma\smin\alpha}A$. This implies
$f^\gamma_\alpha i'=if^\beta_\alpha$ and $f^\gamma_\beta j'=jf^\gamma_\alpha$.
As before, $G_{\gamma\smin\alpha}A$ is the cone of the component $\sfS G_{\gamma\smin\beta}A\to G_{\beta\smin\alpha}A$ of $\dff_A$.

{\it Step (v):}
By Step (iv) and Lemma \ref{chainsysmorphlemma}(i), the $K$-projective resolutions $f^\beta_\alpha$ commute with all maps $\ell$ in $\CdR$, i.e.\ they form a strict morphism $\bfC(A,\dff_A)\to\bfC$, cf.\ Defn.\ \ref{chainsysmorph}(ii), whose components are quasi-isomorphisms. By Lemma \ref{chainsysmorphlemma}(ii), the strict morphism is an isomorphism $\bfC(A,\dff_A)\xrightarrow[]{{\sim}}\bfC$ of chain systems in $\DdR$. In particular, the diagram
\begin{equation}
    \label{diagforall}
\begin{tikzcd}[column sep=large, row sep=large]
G_{\beta\setminus\alpha}A \arrow[r, "i'"] \arrow[d, "f^{\beta}_{\alpha}"',"\sim"] & G_{\gamma\setminus\alpha}A \arrow[r, "j'"] \arrow[d, "f^{\gamma}_{\alpha}"',"\sim"] & G_{\gamma\setminus\beta}A \arrow[r, "k'"] \arrow[d, "f^{\gamma}_{\beta}","\sim"'] & \sfS G_{\beta\setminus\alpha}A \arrow[d, "\sfS f^{\beta}_{\alpha}","\sim"'] \\
C^{\beta}_{\alpha} \arrow[r, "i"]                                            & C^{\gamma}_{\alpha} \arrow[r, "j"]                                             & C^{\gamma}_{\beta} \arrow[r, "k"]                                           & \sfS C^{\beta}_{\alpha}
\end{tikzcd}
\end{equation}
commutes for all $\alpha\subset\beta\subset\gamma$, where the first two squares commute in $\CdR$, and the third square commutes in $\DdR$ by the naturality of the connecting homomorphisms of homotopy-exact sequences, cf.\ Lem.\ \ref{natforhomexact}.
    \qed
\end{proof}

\begin{remark}
\label{gradedversion}
Theorem \ref{mainres987655} has a graded version. If $\bfC$ is an excisive, homotopy-exact chain system of chain complexes, i.e.\ $\Z$-graded differential modules with differentials of degree $-1$, and chain maps of degree $0$, then the proof applies verbatim, with $\sfS$ the shift, $(\sfS X)_n=X_{n-1}$ and $\dff_{\sfS X}=-\dff_X$. Moreover, $\sfS^{-1}Q$ in place of $\sfS Q$ as the source of the chain map $\lambda$ in Theorem~\ref{SESinDR}, which is then a map $Q\to P$ of degree $-1$, so that $\cone(\lambda)=P\oplus Q$ (in the ungraded case $\sfS^{-1}=\sfS$).  In particular, $\lambda_q\colon\sfS^{-1}G_qA\to F_{\beta^\dagger}A$ in Step~(ii) of the proof. 
Indeed, $K$-projective resolutions by complexes of projective modules exist, cf.\ \cite{Spalt}, and mapping cones, Proposition \ref{prop:kprojequiv} and Theorem \ref{SESinDR} respect the grading. The result is a projective $\sP$-graded chain complex $A=\bigoplus_{p\in\sP}G_pA$, where each $G_pA$ is a chain complex and $\dff_A$ has degree $-1$, and $K$-projective resolutions $f^\beta_\alpha$ of degree $0$. This is the graded setting of the (C-)connection matrices in \cite{fran}, cf.\ Rmk.~\ref{franzosaresult}.
\end{remark}

\section{Chain generated spectral systems}
\label{specsystems}
The analogue of a chain system on homology, without assuming it is induced by chains, is called a {Cartan--Eilenberg} system and is defined as follows:
\begin{definition}[cf.\ \cite{SpV4}]
    \label{CEsys}
A \emph{{Cartan--Eilenberg} system}
over a poset $\sQ$ consists of a covariant functor $\sE\colon \IIi\to \sRmodR$   and a natural transformation $\sk\colon \sE \spi_2 \Rightarrow \sE\spi_0$ between the composite functors $\sE\spi_2$ and $\sE\spi_0$,
called the \emph{connecting homomorphism},
such that
\[
\begin{tikzcd}[column sep=small]
\sE \spi_0 \arrow[rr, "\sE\imath", Rightarrow] &                               & \sE\spi_1\arrow[ld, "\sE\jmath", Rightarrow] \\
                              & \sE\spi_2 \arrow[lu, "\sk", Rightarrow] &
\end{tikzcd}
\]
is an exact triangle, where the natural transformations $\sE\imath$ and $\sE\jmath$ are the right whiskerings of $\sE$ and $\imath$, and $\sE$ and $\jmath$ respectively.
A {Cartan--Eilenberg} system over $\sQ$ is denoted by $\bfE = \bigl(\IIi,\sE,\sk\bigr)$.
\end{definition}
If we  unpack the above definition we obtain $E$-terms $E^\beta_\alpha = \sE(\alpha,\beta)\in \sRmodR$, $(\alpha,\beta)\in \IIi$, with additional properties:
\begin{enumerate}[label=\textup{(\roman*)},leftmargin=2.2em,itemsep=0.1em,topsep=0pt]
    \item [(i)]   homomorphisms    $\ell\colon E^\beta_\alpha \to E^{\beta'}_{\alpha'}$ {for all $(\alpha,\beta)\le (\alpha',\beta')$, such that $\ell'\ell=\ell''$ for $(\alpha,\beta)\le (\alpha',\beta')\le (\alpha'',\beta'')$, and $\ell=\id$ for $(\alpha,\beta) = (\alpha',\beta')$};
 \item[(ii)] connecting homomorphisms $k\colon E^\gamma_\beta\to E^\beta_\alpha$ and
$k'\colon E^{\gamma'}_{\beta'}\to E^{\beta'}_{\alpha'}$ such that the diagrams
\begin{equation}
\label{unpackeddiag12}
\begin{tikzcd}[column sep=small]
E^\beta_\alpha \arrow[rr, "i"] &                   & E^\gamma_\alpha \arrow[ld, "j"] \\
                  & E^\gamma_\beta \arrow[lu, "k"] &
\end{tikzcd}
\qquad\qquad
\begin{tikzcd}
E^\gamma_\beta \arrow[r, "k"] \arrow[d, "\ell"'] & E^\beta_\alpha \arrow[d, "\ell"] \\
E^{\gamma'}_{\beta'} \arrow[r, "k'"]                 & E^{\beta'}_{\alpha'}
\end{tikzcd}
\end{equation}
are exact and commutative respectively for all   ordered triples $(\alpha,\beta,\gamma)\le (\alpha',\beta',\gamma')$;
\item[(iii)] {for $\sQ=\sO(\sP)$,} a Cartan--Eilenberg system is \emph{excisive} if $\ell\colon E^\alpha_{\alpha\cap\beta}\xrightarrow[]{\cong} E^{\alpha\cup \beta}_{\beta}$, for all $\alpha,\beta\in \sO(\sP)$.
\end{enumerate}
By construction $\ell\colon E^\beta_\alpha \xrightarrow[\id]{\ell} E^\beta_\alpha$ is the identity homomorphism and the exactness of \eqref{unpackeddiag12}[left] shows that $E^\alpha_\alpha =0$ for all $\alpha \in \sO(\sP)$, cf.\ \cite{HelleRognes}.
The $E$-terms $E^\beta_\alpha$ with $\beta\smin\alpha=\{p\}$ are called \emph{basic} $E$-terms. For an excisive Cartan--Eilenberg system the basic $E$-terms with the same $p\in\sP$ are isomorphic via excision, cf.\ \cite[Lem.\ 4.3]{SpV4}.
A morphism between  {Cartan--Eilenberg} systems $\bfE'$ and $\bfE$ is given by homomorphisms  $h^{\beta}_{\alpha}\colon E'^\beta_\alpha \to E^\beta_\alpha$ which
satisfy the commutative diagram:
\begin{equation}
\label{dia:exacttriangles}
    \begin{tikzcd}[row sep=large, column sep=large]
{} \arrow[r] & E'^\beta_\alpha \arrow[r, "i"] \arrow[d, "h^\beta_\alpha"] & E'^\gamma_\alpha \arrow[r, "j"] \arrow[d, "h^{\gamma}_{\alpha}"] & E'^\gamma_\beta \arrow[r, "k"] \arrow[d, "h^{\gamma}_{\beta}"] & E'^\beta_\alpha \arrow[r] \arrow[d, "h^{\beta}_{\alpha}"] & {} \\
{} \arrow[r] & E^\beta_\alpha \arrow[r, "i'"]                           & E^\gamma_\alpha \arrow[r, "j'"]                                & E^\gamma_\beta \arrow[r, "k'"]                               & E^\beta_\alpha \arrow[r]                                & {}
\end{tikzcd}
\end{equation}
for any $\alpha\subset\beta\subset\gamma$. In this case we write $h\colon \bfE'\to \bfE$ as a morphism of  {Cartan--Eilenberg} systems.

A spectral system $\bfS = (\IIfouri,\sS,\dff)$ is defined by the functor
\[
\sS\colon \IIfouri \longrightarrow \sRmodR, \quad {(\alpha,\beta,\gamma,\delta)} \mapsto \sS(\alpha,\beta,\gamma,\delta)\footnote{Our indexing of the S-terms deviates from the indexing used in  \cite{Matschke}.}={S^{\gamma\delta}_{\alpha\beta}}\in \sRmodR,
\]
and natural transformation $\dff$, satisfying
\begin{enumerate}[label=\textup{(\roman*)},leftmargin=2.2em,itemsep=0.1em,topsep=0pt]
    \item[(i)] differentials
$\dff\colon S_{\gamma\delta}^{\epsilon\zeta} \to S_{\alpha\beta}^{\gamma\delta}$ for all sextuples
$(\alpha,\beta,\gamma,\delta,\epsilon,\zeta)\in \IIsixi$,
and such that
\begin{equation}
    \label{homdef1}
S_{\beta\delta}^{\epsilon\eta} \cong H\bigl(S_{\gamma\delta}^{\epsilon\zeta}\bigr) = \frac{\ker\colon S_{\gamma\delta}^{\epsilon\zeta}\xrightarrow[]{\dff} S_{\alpha\beta}^{\gamma\delta} }{\image\colon S_{\epsilon\zeta}^{\eta\theta}\xrightarrow[]{\dff} S_{\gamma\delta}^{\epsilon\zeta}},\quad\text{(homology)},
\end{equation}
for all octuples  $(\alpha,\beta,\gamma,\delta,\epsilon,\zeta,\eta,\theta)\in \IIeighti$, where $S_{\epsilon\zeta}^{\eta\theta}\xrightarrow[]{\dff} S_{\gamma\delta}^{\epsilon\zeta} \xrightarrow[]{\dff} S_{\alpha\beta}^{\gamma\delta}$ satisfies $\dff^2=0$. The isomorphism in \eqref{homdef1} is induced by the homomorphisms $\ell$ in (ii), i.e.\ $\ell\colon S^{\epsilon\zeta}_{\beta\delta}\to S^{\epsilon\zeta}_{\gamma\delta}$ is injective with image $\ker\dff$, and $\ell'\colon S^{\epsilon\zeta}_{\gamma\delta}\to S^{\epsilon\eta}_{\gamma\delta}$ is surjective with kernel $\image\dff$;
\item [(ii)]  commutative diagrams
\begin{equation}
   \label{commdwithl}
\begin{tikzcd}
S^{\epsilon\zeta}_{\gamma\delta} \arrow[r, "\dff"] \arrow[d, "\ell"'] & S^{\gamma\delta}_{\alpha\beta} \arrow[d, "\ell"] \\
S^{\epsilon'\zeta'}_{\gamma'\delta'} \arrow[r, "\dff"]                 & S^{\gamma'\delta'}_{\alpha'\beta'}
\end{tikzcd}
\end{equation}
for all sextuples $(\alpha,\beta,\gamma,\delta,\epsilon,\zeta)\le (\alpha',\beta',\gamma',\delta',\epsilon',\zeta')$ in $\IIsixi$;
\item [(iii)] short exact sequences
\begin{equation}
    \label{shortexactSt}
\begin{tikzcd}
0 \arrow[r] & S_{\alpha\beta}^{\gamma\epsilon} \arrow[r, "i"] & S_{\alpha\beta}^{\delta\epsilon} \arrow[r, "j"] & S_{\alpha\gamma}^{\delta\epsilon} \arrow[r] & 0,
\end{tikzcd}
\end{equation}
for all quintuples $(\alpha,\beta,\gamma,\delta,\epsilon)\in \IIfivei$.\footnote{In the special case $S_{\alpha\beta}^{\gamma\epsilon}\to S_{\alpha\beta}^{\delta\epsilon}$ and $S_{\alpha\beta}^{\delta\epsilon} \to S_{\alpha\gamma}^{\delta\epsilon}$ the homomorphisms $\ell$ are denoted by $i$ and $j$ respectively.}
\end{enumerate}
The functoriality implies that there are homomorphisms $\ell\colon S^{\gamma\delta}_{\alpha\beta} \to S^{\gamma'\delta'}_{\alpha'\beta'}$ for all quadruples $(\alpha,\beta,\gamma,\delta)\le (\alpha',\beta',\gamma',\delta')$ in $\IIfouri$, such that $\ell'\ell=\ell''$ for $(\alpha,\beta,\gamma,\delta)\le (\alpha',\beta',\gamma',\delta')\le (\alpha'',\beta'',\gamma'',\delta'')$, and $\ell=\id$ for $(\alpha,\beta,\gamma,\delta) = (\alpha',\beta',\gamma',\delta')$.
Every spectral system is generated by a {Cartan--Eilenberg} system via:
\begin{equation}
    \label{Sterms}
S_{\alpha\beta}^{\gamma\delta}=\displaystyle  \frac{\ker \colon E^{\gamma}_{\beta} \xrightarrow{k} E^{\beta}_{\alpha}}{\image \colon E^{\delta}_{\gamma} \xrightarrow{k} E^{\gamma}_{\beta}},\quad \forall (\alpha,\beta,\gamma,\delta)\in \bm{\sQ_4},
\end{equation}
i.e. $\bfS = \bfS(\bfE)$, v.\ \cite{CastVander1} for more details.

\subsection{Chain generated systems}
In general a {Cartan--Eilenberg} system is not a priori given by the homology of differential modules. In the special {case that} the $E$-terms are given by homology we use the following definition:
\begin{definition}
\label{def:chaingenerated}
A {Cartan--Eilenberg} system $\bfE$ over $\sO(\sP)$  is \emph{chain generated} if there exists a homotopy-exact chain system $\bfC$ over
$\sO(\sP)$, cf.\ Defn.\ \ref{hexrows}, such that $\bfE\cong\bfE(\bfC)$, where $\bfE(\bfC)$ is the Cartan--Eilenberg system with $E$-terms $H\big(C^{\beta}_{\alpha}\big)$, $(\alpha,\beta)\in \IIi$, and with the structure maps $i,j,k$ induced by those of $\bfC$ on
homology. We call  $\bfC$ a \emph{chain model} for $\bfE$.
\end{definition}
We require homotopy-exact chain systems in order to apply Theorem \ref{mainres987655}. Chain systems given by short exact or weakly exact sequences, by lattice-filtered differential modules, and by chain complexes of filtered spaces are of this type.
If we use Theorem \ref{mainres987655} we obtain a representation theorem for chain generated {Cartan--Eilenberg} systems and therefore for spectral systems.
A $\sP$-graded differential module $(A,\dff_A)$ defines a $\sP$-split chain system $\bfC(A,\dff_A)$. A {Cartan--Eilenberg} system generated by $(A,\dff_A)$ is given by $\bfE(A,\dff_A) := \bfE\big(\bfC(A,\dff_A)\big)$. The $\sP$-graded differential module $(A,\dff_A)$ is called a \emph{differential model} for $\bfE$ if $\bfE\cong\bfE(A,\dff_A)$.
\begin{theorem}
\label{thm:chaingeneratedmodel}
Let $R$ be an   associative unital ring, $\sP$ a finite poset, and let $\bfE$ be an excisive,
chain generated {Cartan--Eilenberg} system over $\sO(\sP)$. Then,
there exists a projective, $\sP$-graded
 differential
module $(A,\dff_A)$ --- a projective differential model ---, such that
\[
  \bfE \;\cong\; \bfE(A,{\dff_A}).
\]
{In other words}, every chain model $\bfC$ for $\bfE$ can be replaced by a
projective differential model for $\bfE$.
\end{theorem}
\begin{proof}
By assumption $\bfE \cong \bfE(\bfC)$ for some chain system $\bfC$ with homotopy-exact rows. The latter is excisive in $\DdR$ since $\bfE$ is excisive, cf.\ Lem.\ \ref{isoisqis}. This implies that the conditions in Theorem \ref{mainres987655} are satisfied and therefore there exists a $\sP$-split chain system $\bfC(A,\dff_A)$, for some projective, $\sP$-graded differential module $(A,\dff_A)$, and an isomorphism $\bfC(A,\dff_A)\to\bfC$ of chain systems. On homology it induces an isomorphism of Cartan--Eilenberg systems, and thus $\bfE\cong\bfE(\bfC)\cong\bfE(A,\dff_A)$.
    \qed
\end{proof}
{Conversely, for every $\sP$-graded differential module $(A,\dff_A)$ its associated chain system $\bfC(A,\dff_A)$ is an excisive chain model for $\bfE(A,\dff_A)$. For $\bfC(A,\dff_A)$ the rows are short exact and the excision maps are identity maps, cf.\ Sect.~\ref{chainsystems}. Hence a Cartan--Eilenberg system has a differential model if and only if it is excisive and chain generated. Moreover, it allows a projective model.}

\begin{remark}
    \label{franzosaresult}
    Chain generated  {Cartan--Eilenberg} systems   arise in
practice  in Conley index theory, where  $\bfC$ is the system of (co)chain complexes of
the index pairs of a Morse representation, which has short exact rows,
and $\bfE$ is the induced system
of Conley indices \cite{fran}. {With the conventions of Section~\ref{hslfordiffmod} the connecting homomorphisms of $\bfE(\bfC)$ are $-\delta$, where $\delta$ are the connecting homomorphisms of the long exact sequences. Replacing every $C^\beta_\alpha$ by $\sfS C^\beta_\alpha$ changes their sign. A more detailed discussion is given below \eqref{constructwithr}.}
    The main result in \cite{fran} on the existence of $C$-connection matrices is a corollary of Theorem \ref{mainres987655}, given as Theorem \ref{thm:chaingeneratedmodel}, in the graded form of Remark~\ref{gradedversion}.  The latter extends the result in \cite{fran} by removing
    any restrictions on the coefficient ring. The main ingredients are the derived horseshoe lemma in $\DdR$ and the lattice-ordered induction. The output is a poset-graded differential module $(A,\dff_A)$, where $\dff_A$ plays the role of $C$-connection matrix. If the basic $E$-terms are free one obtains a \emph{connection matrix}, i.e. a strict  differential, cf.\ Cor.\ \ref{mainres987655-22}.
\end{remark}

\begin{theorem}
    \label{chaingensemisimple}
    Let $R$ be a semi-simple ring { and $\sP$ a finite poset}. Then, every excisive, chain generated {Cartan--Eilenberg} system $\bfE$ over $\sO(\sP)$  allows a projective, strict differential model.
\end{theorem}
\begin{proof}
Over a semi-simple ring every module is projective. In particular, the basic $E$-terms are projective and Corollary \ref{mainres987655-22} applies.
    \qed
\end{proof}

\subsection{A counterexample to chain generation}
\label{counterex12}
Consider the  ring $R={\Z/4}$, which is non-hereditary.
Let $\sP={\{1<2\}}$ and let $\bfE$ be a {Cartan--Eilenberg} system over $\sO(\sP)$ in the category of {$\Z/4$}-modules.
Over $\sO(\sP) = \{\bt,\alpha,\tp \}$, with $\alpha=\{1\}$, the {Cartan--Eilenberg} system $\bfE$ consists of one exact triangle and associated truncated short exact sequence for $E^\smtp_\smbt$:
\[
\begin{tikzcd}[ row sep=large]
E^\alpha_\smbt \arrow[rr, "i"] &                               & E^\smtp_\smbt \arrow[ld, "j"] \\
                                              & E^\smtp_\alpha \arrow[lu, "k"] &
\end{tikzcd}\quad \begin{tikzcd}
0 \arrow[r] & \coker~ k \arrow[r, "i'"] & E^\smtp_\smbt \arrow[r, "j'"] & \ker k \arrow[r] & 0,
\end{tikzcd}
\]
where {all $E$-terms are $\Z/4$, and} $i,j$ and $k$ are multiplication by $2$.
Then, $\ker i = \image~k = \langle 2\rangle$, $\ker j = \image~i = \langle 2\rangle$, and $\ker k = \image~j = \langle 2\rangle$, which establishes the exact triangle.
Observe that {$\langle 2\rangle = 2\Z/4 \cong \Z/2$},
$i'([1]) = 2\in {\Z/4}$,
and $j'(b) = 2b\mod 4$.
Let us start with a proof that $\bfE$ is \emph{not} chain generated based on Theorem \ref{thm:chaingeneratedmodel} and Corollary \ref{mainres987655-22}.
Suppose $\bfE$ is chain generated, then Corollary \ref{mainres987655-22}
 claims that there exists a strict, $\sP$-graded differential module $(A,\dff_A)$, with $A=\bigoplus_{p\in \sP}G_pA$ and
$G_1A=E^\alpha_\smbt={\Z/4}$, $G_2A=E^\smtp_\alpha={\Z/4}$ and $E^\smtp_\smbt
={\Z/4}$. All $E$-terms  are finitely generated projective {$\Z/4$}-modules, making $A = {\Z/4\oplus\Z/4}$ a finitely generated projective {$\Z/4$}-module.
{By Corollary \ref{mainres987655-22} the differential is strict, and since $\bfE(A,\dff_A)\cong\bfE$ its off-diagonal entry is multiplication by $2$ up to a unit of $\Z/4$, i.e., up to isomorphism,}
\[
\dff_A = \begin{bmatrix}
    ~0~& ~k~~\\
    ~0~& ~0~~
\end{bmatrix}.
\]
This implies that $E^\smtp_\smbt\cong \coker~k\oplus \ker k ={\Z/2\oplus \Z/2}$, which contradicts   that $E^\smtp_\smbt={\Z/4}$.
If the above given {Cartan--Eilenberg} system is assumed to be chain generated, then all conditions are met and a differential exists, which contradicts our example.
 We conclude that $\bfE$ is a {Cartan--Eilenberg} system of {$\Z/4$}-modules that is \emph{not} chain generated.

The above counterexample can also be obtained directly.
Suppose $\bfE\cong\bfE(\bfC)$, and let $C^\alpha_\smbt\xrightarrow[]{i}C^\smtp_\smbt\xrightarrow[]{j}C^\smtp_\alpha$ be the homotopy-exact sequence of the triple $(\bt,\alpha,\tp)$. Consider the short exact sequence
\[
0\to C^\smtp_\smbt\xrightarrow[]{i'}\cone(i)\xrightarrow[]{j'}\sfS C^\alpha_\smbt\to 0
\]
of the standard inclusion and projection, cf.\ \eqref{SESforcone}, whose connecting homomorphism is $H(i)$, cf.\ \eqref{LESinHom}. Since $\rmTheta i'=j$ and $\rmTheta\colon\cone(i)\to C^\smtp_\alpha$ is a quasi-isomorphism, $H(i')$ corresponds to $j$, and thus $\image H(i')=2H\bigl(\cone(i)\bigr)$. Let $\eta=[z]$ generate $H(C^\smtp_\smbt)={\Z/4}$. Then $i'(z)=2w+\dff_i u$ for a cycle $w$ and some $u$ in $\cone(i)$, and since $4=0$ in {$\Z/4$}, $i'(2z)=\dff_i(2u)$. Hence $f=j'(2u)$ is a cycle, and by the construction of the connecting homomorphism, $H(i)[f]=[2z]=2\eta\neq 0$. On the other hand $2f=j'(4u)=0$, and since $H(i)$ is multiplication by $2$, $H(i)[f]=0$. This contradiction shows again that $\bfE$ is not chain generated. The argument uses neither Theorem \ref{thm:chaingeneratedmodel} nor projectivity of the $E$-terms, only $4=0$ and the homotopy-exactness of one row.
\begin{remark}
    \label{notcgsys}
The counterexample shows that  {over the non-hereditary ring $\bR=\Z/4$ there are Cartan--Eilenberg systems} that are not chain generated.
In \cite{CastVander2} we show that for semi-simple rings, excisive {Cartan--Eilenberg} systems over $\sO(\sP)$, with $\sP$ a linear order,
are always chain generated. For arbitrary poset $\sP$ the answer is still open.
Another class of {Cartan--Eilenberg} systems are those generated by  {excisive} chain systems $\bfC$ in $\DdR$. The question whether such {Cartan--Eilenberg} systems are chain generated, i.e. there exists a homotopy-exact chain system $\bfC'$ such that $\bfE(\bfC)\cong \bfE(\bfC')$, is open.
\end{remark}

\section{Special cases}
\label{specialcases}
In the proof of Theorem \ref{mainres987655} we use the existence of $K$-projective resolutions. In some special cases additional properties of these choices can be utilized to obtain more refined statements.

\subsection{Projective homology}
\label{projbasicE}
Let $(A,\dff_A)$ be a differential module with homology $H(A)$.
\begin{proposition}
    \label{resviahom}
    Suppose $H(A)$ is a projective module. Then, there exists a chain map
    \[
    \pi\colon \big(H(A),0 \big)\to (A,\dff_A),
    \]
which is a  $K$-projective resolution of $A$.
\end{proposition}
\begin{proof}
For homology we have the short exact sequence
\begin{equation}
    \label{sesforimdA}
\begin{tikzcd}[column sep=large, row sep=large]
0 \arrow[r] & \image \dff_A \arrow[r, "i","\subset"'] & \ker \dff_A \arrow[r, "j", two heads] & H(A) \arrow[r] \arrow[l, "\pi"', dotted, bend left] & 0
\end{tikzcd}
\end{equation}
Since $H(A)$ is projective the following diagram yields a lift $\pi\colon {H(A)}\to \ker\dff_A$:
\[
\begin{tikzcd}[column sep=large, row sep=large]
                                                & \ker\dff_A \arrow[d, "j", two heads] \\
{H(A)} \arrow[r, "\id"] \arrow[ru, "\pi", dotted] & {H(A)}
\end{tikzcd}
\]
The short exact sequence in \eqref{sesforimdA} is split and therefore $\ker\dff_A = \image\dff_A \oplus \pi\big(H(A)\big)$. The map $\pi$ is a chain map since $\dff_A\pi=0$.
By definition $H(\pi)[a] = \big[\pi[a]\big] =j\pi[a]= [a]$, which makes $\pi$ a quasi-isomorphism. With the trivial flag {$F_0H(A)=0$ and $F_kH(A)=H(A)$, for $k\ge 1$}, we conclude that  $\big(H(A),0\big)$ is a projective flag, which is $K$-projective by Proposition~\ref{Kprojprop}, and therefore $\pi$ is a $K$-projective resolution for $A$.
    \qed
\end{proof}
Using Proposition \ref{resviahom} we have {the} following special case of Theorem \ref{mainres987655}:

\begin{theorem}
    \label{mainres987655proj}
    Let $R$ be an   associative unital ring and $\sP$ a finite poset.
    Every excisive, homotopy-exact chain system $\bfC$
    {over $\sO(\sP)$},
    with the property that $H(C^\beta_\alpha)$, $\beta\setminus\alpha=\{p\}$, is projective for all $p\in \sP$,
     is isomorphic to a $\sP$-split chain system $\bfC(A,\dff_A)$, for some projective, $\sP$-graded, strict\footnote{For a strict differential, cf.\ Eqn.\ \eqref{strictdiff12}, the diagonal entries are trivial.}  differential module $(A,\dff_A)$, with $G_pA\cong H(C^\beta_\alpha)$, $\beta\setminus\alpha=\{p\}$, for all $p\in\sP$.
\end{theorem}
\begin{proof}
The only alterations of the proof of Theorem \ref{mainres987655} occur  in Step (ii), by choosing the $K$-projective resolutions
$f''\colon\big(H(C^\beta_{\beta^\dagger}),0\big)\to C^\beta_{\beta^\dagger}$ given by Proposition \ref{resviahom}.
Only the pairs $\beta=\big\downarrow q$, $\beta^\dagger=\beta\smin\{q\}$ are used, and by excision $H(C^\beta_{\beta^\dagger})\cong H(C^{\beta'}_{\alpha'})$ for every pair with $\beta'\smin\alpha'=\{q\}$, so it suffices to assume projectivity for these pairs.
This implies in particular that $(A,\dff_A)$ constructed in the proof of Theorem \ref{mainres987655} is given by
\[
A = \bigoplus_{p\in \sP} G_p A, \quad G_pA = H(C^\beta_\alpha), ~~\beta=\big\downarrow p, ~~\alpha = \beta\setminus\{p\},~~\forall p\in \sP,
\]
and $\dff_A$ is a strict differential, i.e. $\dff_A(p,p)=0$ for all $p\in \sP$.
    \qed
\end{proof}

\begin{corollary}
    \label{mainres987655-22}
    {Let $R$ be an associative unital ring and $\sP$ a finite poset.}
    Every excisive, chain generated {Cartan--Eilenberg} system $\bfE$,
    with the property that the basic $E$-terms $E^\beta_\alpha$, $\beta\setminus\alpha=\{p\}$, are projective for all $p\in \sP$, is isomorphic to a Cartan--Eilenberg system $\bfE(A,\dff_A)$, for some projective,
$\sP$-graded, strict differential module $(A,\dff_A)$, with $G_pA\cong E^\beta_\alpha$, $\beta\setminus\alpha=\{p\}$, for all $p\in\sP$.
\end{corollary}
\begin{proof}
Let $\bfC$ be a chain model for $\bfE$. As in the proof of Theorem \ref{thm:chaingeneratedmodel}, the chain model $\bfC$ is excisive in $\DdR$, and $H(C^\beta_\alpha)\cong E^\beta_\alpha$. {Theorem \ref{mainres987655proj} yields an isomorphism $\bfC(A,\dff_A)\to\bfC$ of chain systems, and, as in the proof of Theorem \ref{thm:chaingeneratedmodel}, $\bfE\cong\bfE(\bfC)\cong\bfE(A,\dff_A)$, with $G_pA\cong H(C^\beta_\alpha)\cong E^\beta_\alpha$, $\beta\setminus\alpha=\{p\}$.}
    \qed
\end{proof}

\subsection{Hereditary rings}
\label{heredrings12}
A case in which projective flags are easily found is when the ring $R$ is (left) hereditary, i.e. rings for which all submodules of projective modules are again projective, cf.\ \cite[Ch.~I, \S 5]{CE}. In this text hereditary means left hereditary, and Noetherian means left Noetherian.
\begin{lemma}
    \label{heredinderived}
    Let $R$ be a hereditary ring. Then, $(A,\dff_A) \sim \big(H(A),0\big)$ in $\DdR$.
\end{lemma}
\begin{proof}
Let $\pi\colon (P,\dff_P)\to (A,\dff_A)$ be a $K$-projective resolution for $A$, cf.\ Prop.~\ref{existcofrepl}.
Consider the short exact sequence:
\[
\begin{tikzcd}[column sep=large, row sep=large]
0 \arrow[r] & \ker \dff_P \arrow[r, "\subset"] & P \arrow[r, "\dff_P", two heads] & \image \dff_P \arrow[r] \arrow[l, "\psi"', dotted, bend left]  & 0
\end{tikzcd}
\]
The image $\image\dff_P\subset P$ is a submodule of the projective module $P$, and since $R$ is hereditary, also $\image\dff_P$ is projective. Therefore, the sequence splits with a section $\psi\colon\image\dff_P\to P$, and
$P=\ker\dff_P \oplus \psi(\image\dff_P)$. Let $\varphi$ be the composition
$P \xrightarrow[]{\cong}\ker\dff_P\oplus \image\dff_P \xrightarrow[]{\varpi} H(P)$, where $\varpi(z,b)=[z]$. Since $\dff_P x\in\image\dff_P\subset\ker\dff_P$ for all $x\in P$, we have $\varphi\dff_P=0$, and thus
$\varphi\colon (P,\dff_P)\to \big(H(P),0\big)$ is a chain map with $H(\varphi)=\id$. The map $\varphi$ is a quasi-isomorphism. Consequently, in $\DdR$, $(A,\dff_A) \sim (P,\dff_P) \sim \big(H(P),0\big) \cong \big(H(A),0\big)$, which proves the statement.
    \qed
\end{proof}
The property $(A,\dff_A) \sim \big(H(A),0\big)$ in $\DdR$ means that the differential module $(A,\dff_A)$ is \emph{formal}.
For hereditary rings one may look for projective resolutions of the homology as opposed to {resolutions of the differential module itself}.
\begin{proposition}
    \label{smallresol12}
    Let $R$ be a hereditary ring, and let
    \[
\begin{tikzcd}[column sep=large, row sep=large]
0 \arrow[r] & P_1 \arrow[r, "\dff_1"] & P_0 \arrow[r, "\epsilon", two heads] & H(A) \arrow[r]  & 0
\end{tikzcd}
\]
be a projective resolution for $H(A)$. If we define $P=P_0\oplus P_1$, with differential $\dff_P(x,y)=\big( \dff_1 y,0\big)$,
then,
\begin{enumerate}[label=\textup{(\roman*)},leftmargin=2.2em,itemsep=0.1em,topsep=0pt]
    \item [{\rm (i)}] $(P,\dff_P)$ is a projective flag with $H(P)\cong H(A)$;
    \item[{\rm (ii)}] there exists a chain map $\pi\colon (P,\dff_P)\to (A,\dff_A)$, which is a quasi-isomorphism,
\end{enumerate}
and therefore $\pi\colon P\to A$ is a $K$-projective resolution for $A$.
\end{proposition}
\begin{proof}
The stated resolution exists: for a surjection $\epsilon\colon P_0\to H(A)$ with $P_0$ projective, the kernel $P_1=\ker\epsilon$ is projective, since $R$ is hereditary.
The filtering $F_0P=0$, $F_1P = P_0$ and $F_kP= P$, $k\ge 2$, has projective quotients $F_1P/F_0P=P_0$ and $F_2P/F_1P= P_1$, and $\dff_P F_kP\subset F_{k-1}P$ for all $k\ge 1$. This proves that $(P,\dff_P)$ is a projective flag. Moreover, $\ker\dff_P=P_0$, since $\dff_1$ is injective, and $\image\dff_P=\image\dff_1$. Hence $H(P)= P_0/\image \dff_1\cong H(A)$, which establishes (i).

The map $\epsilon^\dagger\colon (P,\dff_P)\to \big(H(A),0\big)$, $\epsilon^\dagger(x,y) = \epsilon(x)$, is a chain map because $\epsilon\dff_1=0$, and $H(\epsilon^\dagger)$ is the isomorphism $P_0/\image\dff_1\cong H(A)$ induced by $\epsilon$. Consequently, $\epsilon^\dagger$ is a quasi-isomorphism. By Lemma \ref{heredinderived}, $\big(H(A),0\big)\sim (A,\dff_A)$ in $\DdR$, and thus $P\xrightarrow[]{\epsilon^\dagger}H(A)\sim A$ is an isomorphism in $\DdR$. Since $P$ is a projective flag, it is $K$-projective, cf.\ Prop.\ \ref{Kprojprop}. In combination with Corollary \ref{genchainmaps} this yields a chain map $\pi\colon P\to A$, which represents the isomorphism. Then, $H(\pi)$ is an isomorphism, and since $P$ is projective, $\pi$ is a $K$-projective resolution, which establishes (ii). By Proposition \ref{prop:kprojequiv}(iii), $\pi$ is unique up to homotopy as a representative of this isomorphism in $\DdR$.
\qed
\end{proof}

\begin{theorem}
    \label{mainres987655chain33}
    Let $R$ be a hereditary, Noetherian ring  and $\sP$ a finite poset.
    Every   excisive, homotopy-exact chain system $\bfC$ {over $\sO(\sP)$},
    with the property that   $H(C^\beta_\alpha)$, $\beta\setminus\alpha=\{p\}$, is finitely generated for all $p\in \sP$, is isomorphic to a $\sP$-split chain system $\bfC(A,\dff_A)$, for some projective, finitely generated,  $\sP$-graded differential module $(A,\dff_A)$, where $G_pA=P_0\oplus P_1$ is the differential module of a projective resolution $0\to P_1\to P_0\to H(C^\beta_\alpha)\to 0$, $\beta\setminus\alpha=\{p\}$, by finitely generated modules, cf.\ Prop.~\ref{smallresol12}.
\end{theorem}
\begin{proof}
As before we alter
 the proof of Theorem \ref{mainres987655}   in Step (ii) by choosing  the $K$-projective resolutions $f''\colon G_qA\to C^\beta_{\beta^\dagger}$, $\beta=\big\downarrow q$. Take $G_qA=P_0\oplus P_1$ and $f''=\pi$ from Proposition \ref{smallresol12}, applied to $C^\beta_{\beta^\dagger}$ and a projective resolution $0\to P_1\to P_0\xrightarrow[]{\epsilon} H(C^\beta_{\beta^\dagger})\to 0$. As in the proof of Theorem \ref{mainres987655proj}, only the pairs $\beta^\dagger\subset \beta$ are relevant. Since $H(C^\beta_{\beta^\dagger})$ is finitely generated, we can choose a surjection $\epsilon\colon R^n\to H(C^\beta_{\beta^\dagger})$. Since $R$ is Noetherian, $\ker\epsilon$ is finitely generated, and since $R$ is hereditary, $\ker\epsilon$ is projective. Hence, $P_0=R^n$ and $P_1=\ker\epsilon$ are finitely generated projective modules, and so is $A=\bigoplus_{q\in\sP}G_qA$ ($\sP$ is finite).
\qed
\end{proof}

\begin{corollary}
    \label{mainres987655-33}
    Let $R$ be a hereditary, Noetherian ring and $\sP$ a finite poset.
    Every excisive, chain generated {Cartan--Eilenberg} system ${\bfE}$,
    with the property that the basic $E$-terms $E^\beta_\alpha$, $\beta\setminus\alpha=\{p\}$, are finitely generated for all $p\in \sP$, is isomorphic to a  {Cartan--Eilenberg system} $\bfE(A,\dff_A)$, for some projective, finitely generated,  $\sP$-graded differential module $(A,\dff_A)$, where $G_pA=P_0\oplus P_1$ is the differential module of a projective resolution $0\to P_1\to P_0\to E^\beta_\alpha\to 0$, $\beta\setminus\alpha=\{p\}$, by finitely generated modules.
\end{corollary}
\begin{proof}
As in the proof of Corollary \ref{mainres987655-22} in combination with Theorem \ref{mainres987655chain33}.
    \qed
\end{proof}

\begin{remark}
    Principal ideal domains are automatically hereditary and Noetherian.
For the special case that $R$ is a principal ideal domain, Theorem \ref{mainres987655chain33} and Corollary \ref{mainres987655-33}, in the graded form of Remark \ref{gradedversion}, cf.\ Rmk.\ \ref{franzosaresult}, give the existence result of R. Franzosa, cf.\ \cite{fran}, referred to as a \emph{C-connection matrix}, in which the basic $E$-terms are replaced by free resolutions; over a principal ideal domain projective modules are free. {Lemma \ref{heredinderived} and Proposition \ref{smallresol12} have graded versions with the same proofs, with $P_1$ shifted by one degree.}
\end{remark}

\subsection{Finitely generated differential modules}
\label{finiteforE}
An important  question is when a chain generated {Cartan--Eilenberg} system can be generated by a finitely generated poset-graded differential module. In the previous  subsection we showed that over hereditary, Noetherian rings this is the case for excisive, chain generated  {Cartan--Eilenberg} systems with finitely generated basic $E$-terms, cf.\ Cor.\ \ref{mainres987655-33}. We now consider arbitrary rings $R$. The condition on the basic $E$-terms is then that they admit \emph{finite} resolutions by \emph{finitely generated} projective modules. Over a hereditary, Noetherian ring every finitely generated module admits such a resolution of length at most one, and Theorem \ref{mainres987655chain33} is a special case of Theorem \ref{mainres987655chain44} below. Theorem \ref{mainres987655proj} is not, since there the basic $E$-terms need not be finitely generated.

\begin{theorem}
    \label{mainres987655chain44}
    Let $R$ be  an associative unital ring  and $\sP$ a finite poset.
    Every   excisive, homotopy-exact chain system $\bfC$ {over $\sO(\sP)$},
    with the property that  $H(C^\beta_\alpha)$, $\beta\setminus\alpha=\{p\}$,  {admits a finite resolution by}   finitely generated projective modules, for all $p\in \sP$, is isomorphic to a $\sP$-split chain system $\bfC(A,\dff_A)$, for some projective, finitely generated,  $\sP$-graded differential module $(A,\dff_A)$.
\end{theorem}

\begin{proof}
As in the proof of Theorem \ref{mainres987655proj} the
alteration to the proof of Theorem \ref{mainres987655} is   in Step (ii), where we choose the $K$-projective resolutions $f''\colon G_qA\to C^\beta_{\beta^\dagger}$, $\beta=\big\downarrow q$. Take the finite projective flag resolution of $C^\beta_{\beta^\dagger}$ consisting of finitely generated projective modules given by Proposition \ref{prop:semifreeexistfinite}.
A finite projective flag $P$ with finitely generated projective quotients $F_kP/F_{k-1}P$ is a finitely generated projective module, because each inclusion $F_{k-1}P\subset F_kP$ splits. Therefore, every $G_qA$ is finitely generated and projective, and so is $A=\bigoplus_{q\in\sP}G_qA$, since $\sP$ is finite.
    \qed
\end{proof}

\begin{corollary}
    \label{mainres987655-44}
        Let $R$ be  an associative unital ring  and $\sP$ a finite poset.
    Every  excisive, chain generated {Cartan--Eilenberg} system $\bfE$,
    with the property that the basic $E$-terms $E^\beta_\alpha$, $\beta\setminus\alpha=\{p\}$, allow finite resolutions  consisting of   finitely generated projective modules for all $p\in \sP$, is isomorphic to a Cartan--Eilenberg system $\bfE(A,\dff_A)$, for some projective, finitely generated,  $\sP$-graded differential module $(A,\dff_A)$.
\end{corollary}
\begin{proof}
As in the proof of Corollary \ref{mainres987655-22} in combination with Theorem \ref{mainres987655chain44}.
    \qed
\end{proof}

\section{An example from dynamical systems}
\label{dynexample}
Consider the logistic map $\bm f(u)=\mu u(1-u)$ on the compactified space $X=[-\infty,1]$, with $\bm f(-\infty)=-\infty$. In this example the parameter $\mu$ satisfies $3<\mu<1+\sqrt5$.
The map $\bm f$ defines a dynamical system $\phi$ on the compact metric space $X$ via $u_{n+1}= \bm f(u_n)$.
The fixed points of $\bm f$ are $-\infty$, $0$ and $w=1-1/\mu$. In the given parameter range for $\mu$ there exists a $2$-periodic orbit $\{v,v'\}$, with $v<w<v'$, i.e. $\bm f(v)=v'$ and $\bm f(v')=v$. The fixed point $-\infty$ and the periodic orbit $\{v,v'\}$ are attracting, and the fixed points $0$ and $w$ are repelling, with $\bm f'(0)>1$ and $\bm f'(w)<-1$, cf.\ \cite{Devaney}.
Choose $c>0$, $\mu/4<r_2<1$, $0<r_1<\min\{\bm f(r_2),1/2\}$, and disjoint closed intervals $J_1\ni v$, $J_2\ni v'$ in $(r_1,r_2)\smin\{w\}$ with $\bm f(J_1)\subseteq\mathrm{int}\,J_2$, $\bm f(J_2)\subseteq\mathrm{int}\,J_1$ and $|(\bm f^2)'|<1$ on $J_1\cup J_2$, which is possible since $(\bm f^2)'(v)=4+2\mu-\mu^2\in(0,1)$. The sets
\[
F_\alpha U=[-\infty,-c],\quad F_\beta U=J_1\cup J_2,\quad F_\delta U=[r_1,r_2],\quad
\]
with $F_\gamma U=F_\alpha U\cup F_\beta U$,  $F_\epsilon U=F_\alpha U\cup F_\delta U$,
$F_\bt U=\varnothing$ and $F_\tp U=X$, form a lattice of attracting blocks $\sU$,\footnote{We use that  $F_\alpha U\cap F_\delta U=\varnothing$ and $F_\beta U\subseteq F_\delta U$.}
which is indexed by the down-set lattice $\sO(\sP)= \{\bt,\alpha,\beta,\gamma,\delta,\epsilon,\tp\}$, with $\sP=\{1,2,3,4\}$ with $2<3<4$ and $1<4$, as displayed by the diagrams:
\begin{center}
\begin{tikzpicture}[scale=1.25,baseline=(current bounding box.center),x=1cm,y=1cm,every node/.style={draw,circle,inner sep=0pt,minimum size=1.45em}]
\node (p1) at (-0.95,0) {$1$};
\node (p2) at (0.95,0) {$2$};
\node (p3) at (0.95,1.25) {$3$};
\node (p4) at (0,2.5) {$4$};
\draw (p1)--(p4) (p2)--(p3) (p3)--(p4);
\end{tikzpicture}
\hspace{3em}
\begin{tikzpicture}[scale=1.25,baseline=(current bounding box.center),x=1.6cm,y=0.95cm,every node/.style={draw,rounded corners=3pt,inner sep=2.5pt}]
\node (z) at (-1,0) {$\bt=\varnothing$};
\node (a) at (-2,1) {$\alpha=\{1\}$};
\node (b) at (0,1) {$\beta=\{2\}$};
\node (g) at (-1,2) {$\gamma=\{1,2\}$};
\node (d) at (1,2) {$\delta=\{2,3\}$};
\node (e) at (0,3) {$\epsilon=\{1,2,3\}$};
\node (o) at (0,4) {$\tp=\{1,2,3,4\}$};
\draw (z)--(a) (z)--(b) (a)--(g) (b)--(g) (b)--(d) (g)--(e) (d)--(e) (e)--(o);
\end{tikzpicture}
\end{center}
\vspace{1.5cm}
For the join-irreducible down-sets we define
$M_p=\mathrm{Inv}(F_\eta U\smin F_\zeta U)$, $\eta = \bigl\downarrow p$ and $\zeta=\eta^\pred$, its unique immediate predecessor in the order.
By \cite{lsa} the poset $(\sM,\le)$, with $\sM=\{M_p\mid p\in\sP\}$, is a Morse representation of $\bm f$.
The Morse sets $M_p$ are $M_1=\mathrm{Inv}(F_\alpha U)=\{-\infty\}$, $M_2 =
\mathrm{Inv}(F_\beta U)=\{v,v'\}$ ($\bm f^2$ is a contraction on $J_1$ and on $J_2$),  $M_3= \mathrm{Inv}(F_\delta U\smin F_\beta U)=\{w\}$ (every full orbit in $[r_1,r_2]$ lies in $K=[\bm f^2(1/2),\mu/4]$), and $M_4=\mathrm{Inv}(X\smin F_\epsilon U)=\{0\}$ (monotonicity of the orbits).

The relative singular chain complexes $C^\eta_\zeta=C(F_\eta U,F_\zeta U)$ are chain complexes of $\Z[z]$-modules, with the $z$-action defined by $z\cdot \bm s=\bm f\circ \bm s$, where $\bm s$ is a singular simplex. Then, $z$ acts on $H(F_\eta U,F_\zeta U)$ via the induced homomorphism $\bm f_*$. This yields a homotopy-exact chain system $\bfC$ over $\sO(\sP)$ in $\Cd\bigl(\Z[z]\bigr)$, with short exact rows. In this special case $\sC$ is excisive in $\Cd\bigl(\Z[z]\bigr)$, since for all pairs of down-sets $\zeta,\eta\in\sO(\sP)$ every component of $F_\zeta U\cup F_\eta U$ is contained in $F_\zeta U$ or in $F_\eta U$.
Let $[u]$ denote the singular $0$-simplex at $u$, let $\bm\sigma$ be the segment from $v$ to $v'$ and $\bm\rho$ the segment from $-\infty$ to $v$.
The homological Conley indices of the Morse sets with $z$ as the index map are given by:
\begin{equation}\label{zyindices}
\def\arraystretch{1.2}
\begin{array}{c|c|c|c|c}
p~~ & ~~(F_\eta U,F_\zeta U)~~ & ~~E^\eta_\zeta=H(F_\eta U,F_\zeta U)~~ & \text{~~degree~~} & \text{~~generator}\\ \hline
1~ & (F_\alpha U,\varnothing) & \Z[z]/(z-1) & 0 & [-\infty]\\
2~ & (F_\beta U,\varnothing) & \Z[z]/(z^2-1) & 0 & [v]\\
3~ & (F_\delta U,F_\beta U) & \Z[z]/(z+1) & 1 & \bm\sigma\\
4~ & (X,F_\epsilon U) & \Z[z]/(z-1) & 1 & \bm\rho
\end{array}
\end{equation}
None of these modules is projective, and thus there is no strict differential model.

The pair $(F_\eta U,F_\zeta U)$ is an index pair for the Morse set $M_{\eta\smin\zeta}$, with index map $z$. The homological Conley index of $M_{\eta\smin\zeta}$ is the shift equivalence class of the $\Z[z]$-module $E^\eta_\zeta$, cf.\ \cite{Mrozek90,FranksRich,Szymczak1995}, where $E$ and $E'$ are \emph{shift equivalent} if there are homomorphisms $r\colon E\to E'$ and $s\colon E'\to E$ with $sr=z^m$ and $rs=z^m$ for some $m\ge0$.
The modules given in \eqref{zyindices} are the Conley indices of the Morse sets in a canonical form. Every module $E$ in the same shift equivalence class satisfies $E\otimes \Z[z,z^{-1}]\cong E^\eta_\zeta$. %, where $L=\Z[z,z^{-1}]$.
The same holds for the Cartan--Eilenberg system $\bfE(\bfC)$. Let $\zeta\mapsto F_\zeta U'$ be another lattice of attracting blocks with the same attractors, and with chain system $\bfC'$. Then $\bm f^m(F_\zeta U)\subseteq F_\zeta U'$ and $\bm f^m(F_\zeta U')\subseteq F_\zeta U$ for some $m\ge0$ and all $\zeta\in\sO(\sP)$. The homomorphisms $\bm f^m_*$ induced by $\bm f^m$, in both directions, commute with the maps $i$, $j$ and $k$ of $\bfE(\bfC)$ and $\bfE(\bfC')$, and their compositions are $z^{2m}$. Since $\Z[z,z^{-1}]$ is flat over $\Z[z]$, they yield an isomorphism $z^{-m}\bm f^m_*\colon\bfE(\bfC)\otimes \Z[z,z^{-1}]\to\bfE(\bfC')\otimes \Z[z,z^{-1}]$ of Cartan--Eilenberg systems of $\Z[z,z^{-1}]$-modules, which does not depend on $m$. Since $z$ acts as an automorphism of all $E$-terms of $\bfE(\bfC)$, the homomorphisms $\xi\mapsto\xi\otimes1$ form an isomorphism $\bfE(\bfC)\cong\bfE(\bfC)\otimes \Z[z,z^{-1}]$. Thus $\bfE(\bfC)$ is the Cartan--Eilenberg system of the Morse representation in the sense of shift equivalence, and its basic $E$-terms are the Conley indices of the Morse sets. %; for arbitrary attracting blocks see Appendix~\ref{ciappendix}.

We now define a differential model $(A,\dff_A)$ in the category $\Cd\bigl(\Z[z]\bigr)$.  In doing so we treat the chain complexes as differential modules.
We take $G_pA\cong \Z[z]^2$ for all $p\in \sP$. For example
$G_2A=\Z[z]\nu_0^2\oplus\Z[z]\nu_1^2\cong\Z[z]^2$, with $\nu_0^2$ and $\nu_1^2$ corresponding to the unit vectors $(1,0)$ and $(0,1)$, and similarly for $G_1A$, $G_3A$ and $G_4A$. An element of $G_2A$ is a finite sum $\sum_k(a_kz^k\nu^2_0+b_kz^k\nu^2_1)$, with $a_k,b_k\in\Z$.\footnote{We may regard $G_2A$ as free abelian group  with basis $z^k\nu^2_0$, $z^k\nu^2_1$, $k\ge0$.}
We define the differential $\dff_A$ by:
\begin{equation}
\label{zydiff}
\dff_A=\left[\arraycolsep=4pt\def\arraystretch{1.2}\begin{array}{cc:cc:cc:cc}
0 & z-1 & 0 & 0 & 0 & 0 & -1 & 0\\
0 & 0 & 0 & 0 & 0 & 0 & 0 & 1\\ \hdashline
0 & 0 & 0 & z^2-1 & z-1 & 0 & 1 & 0\\
0 & 0 & 0 & 0 & 0 & -1 & 0 & 0\\ \hdashline
0 & 0 & 0 & 0 & 0 & z+1 & 0 & -1\\
0 & 0 & 0 & 0 & 0 & 0 & 0 & 0\\ \hdashline
0 & 0 & 0 & 0 & 0 & 0 & 0 & z-1\\
0 & 0 & 0 & 0 & 0 & 0 & 0 & 0
\end{array}\right],
\end{equation}
with the blocks $\dff_A(p,q)\colon G_qA\to G_pA$ of \eqref{diffasmatrix} indicated. It holds that $\dff_A^2=0$ and $\dff_A$ has degree $-1$, making $(A,\dff_A)$ a $\sP$-graded chain complex and therefore a $\sP$-graded differential module. Moreover, $\dff_A(1,2)=\dff_A(1,3)=0$, so that $\dff_A$ is upper-triangular with respect to the poset $\sP$. The diagonal blocks are the resolutions $0\to\Z[z]\to\Z[z]\to E^\eta_\zeta\to0$ of the modules in \eqref{zyindices}, cf.\ Prop.\ \ref{smallresol12}.
For the module $G_2A$
the differential $\dff_A$ maps to elements of the form $\sum_kb_kz^k(z^2-1)\nu_0^2$.
The diagonal block on $G_2A$ is the free resolution
\[
0\to\Z[z]\nu_1^2\xrightarrow{~z^2-1~}\Z[z]\nu_0^2\to\Z[z]/(z^2-1)\to0,
\]
so that $H_0(G_2A)=\Z[z]\nu_0^2/(z^2-1)\Z[z]\nu_0^2\cong\Z[z]/(z^2-1)$
and $H_1(G_2A)=0$, since $\Z[z]$ has no zero divisors. This gives the Conley index for $M_2$. The same holds for all Morse sets $M_p$ with $p\in \sP$.

\begin{proposition}\label{zymodel}
There is a strict  morphism (not roofs) $g\colon\bfC(A,\dff_A)\to\bfC$, whose components $g^\eta_\zeta\colon G_{\eta\smin\zeta}A\to C(F_\eta U,F_\zeta U)$, $\zeta\subset\eta$, are $K$-projective resolutions. Consequently, $(A,\dff_A)$ is a projective differential model for $\bfC$, cf.\ Thm.\ \ref{mainres987655}, and $\bfE(A,\dff_A)\cong\bfE(\bfC)$. In dynamical systems terminology $\dff_A$ is a $C$-connection matrix for $\sM$.
\end{proposition}
\begin{proof}
Since $\partial\bm\sigma=[v']-[v]$, $\partial \bm f_\#\bm\sigma=[v]-[v']$, $\partial\bm\rho=[v]-[-\infty]$ and $\partial \bm f_\#\bm\rho=[v']-[-\infty]$, the chains $\bm f_\#\bm\sigma+\bm\sigma$ and $\bm f_\#\bm\rho-\bm\rho-\bm\sigma$ are cycles in the intervals $F_\delta U$ and $X$, respectively. Hence, there are chains $\bm\tau\in C_2(F_\delta U)$ and $\bm\omega\in C_2(X)$ with $\partial\bm\tau=\bm f_\#\bm\sigma+\bm\sigma$ and $\partial\bm\omega=\bm f_\#\bm\rho-\bm\rho-\bm\sigma$. Define the $\Z[z]$-linear map $g\colon A\to C(X)$ by
\[
\begin{aligned}
g(\nu_0^1)&=[-\infty],\ \ g(\nu^2_0)=[v],\ \ g(\nu^3_0)=\bm\sigma,\ \ g(\nu^4_0)=\bm\rho,\\
g(\nu^1_1)&=g(\nu^2_1)=0,\ \ g(\nu^3_1)=\bm\tau,\ \ g(\nu^4_1)=\bm\omega.
\end{aligned}
\]
Using \eqref{zydiff} we verify that $g$ is a chain map. Indeed, $g(\dff_A\nu^1_1)=(z-1)[-\infty]=0$, $g(\dff_A\nu^2_1)=(z^2-1)[v]=0$, $g(\dff_A\nu^3_0)=\partial\bm\sigma$, $g(\dff_A\nu^3_1)=\bm f_\#\bm\sigma+\bm\sigma=\partial\bm\tau$, $g(\dff_A\nu^4_0)=\partial\bm\rho$ and $g(\dff_A\nu^4_1)=\bm f_\#\bm\rho-\bm\rho-\bm\sigma=\partial\bm\omega$. Since $g(F_\zeta A)\subseteq C(F_\zeta U)$ for all $\zeta\in\sO(\sP)$, it induces chain maps $g^\eta_\zeta\colon G_{\eta\smin\zeta}A\to C(F_\eta U,F_\zeta U)$ which commute with the maps $i$ and $j$ in $\Cd\bigl(\Z[z]\bigr)$, i.e.\ a strict morphism, cf.\ Defn.\ \ref{chainsysmorph}(ii). For $\eta=\big\downarrow p$ and $\zeta=\eta^\pred$ the homologies of $G_pA$ are the cyclic modules of \eqref{zyindices}, generated by the class of $\nu^1_0$, $\nu^2_0$, $\nu^3_0$ and $\nu^4_0$ respectively, which $g^\eta_\zeta$ maps to the generator in \eqref{zyindices}. Therefore, the maps $g^\eta_\zeta$ are quasi-isomorphisms.
By excision, the maps $g^\eta_\zeta$ are quasi-isomorphisms for all pairs with $\eta\smin\zeta=\{p\}$, and hence for all pairs $\zeta\subset\eta$, by Theorem\ \ref{SESinDR}(i) and induction on the number of elements of $\eta\smin\zeta$. The modules $G_{\eta\smin\zeta}A$ are free and are finite projective flags, hence $K$-projective, cf.\ Prop.\ \ref{Kprojprop}.
\qed
\end{proof}

Two generators are needed for each $p$, since the Conley indices are not free $\Z[z]$-modules. Under the strict morphism $g$ of Proposition \ref{zymodel}, multiplication by $z^k$ corresponds to the $k$-th iterate: $g(z^k\nu^1_0)=[\bm f^k(-\infty)]=[-\infty]$, $g(z^k\nu^2_0)=[\bm f^k(v)]$, which is $[v]$ for $k$ even and $[v']$ for $k$ odd, $g(z^k\nu^3_0)=\bm f^k_\#\bm\sigma$ and $g(z^k\nu^4_0)=\bm f^k_\#\bm\rho$. The generators of degree $1$ impose the relations that hold for the dynamics. For instance, $\dff_A\nu^2_1=(z^2-1)\nu^2_0$ expresses $\bm f^2(v)=v$, and $g(\nu^2_1)=0$, since $(z^2-1)[v]=[\bm f^2(v)]-[v]=0$ already holds for chains. Likewise, the off-diagonal entry $\dff_A\nu^3_0=(z-1)\nu^2_0$ corresponds to $\partial\bm\sigma=[v']-[v]=(z-1)[v]$.

\begin{remark}\label{zyremark}
The construction of $g$ is Steps (ii) and (iii) of the proof of Theorem \ref{mainres987655} carried out by hand: the off-diagonal entries of $\dff_A$ are forced by the boundaries of $\bm\sigma$, $\bm\rho$, $\bm\tau$ and $\bm\omega$. For $q$ covering $p$ the Conley index of the union of $M_p$, $M_q$ and the orbits connecting $M_q$ to $M_p$ is $H(G_{\{p,q\}}A)$, which is $0$ for $\{1,4\}$, $\Z[z]/(z-1)$ in degree $0$ for $\{2,3\}$, and $\Z[z]/(z^2-1)$ in degree $1$ for $\{3,4\}$. In no case is it the direct sum of the modules in \eqref{zyindices}, in accordance with the connecting orbits from $0$ to $-\infty$, from $w$ to $\{v,v'\}$, and from $0$ to $w$. For $\{3,4\}$ both indices are concentrated in degree $1$, the connecting homomorphism vanishes on homology, and the entry $\dff_A(3,4)$, given by $\nu^4_1\mapsto-\nu^3_0$, records the non-split extension $0\to\Z[z]/(z+1)\to\Z[z]/(z^2-1)\to\Z[z]/(z-1)\to0$, whose class is the non-zero element of ${\rm Ext}^1_{\Z[z]}\bigl(\Z[z]/(z-1),\Z[z]/(z+1)\bigr)\cong\Z/2$. The orbit from $0$ to $w$ is therefore not detected by a connection matrix over $\Z$, nor over $\K[z]$ for a field $\K$ of characteristic $\neq2$, where the extension splits, but it is detected over $\Z[z]$.
\end{remark}
\begin{remark}\label{zylonger}
The diagonal blocks of $\dff_A$ are $2\times2$, since the $E$-terms are free abelian. If for example the map $\mathrm{id}\times \bm f$ on $\mathbb{R}P^2\times X$ is considered,
  with isolating blocks $\mathbb{R}P^2\times F_\zeta U$, then $4\times4$  blocks enter in the differential $\dff_A$.
\end{remark}

\section{A chain system in the derived category that is not $\sP$-split}
\label{csysnothe}
Theorem \ref{mainres987655} states that every excisive, \emph{homotopy-exact} chain system $\bfC$  over $\sO(\sP)$ is isomorphic   to a $\sP$-split chain system $\bfC(A,\dff_A)$, for some projective, $\sP$-graded differential module $(A,\dff_A)$, via {$K$-projective resolutions} $f_\alpha^\beta\colon G_{\beta\smin\alpha}A\to C_\alpha^\beta$.
This raises the following question:
\begin{enumerate}[label=\textup{(\roman*)},leftmargin=2.2em,itemsep=0.1em,topsep=0pt]
    \item [(Q)]
    Let $\bfC$ be an excisive chain system in $\DdR$ over $\sO(\sP)$. Do there exist a $\sP$-graded differential module $(A,\dff_A)$ and isomorphisms $f^\beta_\alpha\colon G_{\beta\smin\alpha}A\to C^\beta_\alpha$ in $\DdR$, for all $\alpha\subset\beta$, such that the three squares of \eqref{diagforall} commute in $\DdR$ for all $\alpha\subset\beta\subset\gamma$?
\end{enumerate}
In this section we show via a counterexample that (Q) does not have an affirmative answer in general, i.e. we construct an excisive chain system over $\sO(\sP)$, with the poset $\sP=\{1<2<3\}$, that is not isomorphic in $\DdR$ to any $\sP$-split chain system. Its Cartan--Eilenberg system is nevertheless chain generated, cf.\ Thm.\ \ref{absnonreal}.
\begin{lemma}
\label{absequiv}
For an excisive chain system $\bfC$ in $\DdR$ over $\sO(\sP)$ the following statements are equivalent:
\begin{enumerate}[label=\textup{(\roman*)},leftmargin=2.2em,itemsep=0.1em,topsep=0pt]
\item[{\rm (i)}]
Question (Q) has an affirmative answer for $\bfC$;
\item[{\rm (ii)}] $\bfC\cong\bfC(A,\dff_A)$ as chain systems in $\DdR$, cf.\ Defn.\ \ref{chainsysmorph}(i), for some $\sP$-graded differential module $(A,\dff_A)$, not necessarily projective;
\item[{\rm (iii)}] $\bfC$ is isomorphic in $\DdR$ to an excisive, homotopy-exact chain system;
\item[{\rm (iv)}] the differential module $(A,\dff_A)$ in {\rm (i)} can be chosen projective, with $K$-projective $G_pA$, and the maps $f^\beta_\alpha$  can be chosen to be $K$-projective resolutions.
\end{enumerate}
\end{lemma}

\begin{proof}
By Definition \ref{chainsysmorph}(i), statements in (i) and (ii) say the same. The first two squares of \eqref{diagforall} imply that the maps $f^\beta_\alpha$ commute with the maps $i$ and $j$, and hence with all maps $\ell$, cf.\  Lem.\ \ref{chainsysmorphlemma}(i).
The third square implies that maps commute with the connecting homomorphisms. Therefore, maps $f_\alpha^\beta$ as  given in (i), yield an isomorphism of chain systems $\bfC(A,\dff_A)\to\bfC$.

(ii)$\Rightarrow$(iii): $\bfC(A,\dff_A)$   is excisive and homotopy-exact. The rows are short exact sequences and the excision maps are identities, cf.\ Sect.\ \ref{chainsystems}.

(iii)$\Rightarrow$(iv): Suppose $\bfC\cong {\tilde\bfC}$, with ${\tilde\bfC}$ excisive and homotopy-exact. By Theorem \ref{mainres987655}, there exists a projective $\sP$-graded differential module  $(A,\dff_A)$, with  $K$-projective $G_pA$, such that $\bfC(A,\dff_A)\cong\tilde\bfC$.
Composing, $\bfC(A,\dff_A)\cong\bfC$. Each component $G_{\beta\smin\alpha}A\to C^\beta_\alpha$ is an isomorphism in $\DdR$ with $K$-projective source, so it is represented by a chain map, cf.\ Prop.\ \ref{prop:kprojequiv}(iii). Since this chain map represents an isomorphism in $\DdR$, it is a quasi-isomorphism, cf.\ Lem.\ \ref{isoisqis}, which yields a $K$-projective resolution.

(iv)$\Rightarrow$(i): A projective $\sP$-graded differential module is $\sP$-graded, and quasi-isomorphisms are isomorphisms in $\DdR$.
\qed
\end{proof}

 Lemma \ref{absequiv} shows that  the answer to Question (Q)  is affirmative exactly for the chain systems $\bfC$ that are isomorphic to an excisive, homotopy-exact chain system, i.e.\ for the chain systems covered by Theorem
 \ref{mainres987655}, up to isomorphism.
In this sense Theorem \ref{mainres987655} is sharp for Question (Q). Theorem \ref{absnonreal}   shows that not every excisive chain system in $\DdR$ is of this type.
The obstruction is discussed in Remark \ref{absnonrealrem}.

For a chain system $\bfC$ in $\DdR$ over $\sO(\sP)$ the functor $\beta\mapsto C^\beta_\smbt$, $\beta\in\sO(\sP)$, with the maps $i^\beta_\alpha\colon C^\alpha_\smbt\to C^\beta_\smbt$, $\alpha\subset\beta$, is called the \emph{base} of $\bfC$, and the rows of the triples $(\bt,\alpha,\beta)$, with the maps $i^\beta_\alpha$, $j^\beta_\alpha$ and $k^\beta_\alpha$, are called the \emph{bottom triangles} of $\bfC$. The base of $\bfC(A,\dff_A)$ consists of the differential modules $F_\beta A$ and the inclusions $F_\alpha A\to F_\beta A$. For  a linearly ordered poset $\sP$ an isomorphism between the bases of two $\sP$-split chain systems extends to an isomorphism of chain systems:

\begin{lemma}
\label{absbaselemma}
Let $\sP$ be a finite, linearly ordered poset, and let $(A,\dff_A)$ and $(B,\dff_B)$ be $\sP$-graded differential modules such that the differential modules $F_\beta B$, $\beta\in\sO(\sP)$, are $K$-projective. Let $f^\beta_\smbt\colon F_\beta B\to F_\beta A$, $\beta\in\sO(\sP)$, be isomorphisms in $\DdR$ which commute in $\DdR$ with the inclusions $F_\alpha B\to F_\beta B$ and $F_\alpha A\to F_\beta A$ for all $\alpha\subset\beta$. Then, there exists an isomorphism of chain systems $f\colon\bfC(B,\dff_B)\to\bfC(A,\dff_A)$ whose components on the base are the given $f^\beta_\smbt$.
\end{lemma}
\begin{proof}
Let $\sP=\{1<\cdots<n\}$. The down-sets are $\alpha_m=\{1,\dots,m\}$, $0\le m\le n$, with $\alpha_0=\bt$, and we denote the inclusions $F_{\alpha_{m-1}}B\to F_{\alpha_m}B$ and $F_{\alpha_{m-1}}A\to F_{\alpha_m}A$ by $\imath$ and $i$. Since the $F_{\alpha_m}B$ are $K$-projective, the morphisms $f^{\alpha_m}_\smbt$ are represented by chain maps, cf.\ Prop.\ \ref{prop:kprojequiv}(iii). We show, by induction on $m$, that these chain maps can be chosen such that $f^{\alpha_m}_\smbt\imath=i f^{\alpha_{m-1}}_\smbt$ in $\CdR$, starting with $f^{\alpha_0}_\smbt=0$. Let $f^{\alpha_{m-1}}_\smbt$ be chosen, and let $g$ be a chain map which represents $f^{\alpha_m}_\smbt$. The chain maps $g\imath$ and $i f^{\alpha_{m-1}}_\smbt$ both represent $f^{\alpha_m}_\smbt\imath=i f^{\alpha_{m-1}}_\smbt$ in $\DdR$, and since $F_{\alpha_{m-1}}B$ is $K$-projective they are chain-homotopic, cf.\ Prop.\ \ref{prop:kprojequiv}(iii), i.e.\ $g\imath-i f^{\alpha_{m-1}}_\smbt=\dff h+h\dff$ for a homomorphism $h\colon F_{\alpha_{m-1}}B\to F_{\alpha_m}A$. As a module, $F_{\alpha_m}B=F_{\alpha_{m-1}}B\oplus G_m B$, and extending $h$ by zero on $G_m B$ gives a homomorphism $H\colon F_{\alpha_m}B\to F_{\alpha_m}A$ with $H\imath=h$. The chain map $g-\dff H-H\dff$ is chain-homotopic to $g$, and therefore represents $f^{\alpha_m}_\smbt$, and $(g-\dff H-H\dff)\imath=g\imath-\dff h-h\dff=i f^{\alpha_{m-1}}_\smbt$. We choose $f^{\alpha_m}_\smbt=g-\dff H-H\dff$.

With these choices, the restriction of the chain map $f^{\alpha_n}_\smbt\colon B\to A$ to $F_{\alpha_m}B$ is $f^{\alpha_m}_\smbt$, so that $f^{\alpha_n}_\smbt(F_\beta B)\subset F_\beta A$ for all $\beta\in\sO(\sP)$. Hence $f^{\alpha_n}_\smbt$ induces chain maps $f^\beta_\alpha\colon G_{\beta\smin\alpha}B\to G_{\beta\smin\alpha}A$ on all quotients, which commute with the maps $i$ and $j$ of $\bfC(B,\dff_B)$ and $\bfC(A,\dff_A)$ in $\CdR$. Thus $f=(f^\beta_\alpha)$ is a strict morphism, cf.\ Defn.\ \ref{chainsysmorph}(ii), whose components on the base are the given $f^\beta_\smbt$. Since these represent isomorphisms in $\DdR$, they are quasi-isomorphisms, cf.\ Lem.\ \ref{isoisqis}. For $\alpha\subset\beta$ the chain maps $f^\alpha_\smbt$, $f^\beta_\smbt$ and $f^\beta_\alpha$ form a morphism between the short exact sequences $0\to F_\alpha\to F_\beta\to G_{\beta\smin\alpha}\to0$ of $B$ and $A$, and the five lemma, applied to the long exact sequences in homology, shows that $f^\beta_\alpha$ is a quasi-isomorphism. By Lemma \ref{chainsysmorphlemma}(ii), $f$ is an isomorphism of chain systems.
\qed
\end{proof}

Let $R$ be a commutative ring with an element $z\in\bR$ such that $z^3=0$ and $z^2\neq 0$. Examples are $\bR=\Z/8$, with $z=2$, $\bR=\Z/p^3$, $p$ a prime, with $z=p$, $\bR=\Z/32$, with $z=4$ and $\bR=\K[z]/(z^3)$, $\K$ a field. Since $z$ is nilpotent, the element $1+z$ is a unit, with inverse $(1+z)^{-1}=1-z+z^2$, and $z^2(1+z)^{-1}=z^2$.
We consider the poset $\sP=\{1<2<3\}$. The lattice of down-sets  is given by
$\sO(\sP)=\{\bt,\alpha,\beta,\tp\}$, with $\alpha=\{1\}$ and $\beta=\{1,2\}$.
Let $A$ be the free $\bR$-module with basis  $\bm\xi,\bm\eta,\bm\zeta$, graded by $\sP$ via $G_1A=\langle {\bm\xi}\rangle$, $G_2A=\langle {\bm\eta}\rangle$ and $G_3A=\langle {\bm\zeta}\rangle$, with upper-triangular differential $\dff_A$ defined by  $\dff_A\bm\xi=0$, $\dff_A\bm\eta=z^2\bm\xi$ and $\dff_A\bm\zeta=z\bm\xi$.
With respect to the basis  $\bm\xi,\bm\eta,\bm\zeta$,
we have the blocks $\dff_A(p,q)\colon G_qA\to G_pA$ of \eqref{diffasmatrix} indicated,
\[
 \dff_A=\left[\arraycolsep=5pt\def\arraystretch{1.3}\begin{array}{c:c:c}
0 & z^2 & z\\ \hdashline
0 & 0 & 0\\ \hdashline
0 & 0 & 0
\end{array}\right],
\]
i.e. $\dff_A(1,2)=z^2$, $\dff_A(1,3)=z$ and $\dff_A(p,q)=0$ otherwise.
Clearly, $\dff_A^2=0$ and $\dff_A$ is strict.
The terms of the chain system $\bfC(A,\dff_A)$ are
\[
\begin{aligned}
C^\alpha_\smbt&=F_\alpha A=\langle\bm\xi\rangle, & C^\beta_\smbt&=F_\beta A=\langle\bm\xi,\bm\eta\rangle, & C^\tp_\smbt&=F_\tp A=A,\\
C^\beta_\alpha&=G_{\beta\smin\alpha}A=\langle\bm\eta\rangle, & C^\tp_\beta&=G_{\tp\smin\beta}A=\langle\bm\zeta\rangle, & C^\tp_\alpha&=G_{\tp\smin\alpha}A=\langle\bm\eta,\bm\zeta\rangle.
\end{aligned}
\]
The only non-zero differentials are $\dff\bm\eta=z^2\bm\xi$ in $F_\beta A$ and in $A$, and $\dff\bm\zeta=z\bm\xi$ in $A$. The other four terms have zero differential and coincide with their homology.
The maps in the chain system $\bfC(A,\dff_A)$ are   $i^\beta_\alpha$,
$j^\beta_\alpha$, $k^\beta_\alpha$ for $(\bt,\alpha,\beta)$, $i^\tp_\alpha$, $j^\tp_\alpha$, $k^\tp_\alpha$ for $(\bt,\alpha,\tp)$, $i^\tp_\beta$, $j^\tp_\beta$, $k^\tp_\beta$ for $(\bt,\beta,\tp)$, and
\[
i\colon G_{\beta\smin\alpha}A\to G_{\tp\smin\alpha}A,\qquad j\colon G_{\tp\smin\alpha}A\to G_{\tp\smin\beta}A,\qquad k\colon G_{\tp\smin\beta}A\to\sfS G_{\beta\smin\alpha}A
\]
for the inner triple $(\alpha,\beta,\tp)$.
The maps $i$  and  $j$ are the inclusions and projections of basis vectors. In the rest of this section the maps of $\bfC(A,\dff_A)$ are written without primes; primes refer to the chain system $\bfC'$ below. The connecting homomorphisms are obtained as follows, cf.\ Sect.\ \ref{hslfordiffmod}. Let $\delta'\subset\delta$ be down-sets. For $F_\delta A=F_{\delta'}A\oplus G_{\delta\smin\delta'}A$, the differential $\dff_A$ is given by
\[
\dff_{F_\delta A}=\begin{bmatrix}\dff_{F_{\delta'}A}&\lambda^\delta_{\delta'}\\ 0&\dff_{G_{\delta\smin\delta'}A}\end{bmatrix},
\]
so that the matrix of $\lambda^\delta_{\delta'}$ is the block of $\dff_A$ with rows in $\delta'$ and columns in $\delta\smin\delta'$.
Since $\dff_{F_\delta A}^2=0$ in $F_{\delta'}A$ it holds that $\dff_{F_{\delta'}A}\lambda^\delta_{\delta'}+\lambda^\delta_{\delta'}\dff_{G_{\delta\smin\delta'}A}=0$, and therefore $\lambda^\delta_{\delta'}\colon\sfS G_{\delta\smin\delta'}A\to F_{\delta'}A$ is a chain map. Consequently, $F_\delta A=\cone(\lambda^\delta_{\delta'})$,
and the connecting homomorphism is given by the chain map $-\sfS\lambda^\delta_{\delta'}$, cf.\
\eqref{firstrot12}.
Since $\sfS$ does not change the underlying homomorphism,
\[
{k}^\delta_{\delta'}(a)=-\lambda^\delta_{\delta'}(a),\qquad a\in G_{\delta\smin\delta'}A.
\]
Similarly, $G_{\tp\smin\alpha}A=G_{\beta\smin\alpha}A\oplus G_{\tp\smin\beta}A=\cone(\lambda)$, with $\lambda=\dff_A(2,3)$, and $k=-\sfS\lambda$.
Summarizing the connecting homomorphisms for the four triples we have:
\begin{center}
\renewcommand{\arraystretch}{1.25}
\begin{tabular}{llll}
\hline
triple & from & off-diagonal block & connecting homomorphism\\
\hline
$(\bt,\alpha,\beta)$ & $\dff_A {\bm\eta}$ & $\lambda^\beta_\alpha( {\bm\eta)=z^2\bm\xi}$ & ${k}^\beta_\alpha( {\bm\eta)=-z^2\bm\xi}$\\
$(\bt,\alpha,\tp)$ & $\dff_A {\bm\eta}$ & $\lambda^\tp_\alpha( {\bm\eta)=z^2\bm\xi}$ & ${k}^\tp_\alpha( {\bm\eta)=-z^2\bm\xi}$\\
 & $\dff_A {\bm\zeta}$ & $\lambda^\tp_\alpha( {\bm\zeta)=z\bm\xi}$ & ${k}^\tp_\alpha( {\bm\zeta)=-z\bm\xi}$\\
$(\bt,\beta,\tp)$ & $\dff_A {\bm\zeta}$ & $\lambda^\tp_\beta( {\bm\zeta)=z\bm\xi}$ & ${k}^\tp_\beta( {\bm\zeta)=-z\bm\xi}$\\
$(\alpha,\beta,\tp)$ & $\dff_{G_{\tp\smin\alpha}A} {\bm\zeta}$ & $\lambda( {\bm\zeta)=0}$ & $k {=0}$\\
\hline
\end{tabular}
\end{center}
For  the chain system $\bfC(A,\dff_A)$  the base is $F_\alpha A\to F_\beta A\to A$, with the bottom triangles of $(\bt,\alpha,\beta)$, $(\bt,\alpha,\tp)$ and $(\bt,\beta,\tp)$, and the inner triangle of $(\alpha,\beta,\tp)$. These form a Verdier octahedron,  i.e.\ a diagram as in the octahedral axiom (TR4), cf.\ \cite[Defn.~10.2.1]{weibel} and \cite{Kunzer}, v.\
Fig.\ \ref{octafig}. Figure \ref{octafigc} shows the octahedron as braid diagram.
\begin{figure}[ht]
\centering
\begin{tikzcd}[column sep=small]
 & & G_{\beta\smin\alpha}A \arrow[rr, "k^\beta_\alpha", bend left] \arrow[rd, "i"', dashed] & & \sfS F_\alpha A \arrow[rd, "\sfS i^\beta_\alpha"] & & \\
 & F_\beta A \arrow[ru, "j^\beta_\alpha"] \arrow[rd, "i^\tp_\beta"'] & & G_{\tp\smin\alpha}A \arrow[ru, "k^\tp_\alpha"'] \arrow[rd, "j"', dashed] & & \sfS F_\beta A \arrow[rd, "\sfS j^\beta_\alpha"] & \\
F_\alpha A \arrow[rr, "i^\tp_\alpha"', bend right] \arrow[ru, "i^\beta_\alpha"] & & A \arrow[rr, "j^\tp_\beta"', bend right] \arrow[ru, "j^\tp_\alpha"] & & G_{\tp\smin\beta}A \arrow[rr, "k"', dashed, bend right] \arrow[ru, "k^\tp_\beta"] & & \sfS G_{\beta\smin\alpha}A
\end{tikzcd}
\caption{ The four strands are the distinguished triangles of $(\bt,\alpha,\beta)$, $(\bt,\alpha,\tp)$, $(\bt,\beta,\tp)$ and $(\alpha,\beta,\tp)$;  the strand of $(\bt,\alpha,\beta)$ is $F_\alpha A\to F_\beta A\to G_{\beta\smin\alpha}A\to\sfS F_\alpha A$, continued by its suspension $\sfS F_\alpha A\to\sfS F_\beta A\to\sfS G_{\beta\smin\alpha}A$. The dashed arrows form the row of $(\alpha,\beta,\tp)$. The four triangles and the two squares enclosed by the strands commute in $\DdR$; the four triangles are the four commutative faces of Figure \ref{octafig}, including the outer face.
For the chain maps of $\bfC(A,\dff_A)$ they commute in $\CdR$, except for the square with vertices $G_{\tp\smin\alpha}A$ and $\sfS F_\beta A$, which commutes up to homotopy. For the chain system $\bfC'$ of Theorem \ref{absnonreal} the  maps $i$ and $j$ are replaced by $i'=\varphi\,i$ and $j'=j\,\varphi^{-1}$.}
\label{octafigc}
\end{figure}
The  distinguished triangles are given by the triples $(\bt,\alpha,\beta)$, $(\bt,\beta,\tp)$, $(\bt,\alpha,\tp)$ and $(\alpha,\beta,\tp)$, and the commutative faces are
\[
i^\tp_\alpha=i^\tp_\beta\, i^\beta_\alpha,\qquad k=(\sfS j^\beta_\alpha)\,k^\tp_\beta,\qquad k^\tp_\alpha\, i=k^\beta_\alpha,\qquad j\,j^\tp_\alpha=j^\tp_\beta .
\]
In addition, the two squares $j^\tp_\alpha i^\tp_\beta=i\,j^\beta_\alpha$ and $k^\tp_\beta\, j=(\sfS i^\beta_\alpha)\,k^\tp_\alpha$ commute; they correspond to the two ways from $F_\beta A$ to $G_{\tp\smin\alpha}A$ through the octahedron, and from $G_{\tp\smin\alpha}A$ back to $F_\beta A$. For the chain maps of $\bfC(A,\dff_A)$ all these identities hold strictly, except for the second square, which holds only up to homotopy, since $k^\tp_\beta j(\bm\eta)=0$ and $(\sfS i^\beta_\alpha)k^\tp_\alpha({\bm\eta})= {-z^2\bm\xi=\dff_{\sfS F_\beta A}\bm\eta}$.

We now define a
chain system $\bfC'$ by using the same octahedron as $\bfC(A,\dff_A)$, except that the inner triangle is transported along an automorphism of $G_{\tp\smin\alpha}A$. Let $\varphi$ be the automorphism of $G_{\tp\smin\alpha}A=\langle\bm\eta,\bm\zeta\rangle$ given by $\varphi(\bm\eta)=\bm\eta-z^2\bm\zeta$ and $\varphi(\bm\zeta)=(1+z)^{-1}\bm\zeta$, with inverse $\varphi^{-1}(\bm\eta)=\bm\eta+z^2\bm\zeta$ and $\varphi^{-1}(\bm\zeta)=(1+z)\bm\zeta$; it is a chain map, since $G_{\tp\smin\alpha}A$ has zero differential. The chain system $\bfC'$ has the same terms, the same connecting homomorphisms and the same base and bottom triangles as $\bfC(A,\dff_A)$, and the maps $i$ and $j$ of the triple $(\alpha,\beta,\tp)$ are replaced by
\[
{i'=\varphi\, i,\quad j'=j\,\varphi^{-1},\quad\text{i.e.}\quad i'(\bm\eta)=\bm\eta-z^2\bm\zeta,\quad j'(\bm\eta)=z^2\bm\zeta,\quad j'(\bm\zeta)=(1+z)\bm\zeta.}
\]
We show that $\bfC(A,\dff_A)$ and $\bfC'$ are two octahedra on the same base with the same bottom triangles, that are not isomorphic as chain systems.  For a chain system $\bfC$ in $\DdR$ the modules $H(C^\beta_\alpha)$
form a Cartan--Eilenberg system $\bfE(\bfC)$. Homology is a functor on $\DdR$ that maps distinguished triangles to exact triangles, cf.\ Defn.\ \ref{defndisttri} and \eqref{SESforcone2}. A (iso)morphism $f$ of chain systems yields a (iso)morphism $H(f)$ of Cartan--Eilenberg systems, cf.\ \eqref{dia:exacttriangles}.

\begin{theorem}
\label{absnonreal}
Let $\bR$ be a  commutative ring with an element $z$ such that $z^3=0$ and $z^2\neq0$, let $(A,\dff_A)$ {be the} $\sP$-graded differential module {defined above,} with $\sP=\{1<2<3\}$, and let
 $\bfC'$ be  the chain system $\bfC(A,\dff_A)$ with the maps $i$ and $j$ of the triple $(\alpha,\beta,\tp)$ replaced by $i'=\varphi\,i$ and $j'=j\,\varphi^{-1}$. Then:
 \begin{enumerate}[label=\textup{(\roman*)},leftmargin=2.2em,itemsep=0.1em,topsep=0pt]
\item[{\rm (i)}] $\bfC'$ is an excisive chain system in $\DdR$ over $\sO(\sP)$, with the same base and the same bottom triangles as $\bfC(A,\dff_A)$, and its rows are short exact sequences of differential modules;
\item[{\rm (ii)}] $\bfC'$ is not isomorphic to $\bfC(B,\dff_B)$ for any $\sP$-graded differential module $(B,\dff_B)$;
\item[{\rm (iii)}] $\bfE(\bfC')\cong\bfE(A,\dff_A)$; in particular, $\bfE(\bfC')$ is chain generated.
\end{enumerate}
\end{theorem}
Theorem \ref{absnonreal}(ii) provides a negative answer to Question (Q) for $\bfC'$, and $\bfC'$ is not isomorphic to any homotopy-exact chain system.

\begin{proof}[of Theorem \ref{absnonreal}]
 (i): The conditions of Definition \ref{Csys} that do not involve $i'$ or $j'$ hold, since they hold for $\bfC(A,\dff_A)$. We check those that involve $i'$ or $j'$. Chain-homotopic chain maps
 are equal in $\DdR$. To show that $f\pmb{\sim}g$ for chain maps $f,g\colon X\to Y$, we show that the chain map $f-g$ is null-homotopic.
\begin{enumerate}[label=\textup{(\roman*)},leftmargin=2.2em,itemsep=0.1em,topsep=0pt]
\item[(a)] {\it Functoriality for $(\bt,\beta)\le(\alpha,\tp)$:} $i'j^\beta_\alpha\pmb{\sim} j^\tp_\alpha i^\tp_\beta$. The chain map $i'j^\beta_\alpha-j^\tp_\alpha i^\tp_\beta\colon F_\beta A\to G_{\tp\smin\alpha}A$ maps $\bm\xi\mapsto0$ and $\bm\eta\mapsto\varphi(\bm\eta)-\bm\eta=-z^2\bm\zeta$. Since $G_{\tp\smin\alpha}A$ has zero differential, it equals $h\,\dff_{F_\beta A}$ for $h(\bm\xi)=-\bm\zeta$ and $h(\bm\eta)=0$; indeed, $h\,\dff_{F_\beta A}(\bm\eta)=h(z^2\bm\xi)=-z^2\bm\zeta$.
\item[(b)] {\it Functoriality for $(\bt,\tp)\le(\beta,\tp)$:} $j'j^\tp_\alpha\pmb{\sim} j^\tp_\beta$. The chain map $j'j^\tp_\alpha-j^\tp_\beta\colon A\to G_{\tp\smin\beta}A$ maps $\bm\xi\mapsto0$, $\bm\eta\mapsto z^2\bm\zeta$ and $\bm\zeta\mapsto(1+z)\bm\zeta-\bm\zeta=z\bm\zeta$. Since $G_{\tp\smin\beta}A$ has zero differential, it equals $h\,\dff_A$ for $h(\bm\xi)=\bm\zeta$ and $h(\bm\eta)=h(\bm\zeta)=0$; indeed, $h\,\dff_A(\bm\eta)=h(z^2\bm\xi)=z^2\bm\zeta$ and $h\,\dff_A(\bm\zeta)=h(z\bm\xi)=z\bm\zeta$.
\item[(c)] {\it Naturality of $\sk$ for $(\bt,\alpha,\beta)\le(\bt,\alpha,\tp)$:} $k^\tp_\alpha i'=k^\beta_\alpha$, since $k^\tp_\alpha i'(\bm\eta)=k^\tp_\alpha(\bm\eta-z^2\bm\zeta)=-z^2\bm\xi+z^3\bm\xi=k^\beta_\alpha(\bm\eta)$.
\item[(d)] {\it Naturality of $\sk$ for $(\bt,\alpha,\tp)\le(\bt,\beta,\tp)$:} $k^\tp_\beta j'\pmb{\sim}(\sfS i^\beta_\alpha)k^\tp_\alpha$. The chain map $k^\tp_\beta j'-(\sfS i^\beta_\alpha)k^\tp_\alpha\colon G_{\tp\smin\alpha}A\to\sfS F_\beta A$ maps $\bm\eta\mapsto k^\tp_\beta(z^2\bm\zeta)+z^2\bm\xi=-z^3\bm\xi+z^2\bm\xi=z^2\bm\xi$ and $\bm\zeta\mapsto-z(1+z)\bm\xi+z\bm\xi=-z^2\bm\xi$. Since $G_{\tp\smin\alpha}A$ has zero differential and $\dff_{\sfS F_\beta A}\bm\eta=-z^2\bm\xi$, it equals $\dff_{\sfS F_\beta A}h$ for $h(\bm\eta)=-\bm\eta$ and $h(\bm\zeta)=\bm\eta$.
\item[(e)] {\it The row of $(\alpha,\beta,\tp)$ is distinguished:} since $\varphi\,i=i'$ and $j'\varphi=j$, the triple $(\id,\varphi,\id)$ is an isomorphism from the row $(i,j,k)$ of $\bfC(A,\dff_A)$, which is distinguished, onto $(i',j',k)$, cf.\ Defn.\ \ref{defndisttri}. For the same reason $0\to G_{\beta\smin\alpha}A\xrightarrow{i'}G_{\tp\smin\alpha}A\xrightarrow{j'}G_{\tp\smin\beta}A\to0$ is short exact.
\end{enumerate}
The remaining conditions either involve a term $C^\delta_\delta=0$, or follow from (e), since consecutive maps in a distinguished triangle compose to zero; for instance $j'i'=0$, which is functoriality for $(\alpha,\beta)\le(\alpha,\tp)\le(\beta,\tp)$. Hence $\bfC'$ is a chain system in $\DdR$, and its rows are short exact. It is excisive since $\sP$ is linear: any two down-sets are comparable, so that all excision maps are identities. The chain maps of $\bfC'$ do not form a functor $\IIi\to\CdR$, since $i'j^\beta_\alpha\neq j^\tp_\alpha i^\tp_\beta$, and by (ii) this cannot be repaired by passing to an isomorphic chain system.

{(ii):} \emph{Step 1.} There is no morphism $h\colon\bfE(A,\dff_A)\to\bfE(\bfC')$ of Cartan--Eilenberg systems with $h^\delta_\smbt=\id$ for all $\delta\in\sO(\sP)$. The terms $F_\alpha A$, $G_{\beta\smin\alpha}A$ and $G_{\tp\smin\alpha}A$ have zero differential, so that $E^\alpha_\smbt=\langle\bm\xi\rangle$, $E^\beta_\alpha=\langle\bm\eta\rangle$ and $E^\tp_\alpha=\langle\bm\eta,\bm\zeta\rangle$ for both Cartan--Eilenberg systems. Since $z^3=0$, the element $z\bm\eta\in F_\beta A$, and the elements $z^2\bm\zeta$ and $\bm\eta-z\bm\zeta$ of $A$, are cycles. Suppose $h$ exists, and write $h^\beta_\alpha(\bm\eta)=r\bm\eta$ and $h^\tp_\alpha(\bm\zeta)=p\bm\eta+q\bm\zeta$, with $r,p,q\in\bR$. The squares of \eqref{dia:exacttriangles} give:
\begin{enumerate}[label=\textup{(\roman*)},leftmargin=2.2em,itemsep=0.1em,topsep=0pt]
\item[(a)] for the map $j^\beta_\alpha$: $h^\beta_\alpha H(j^\beta_\alpha)=H(j^\beta_\alpha)$; at $[z\bm\eta]$ this reads $rz\bm\eta=z\bm\eta$, i.e.\ $z(r-1)=0$, and hence $rz^2=z^2$;
\item[(b)] for the map $i$ of $(\alpha,\beta,\tp)$: $h^\tp_\alpha H(i)=H(i')h^\beta_\alpha$; at $\bm\eta$ this gives $h^\tp_\alpha(\bm\eta)=r\,i'(\bm\eta)=r\bm\eta-z^2\bm\zeta$, by (a);
\item[(c)] for the map $j^\tp_\alpha$: $h^\tp_\alpha H(j^\tp_\alpha)=H(j^\tp_\alpha)$; at $[z^2\bm\zeta]$ this gives $z^2h^\tp_\alpha(\bm\zeta)=z^2\bm\zeta$, hence $z^2p=0$, and at $[\bm\eta-z\bm\zeta]$ it gives $h^\tp_\alpha(\bm\eta)-z\,h^\tp_\alpha(\bm\zeta)=\bm\eta-z\bm\zeta$, whose coefficient of $\bm\zeta$ reads $z^2+zq=z$, by (b);
\item[(d)] for the map $k^\tp_\alpha$: $H(k^\tp_\alpha)h^\tp_\alpha=H(k^\tp_\alpha)$; at $\bm\zeta$ this gives $-(z^2p+zq)\bm\xi=-z\bm\xi$, hence $zq=z$, by (c).
\end{enumerate}
By (c) and (d), $z^2=z-zq=0$, a contradiction.
\emph{Step 2.} Suppose that $\bfC'\cong\bfC(B,\dff_B)$ for some $\sP$-graded differential module $(B,\dff_B)$. By Lemma \ref{absequiv} we may assume that there is an isomorphism $f\colon\bfC(B,\dff_B)\to\bfC'$ whose components are $K$-projective resolutions; in particular, the differential modules $F_\delta B$ are $K$-projective. Since $\bfC'$ has the base of $\bfC(A,\dff_A)$, the components $f^\delta_\smbt\colon F_\delta B\to F_\delta A$ satisfy the hypotheses of Lemma \ref{absbaselemma}, which yields an isomorphism $g\colon\bfC(B,\dff_B)\to\bfC(A,\dff_A)$ with $g^\delta_\smbt=f^\delta_\smbt$ for all $\delta\in\sO(\sP)$. Then $f\circ g^{-1}\colon\bfC(A,\dff_A)\to\bfC'$ is an isomorphism of chain systems whose components on the base are identities, and $H(f\circ g^{-1})$ is a morphism as excluded in Step 1. Hence $\bfC'$ is not isomorphic to any $\bfC(B,\dff_B)$, and the remaining statements follow from Lemma \ref{absequiv}.

{(iii):} Let $\rmTheta$ be the automorphism of $F_\beta A$ given by $\rmTheta(\bm\xi)=(1+z)^{-1}\bm\xi$ and $\rmTheta(\bm\eta)=\bm\eta$; it is a chain map, since $z^2(1+z)^{-1}=z^2$. Define $h\colon\bfE(A,\dff_A)\to\bfE(\bfC')$ by the components
\[
h^\alpha_\smbt=(1+z)^{-1},\quad h^\beta_\smbt=H(\rmTheta),\quad h^\tp_\smbt=\id,\qquad h^\beta_\alpha=\id,\quad h^\tp_\alpha=\varphi,\quad h^\tp_\beta=\id,
\]
where $(1+z)^{-1}$ acts on $E^\alpha_\smbt=\langle\bm\xi\rangle$ by multiplication. All components are isomorphisms. We check the squares of \eqref{dia:exacttriangles} for the four triples of distinct down-sets; for the other triples they commute trivially.
\begin{enumerate}[label=\textup{(\roman*)},leftmargin=2.2em,itemsep=0.1em,topsep=0pt]
\item[(a)] $(\bt,\alpha,\beta)$: $\rmTheta\, i^\beta_\alpha=i^\beta_\alpha(1+z)^{-1}$, $j^\beta_\alpha\rmTheta=j^\beta_\alpha$ and $(1+z)^{-1}k^\beta_\alpha=k^\beta_\alpha$ hold for the chain maps, the last one since $(1+z)^{-1}z^2=z^2$.
\item[(b)] $(\alpha,\beta,\tp)$: $\varphi\,i=i'$, $j'\varphi=j$, and $k=0$.
\item[(c)] $(\bt,\alpha,\tp)$: $k^\tp_\alpha\varphi=(1+z)^{-1}k^\tp_\alpha$ holds for the chain maps, since $k^\tp_\alpha\varphi(\bm\eta)=-z^2\bm\xi+z^3\bm\xi=-z^2\bm\xi$. The chain maps $i^\tp_\alpha(1+z)^{-1}-i^\tp_\alpha$ and $\varphi\, j^\tp_\alpha-j^\tp_\alpha$ are null-homotopic: the first one maps $\bm\xi\mapsto(-z+z^2)\bm\xi$, and equals $\dff_Ah_1$ for $h_1(\bm\xi)=(z-1)\bm\zeta$; the second one maps $\bm\xi\mapsto0$, $\bm\eta\mapsto-z^2\bm\zeta$ and $\bm\zeta\mapsto(-z+z^2)\bm\zeta$, and equals $h_2\dff_A$ for $h_2(\bm\xi)=(z-1)\bm\zeta$, $h_2(\bm\eta)=h_2(\bm\zeta)=0$.
\item[(d)] $(\bt,\beta,\tp)$: the square for $j^\tp_\beta$ commutes trivially. The chain map $i^\tp_\beta\rmTheta-i^\tp_\beta$ maps $\bm\xi\mapsto(-z+z^2)\bm\xi=\dff_A\bigl((z-1)\bm\zeta\bigr)$ and $\bm\eta\mapsto0$, and $(\sfS\rmTheta)k^\tp_\beta-k^\tp_\beta$ maps $\bm\zeta\mapsto-z(1+z)^{-1}\bm\xi+z\bm\xi=z^2\bm\xi=\dff_{F_\beta A}\bm\eta$. Both map cycles to boundaries, and therefore induce zero on homology.
\end{enumerate}
Hence $h$ is an isomorphism of Cartan--Eilenberg systems. Since $\bfC(A,\dff_A)$ is homotopy-exact, $\bfE(\bfC')$ is chain generated, cf.\ Defn.\ \ref{def:chaingenerated}.
    \qed
\end{proof}
\begin{figure}[ht]
\centering
\begin{tikzpicture}[x=1cm,y=0.88cm,>={cm to},line width=0.4pt,
  every node/.style={inner xsep=1ex,inner ysep=0.85ex},
  lab/.style={font=\scriptsize,inner sep=0.5ex},
  face/.style={font=\small},
  conn/.style={->},
  row/.style={->,dash pattern=on 0.8ex off 0.4ex,dash phase=0.8ex}]
\node (Fb) at (0,5.3) {$F_\beta A$};
\node (Fa) at (-4.9,-3.2) {$F_\alpha A$};
\node (A)  at (4.9,-3.2) {$A$};
\node (Gba) at (-1.35,0.41) {$G_{\beta\smin\alpha}A$};
\node (G1b) at (1.35,0.41) {$G_{\tp\smin\beta}A$};
\node (G1a) at (0,-1.93) {$G_{\tp\smin\alpha}A$};
\draw[->] (Fa) -- node[lab,above left] {$i^\beta_\alpha$} (Fb);
\draw[->] (Fb) -- node[lab,above right] {$i^\tp_\beta$} (A);
\draw[->] (Fa) -- node[lab,below] {$i^\tp_\alpha$} (A);
\draw[->] (Fb) -- node[lab,left,pos=0.5] {$j^\beta_\alpha$} (Gba);
\draw[conn] (G1b) -- node[lab,right,pos=0.5] {$k^\tp_\beta$} (Fb);
\draw[conn] (Gba) -- node[lab,above left,pos=0.4] {$k^\beta_\alpha$} (Fa);
\draw[->] (A) -- node[lab,above right,pos=0.6] {$j^\tp_\beta$} (G1b);
\draw[conn] (G1a) -- node[lab,above,pos=0.42] {$k^\tp_\alpha$} (Fa);
\draw[->] (A) -- node[lab,above,pos=0.58] {$j^\tp_\alpha$} (G1a);
\draw[row] (Gba) -- node[lab,left] {$i$} (G1a);
\draw[row] (G1a) -- node[lab,right] {$j$} (G1b);
\draw[row] (G1b) -- node[lab,above] {$k$} (Gba);
\node[face] at (-2.25,1.05) {$\mathrm{d}$};
\node[face] at (2.25,1.05) {$\mathrm{d}$};
\node[face] at (0,-2.8) {$\mathrm{d}$};
\node[face] at (0,-0.4) {$\mathrm{d}$};
\node[face] at (0,2.2) {$\mathrm{c}$};
\node[face] at (-2.0,-1.45) {$\mathrm{c}$};
\node[face] at (2.0,-1.45) {$\mathrm{c}$};
\end{tikzpicture}
\caption{The octahedron of $\bfC(A,\dff_A)$ over $\sO(\sP)=\{\bt,\alpha,\beta,\tp\}$. The four faces marked $\mathrm{d}$ are distinguished triangles, the three faces marked $\mathrm{c}$ and the outer face are commutative, and the arrows $k^\beta_\alpha$, $k^\tp_\alpha$, $k^\tp_\beta$ and $k$ are connecting homomorphisms, whose targets are the suspensions of the indicated vertices. The dashed arrows form the row of $(\alpha,\beta,\tp)$. For the chain system $\bfC'$ of Theorem \ref{absnonreal} the inner edges $i$ and $j$ are replaced by $i'=\varphi\,i$ and $j'=j\,\varphi^{-1}$
\vspace{-0.9cm}}.
\label{octafig}
\end{figure}
\begin{remark}
\label{absnonrealrem}
The chain systems $\bfC(A,\dff_A)$ and $\bfC'$ are two Verdier octahedra in $\DdR$ on the same base, with the same bottom triangles, which are not isomorphic as chain systems,  cf.\ \cite{Kunzer}.
Their inner rows are isomorphic via $(\id,\varphi,\id)$, as provided by the axiom (TR3). By %but by Step 1 of
the proof of Theorem \ref{absnonreal}(ii) there is no isomorphism $\bfC(A,\dff_A)\to\bfC'$ that is the identity on the base. The obstruction is therefore not due to the difference between commutation in $\CdR$ and in $\DdR$, but to the non-uniqueness of the morphisms provided by (TR3).
\end{remark}
\begin{remark}
\label{absnonrealrem1}
By Theorem \ref{absnonreal}(ii), the isomorphism $h$ of Theorem \ref{absnonreal}(iii) cannot be the identity on the base. Moreover, its base components $(1+z)^{-1}$, $H(\rmTheta)$, $\id$ are not induced by an automorphism of the base of $\bfC(A,\dff_A)$ in $\DdR$. Otherwise Lemma \ref{absbaselemma}, with $B=A$,   extends the automorphism to an automorphism $f$ of $\bfC(A,\dff_A)$, and the choice $h\circ H(f)^{-1}$   contradicts the first step of the proof in Theorem \ref{absnonreal}(ii). Lemma \ref{absbaselemma} applies, since the $F_\delta A$ are $K$-projective, cf.\ Props.\ \ref{Kprojprop} and \ref{prop:kproj}. In particular, $i^\tp_\beta\rmTheta$ and $i^\tp_\beta$ induce the same map on homology, but they are not chain-homotopic.
\end{remark}
\begin{remark}
\label{absnonrealrem2}
The question of Remark \ref{notcgsys}, whether $\bfE(\bfC)$ is chain generated for every excisive chain system $\bfC$ in $\DdR$, remains open.
For $\bfC'$ the answer is affirmative, by Theorem \ref{absnonreal}(iii). By (ii), a proof cannot follow the strategy of the proof of Theorem \ref{mainres987655}, which first realizes the base in $\DdR$ and then adjusts the other terms.
\end{remark}

\section*{Disclosure of AI use}
The author used Claude (Anthropic) Opus 5.5 to assist with the literature search, construct the example in Section~\ref{dynexample}, the counterexample in Section~\ref{csysnothe}, and review the manuscript for mathematical errors and notational inconsistencies. All derivations, proofs, and citations were independently verified by the author, who takes full responsibility for the content.

 \begin{acknowledgements}
 The author wishes to thank Kelly Spendlove for many valued discussions on the subject.
 \end{acknowledgements}

\bibliographystyle{spmpsci}      % mathematics and physical sciences
\bibliography{references}   % name your BibTeX data base

\appendix

\section{Rotating distinguished triangles and $K$-projective cones}
\label{rottrian}
For a chain map $f\colon B\to A$ we have the canonical short exact sequence:
\[
\begin{tikzcd}
0 \arrow[r] & A \arrow[r, "i", tail] & \cone(f) \arrow[r, "j", two heads] & \sfS B \arrow[r] & 0,
\end{tikzcd}
\]
which is split as a sequence of $\bR$-modules, cf.\ \eqref{SESforcone}. Following Section \ref{mappingconesect}, it induces the triangle
\begin{equation}
    \label{induceddtr12}
\begin{tikzcd}
 A \arrow[r, "i"] & \cone(f) \arrow[r, "j"] & \sfS B  \arrow[r, "k"]  & \sfS A
\end{tikzcd}
\end{equation}
with $k=j^\dagger\rmTheta^{-1}$, where $\rmTheta\colon\cone(i)\to\sfS B$, $\rmTheta(y,a)=j(y)$, and {the projection} $j^\dagger\colon\cone(i)\to\sfS A$ are {the maps $\rmTheta$ and $j$ of} \eqref{constructwithr} {for $r=i$ and $s=j$}. {The map} $\rmTheta$ is a  chain-equivalence and \eqref{induceddtr12} is a distinguished triangle in $\KdR$, as stated in Remark \ref{sometimesinKdR}. The following proposition proves this directly and computes $k$; for chain complexes, cf.\ \cite[Prop.~10.2.4]{weibel}.

\begin{proposition}
    \label{firstshift}
    The map $\rmTheta$ is a  chain-equivalence, and in the distinguished triangle \eqref{induceddtr12}, $k=-\sfS f$.
\end{proposition}
\begin{proof}
Let $i^\dagger\colon\cone(f)\to\cone(i)$ and $j^\dagger\colon\cone(i)\to\sfS A$ be the standard inclusion and projection, so that
\[
\begin{tikzcd}
A \arrow[r, "i"] & \cone(f) \arrow[r, "i^\dagger"] & \cone(i) \arrow[r, "j^\dagger"] & \sfS A
\end{tikzcd}
\]
is the standard triangle of $i\colon A\to\cone(f)$. By \eqref{dagmapcone}, $\cone(i)=\cone(f)\oplus A=A\oplus B\oplus A$ with differential
\[
\dff_i = \begin{pmatrix}
    \dff_A~& f~ & \id_A\\
    0~& -\dff_B~& 0\\
    0~ & 0~& -\dff_A
\end{pmatrix},
\]
and $\rmTheta(a,b,a')=b$. By definition $\rmTheta i^\dagger=j$, which yields the commutative diagram
\[
\begin{tikzcd}[column sep=large, row  sep=large]
A \arrow[r, "i"] \arrow[equal]{d} & \cone(f) \arrow[r, "j"] \arrow[equal]{d} & \sfS B \arrow[r, "k"]  & \sfS A \arrow[equal]{d} \\
A \arrow[r, "i"]                    & \cone(f) \arrow[r, "i^\dagger"]                    & \cone(i) \arrow[u, "\rmTheta"]\arrow[r, "j^\dagger"]    & \sfS A
\end{tikzcd}
\]
 cf.\ \eqref{constructwithr}. Define the map $\chi\colon\sfS B\to\cone(i)$ by $\chi(b)=\bigl(0,b,-f(b)\bigr)$, and since $f$ is a chain map,
\[
\dff_i\chi(b)=\bigl(0,-\dff_B b,\dff_A f(b)\bigr)=\chi(-\dff_B b)=\chi(\dff_{\sfS B}b),
\]
which shows that $\chi$ is a chain map. Note that $\rmTheta\chi=\id_{\sfS B}$. With the map $h\colon\cone(i)\to\cone(i)$, given by $h(a,b,a')=(0,0,a)$, we have
\[
\begin{aligned}
(\dff_ih+h\dff_i)(a,b,a') &=(a,0,-\dff_Aa)+\bigl(0,0,\dff_Aa+f(b)+a'\bigr)\\
&=\bigl(a,0,a'+f(b)\bigr)=(\id-\chi\rmTheta)(a,b,a'),
\end{aligned}
\]
which proves that $\rmTheta$ is a  chain-equivalence with homotopy inverse $\chi$. In the homotopy category $\KdR$
\[
k=j^\dagger\rmTheta^{-1}=j^\dagger\chi=-\sfS f,
\]
using that $j^\dagger\chi(b)=-f(b)$ and $\sfS f=f$ as maps.
\qed
\end{proof}

The triangle \eqref{induceddtr12} with $k=-\sfS f$ is the first rotation, cf.\ (TR2), of the standard triangle
\[
\begin{tikzcd}
B \arrow[r, "f"] & A \arrow[r, "i"] & \cone(f) \arrow[r, "j"] & \sfS B,
\end{tikzcd}
\]
and by Proposition \ref{firstshift} it is isomorphic in $\KdR$ to the standard triangle of $i$. This is the direct argument for \eqref{firstrot12}. Since $\sfS^2=\id$, rotating six times returns the standard triangle, i.e.\ rotation is periodic with period $6$.

Mapping cones of chain maps between $K$-projective differential modules are again $K$-projective,
cf.\ Prop.\ \ref{prop:kproj}.
\begin{proposition}
    \label{Kprojcones}
    Let $f\colon B\to A$ be a chain map between $K$-projective differential modules. Then $\sfS B$ and $\cone(f)$ are $K$-projective.
\end{proposition}
\begin{proof}
Let $(N,\dff_{N})$ be acyclic.  A homomorphism $h\colon  \sfS B\to N$ is a chain map if and only if   {$h\colon B\to\sfS N$ is a chain map}, since both conditions read $\dff_{N} h + h\dff_B=0$. Because $\sfS N$ is acyclic and $B$ is $K$-projective, $h\pmb{\sim} 0$ as a chain map $B\to\sfS N$, i.e. $h=\dff_{\sfS N}g+g\dff_B=\dff_{N}(-g)+(-g)\dff_{\sfS B}$ for a homomorphism $g$, and thus $h\pmb{\sim} 0$ as a chain map $\sfS B\to {N}$, {with null-homotopy $-g$}.
Consequently, $\sfS B$ is $K$-projective.

Let $h\colon\cone(f)\to {N}$ be a chain map, and let
$i$ and $j$ be as in \eqref{induceddtr12}.
Then, the  composition $hi\colon A\to {N}$ is a chain map and since
  $A$ is $K$-projective, $hi=\dff_{{N}}g+g\dff_A$ for a homomorphism $g\colon A\to {N}$. Define $G\colon\cone(f)\to {N}$ by $G(a,b)=g(a)$. Then, $h'=h-\dff_{N}G-G\dff_f$ is a chain map with $h'{\pmb{\sim}} h$.  From the differential on $\cone(f)$ we have that $\dff_f(a,0)=(\dff_Aa,0)$, and thus
\[
\begin{aligned}
h'(a,0)&= h(a,0)-\dff_{N} G(a,0) - G\dff_f(a,0)\\
&=hi(a)-\dff_{N}g(a)-g(\dff_Aa)=0.
\end{aligned}
\]
Consequently, $h'$ vanishes on $\ker j= \image i$, and $h'=sj$ for a homomorphism $s\colon\sfS B\to N$. The latter is a chain map since $j$ is a surjective chain map. The fact that  $\sfS B$ is $K$-projective implies that $s\pmb{\sim} 0$ and therefore $s=\dff_{N} r+ r\dff_{\sfS B}$ for some homomorphism $r\colon \sfS B\to N$. This yields, using $j\dff_f=\dff_{\sfS B}j$,
\[
h'=sj = (\dff_{N}r + r\dff_{\sfS B})j = \dff_{N}(rj) + (rj)\dff_f.
\]
Hence, $h'\pmb{\sim} 0$, and thus $h {\pmb{\sim}} 0$, which proves that $\cone(f)$ is $K$-projective.
    \qed
\end{proof}
Proposition \ref{Kprojcones} yields Proposition \ref{prop:kproj}. Closure under $\sfS$ is the first statement, and $\sfS^{-1}=\sfS$, since $\sfS^2=\id$. Closure under mapping cones is the second statement. A   direct sum $P\oplus P'$ of $K$-projective differential modules is the mapping cone of the zero map $\sfS P'\to P$, and is therefore $K$-projective.

Direct summands of $K$-projective differential modules are $K$-projective. Let $P=P'\oplus P''$ as differential modules, i.e.\ $\dff_P=\dff_{P'}\oplus\dff_{P''}$, and assume that $P$ is $K$-projective. The inclusion $i\colon P'\to P$ and the projection $j\colon P\to P'$ are chain maps with $ji=\id_{P'}$. Let $h'\colon P'\to N$ be a chain map, with $N$ acyclic.  Then, $h'j\colon P\to N$, which is $h'$ extended by zero on $P''$, is a chain map. Since $P$ is $K$-projective, $h'j=\dff_{N}g+g\dff_P$ for some homomorphism $g\colon P\to N$. Composing with  $i$ on the right, and using $\dff_Pi=i\dff_{P'}$, yields
\[
h'=h'ji=\dff_{N}(gi)+(gi)\dff_{P'}.
\]
This proves that  $h'\pmb{\sim} 0$, with null-homotopy $gi$,  and thus $P'$ is $K$-projective.

Finally, $K$-projectivity of $P$ is the condition $\Hom_{\KdR}(P,N)=0$ for all acyclic $N$, which is invariant under isomorphisms in $\KdR$.   {If two of the objects of a distinguished triangle in $\KdR$ are $K$-projective, then so is the third: by (TR2), and since $\sfS$ preserves $K$-projectivity, we may assume that these are the first two; the triangle is isomorphic in $\KdR$ to a standard triangle $\bigl(B,A,\cone(f),f,i,j\bigr)$, cf.\ Defn.\ \ref{defndisttri}, so that $B$ and $A$ are $K$-projective, and the third object is isomorphic to $\cone(f)$, which is $K$-projective by Proposition \ref{Kprojcones}.} Hence the $K$-projective differential modules form a thick triangulated subcategory of $\KdR$, cf.\ \cite[Prop.~1.3]{Spalt} for complexes. If, moreover, $A$ and $B$  in Proposition \ref{Kprojcones} are projective modules, then so is $\cone(f)=A\oplus B$.

\section{Projective flags and $K$-projective resolutions}
\label{resolfordiffmod}
Let $\N^* = \N\cup \{\infty\}$. A  \emph{projective flag} is an $\N^*$-filtered
differential module $(P,\dff_P)$,
with  strict differential and projective graded module ${\rm gr} P=\bigoplus_{k\ge 1}F_kP/F_{k-1}P$.  Unpacking the hypotheses
translates to having a filtering  by submodules
\[
  0= F_0P\ \subseteq\ F_1P\ \subseteq\ F_2P\ \subseteq\ F_3P\ \subseteq\ \cdots
  \ \subseteq\ F_\infty P = P,
\]
satisfying
\begin{enumerate} [label=\textup{(\roman*)},leftmargin=2.2em,itemsep=0.1em,topsep=0pt]
\item [(i)]  $\displaystyle\bigcup_{k\ge 0}F_kP=P$;
\item [(ii)] $\dff_P F_kP\subseteq F_{k-1}P$ for every $k\ge 1$;
\item [(iii)]  $G_kP=F_kP/F_{k-1}P$ is a projective $\bR$-module for every
      $k\ge 1$.
\end{enumerate}
Since the modules $G_kP$ are projective, the inclusions $F_{k-1}P\subseteq F_kP$ split as $\bR$-linear maps, i.e.\ $F_kP=F_{k-1}P\oplus Q_k$, with $Q_k\cong G_kP$, and by (i), $P=\bigoplus_{k\ge 1}Q_k\cong\bigoplus_{k\ge 1}G_kP$ as $\bR$-modules. In particular, the modules $F_kP$ and $P$ are projective.
If  all modules $G_kP$ are free, then so are the modules $F_kP$ and $P$, and
$(P,\dff_P)$ is called a \emph{free flag}, cf.\ \cite[Sect.~2]{ABI2007} and \cite[Sect.~2.2]{BrEr}.
\begin{definition}\label{def:semifree}
A \emph{projective flag resolution} of a differential module $A$
is a quasi-isomorphism
\[
  \pi\colon (P,\dff_P)\longrightarrow(A,\dff_A),
\]
with $(P,\dff_P)$ a projective  flag. A \emph{free flag resolution} requires $P$ to be a free flag.
\end{definition}

The $K$-projective resolutions play the same role for differential modules as projective resolutions do for $R$-modules.

\begin{proposition}
    \label{Kprojprop}
    Let $(P,\dff_P)$ be a (free) projective flag. Then, $P$ is $K$-projective.
\end{proposition}
\begin{proof}
Let $f\colon P \to N$ be a chain map with $(N,\dff_{N})$ acyclic.  To prove that $P$ is $K$-projective we construct a module homomorphism $h\colon P \to N$ such that
$f=h\dff_P + \dff_{N} h$. Define $h_0=0$ as map $h_0\colon F_0P\to N$. Suppose $h_{k-1}\colon F_{k-1}P\to N$ is given and which satisfies
\[
f|_{F_{k-1}P} =h_{k-1}\dff_P + \dff_{N} h_{k-1}.
\]
Consider the short exact sequence of differential modules:
\[
\begin{tikzcd}[column sep=large, row sep=large]
0 \arrow[r] & F_{k-1} P \arrow[r, "\subset"] & F_kP \arrow[r, "\pi"] & G_k P \arrow[r] & 0
\end{tikzcd}
\]
By assumption $G_k P$ is projective and therefore $F_kP = F_{k-1}P\oplus Q_k$, with $Q_k\cong G_kP$. Define $h_k|_{F_{k-1}P} = h_{k-1}$ and $h_k|_{Q_k}$ to be constructed.
Let $q\in Q_k$. Since the differential $\dff_P$ is strict it holds that $\dff_P q\in F_{k-1} P$. Moreover, $h_{k-1}$ is defined by the induction hypothesis, and so is $h_{k-1}(\dff_Pq)$. Consider the map $\psi\colon Q_k \to N$ defined by
\[
\psi(q) = f(q) - h_{k-1}(\dff_Pq).
\]
By the induction hypothesis and the fact that $f$ is a chain map we obtain:
\[
\begin{aligned}
\dff_{N} \psi(q) &= \dff_{N} f(q) - \dff_{N} h_{k-1}(\dff_Pq) = f(\dff_P q) - \dff_{N} h_{k-1}(\dff_Pq)\\
&= f(\dff_P q) - \bigl[ f(\dff_P q) - h_{k-1}(\dff_P^2q)\bigr] = h_{k-1}(\dff_P^2q)=0.
\end{aligned}
\]
The differential module $N$ is acyclic and therefore $\ker \dff_{N} = \image\dff_{N}$.
This implies that $\psi(q) \in \ker \dff_{N} = \image\dff_{N}$ and  {therefore} $\psi\colon Q_k \to \image \dff_{N}$. The differential $\dff_{N}\twoheadrightarrow \image \dff_{N}$ is surjective. Since $Q_k$ is projective,  there exists a lift $h_k|_{Q_k}\colon Q_k \to N$ with $\dff_{N} h_k = \psi$:
\[
\begin{tikzcd}[column sep=large, row sep=large]
                                                & N \arrow[d, "\dff_{N}", two heads] \\
Q_k \arrow[ru, "h_k", dotted] \arrow[r, "\psi"] & \image \dff_{N}
\end{tikzcd}
\]
For $p+q\in F_{k-1}P \oplus Q_k$ this yields:
\[
\begin{aligned}
    f(p+q) &= f(p) + f(q) = h_{k-1}(\dff_P p) + \dff_{N} h_{k-1}(p)
    + h_{k-1}(\dff_P q) + \psi(q)\\
    &= h_{k-1}(\dff_P (p+q)) + \dff_{N} h_{k-1}(p) + \dff_{N} h_{k}(q)\\
    &= h_{k}(\dff_P (p+q)) + \dff_{N} h_{k}(p) + \dff_{N} h_{k}(q)\\
    &= h_{k}(\dff_P (p+q)) + \dff_{N} h_{k}(p+q),
\end{aligned}
\]
which proves that $f|_{F_{k}P} =h_{k}\dff_P + \dff_{N} h_{k}$. There exists a unique module homomorphism $h\colon P\to N$ such that $h|_{F_kP} = h_k$, which satisfies
\[
f =h \dff_P + \dff_{N} h.
\]
The latter establishes $f\pmb{\sim} 0$, proving that $P$ is $K$-projective.
    \qed
\end{proof}

\begin{proposition}\label{prop:semifreeexist}
Every differential module $(A,\dff_A)$ admits a (free) projective flag resolution.
\end{proposition}
\begin{proof}
Put $A_0=A$, and suppose the differential module $(A_k,\dff_{A_k})$ is defined.  Choose a set $Z_k$ of cycles in $A_k$ whose classes generate $H(A_k)$, let $P_k=\bigoplus_{z\in Z_k}Re_z$ be the free module with zero differential, and let $f_k\colon(P_k,0)\to A_k$ be the chain map with $f_k(e_z)=z$. Then $H(f_k)\colon P_k\to H(A_k)$ is surjective. Define $A_{k+1}=\cone(f_k)=A_k\oplus P_k$, with differential $\dff_{A_{k+1}}(x,p)=\big(\dff_{A_k}x+f_k(p),0\big)$, cf.\ \eqref{dagmapcone}, and let $i\colon A_k\to A_{k+1}$ be the inclusion. Since $H(f_k)$ is surjective, the long exact sequence \eqref{LESinHom}
\[
P_k\xrightarrow[]{H(f_k)}H(A_k)\xrightarrow[]{H(i)}H(A_{k+1})\to P_k
\]
shows that $H(i)=0$.

The differential modules $A_k$ form an increasing sequence of differential submodules, whose union is
\[
A_\infty=\bigcup_{k\ge 0}A_k=A\oplus Q,\quad Q=\bigoplus_{k\ge 0}P_k.
\]
The differential module $A_\infty$ is acyclic: a cycle $x\in A_\infty$ lies in some $A_k$, and since $H(i)=0$, $x$ is a boundary in $A_{k+1}$. Since $A$ is a differential submodule of $A_\infty$, the differential of $A_\infty$ is of the form
\[
\dff_{A_\infty}(a,q)=\big(\dff_Aa+\pi(q),-\dff_Pq\big),
\]
for homomorphisms $\pi\colon Q\to A$ and $\dff_P\colon Q\to Q$. From $\dff_{A_\infty}^2=0$ we obtain $\dff_P^2=0$ and $\dff_A\pi=\pi\dff_P$. Hence $P=(Q,\dff_P)$ is a differential module, $\pi\colon P\to A$ is a chain map, and $A_\infty=\cone(\pi)$, cf.\ \eqref{dagmapcone}. Since $A_\infty$ is acyclic, $\pi$ is a quasi-isomorphism by Lemma \ref{qisviacone}.

Finally, let $F_0P=0$ and $F_mP=P_0\oplus\cdots\oplus P_{m-1}$, $m\ge 1$. For $p\in P_k$ we have $\dff_{A_\infty}(0,p)=f_k(p)\in A_k=A\oplus P_0\oplus\cdots\oplus P_{k-1}$, and thus $\dff_PF_mP\subset F_{m-1}P$. The quotients $F_mP/F_{m-1}P=P_{m-1}$ are free, and the filtering is exhaustive. Hence $P$ is a free flag, and $\pi\colon P\to A$ is a free flag resolution.
    \qed
\end{proof}

In some cases projective flag resolutions can be constructed that have finiteness properties related to the differential module $(A,\dff_A)$.
\begin{proposition}\label{prop:semifreeexistfinite}
Every differential module $(A,\dff_A)$,
with the property that $H(A)$ allows a finite resolution  consisting of finitely generated projective modules,
admits a finite, projective flag resolution consisting of finitely generated projective modules.
\end{proposition}

\begin{proof}
Let $0\to P_n\xrightarrow[]{\partial_n}\cdots\xrightarrow[]{\partial_1}P_0\xrightarrow[]{\epsilon}H(A)\to 0$ be a resolution of $H(A)$ by finitely generated projective modules, and put $\partial_0=\epsilon$, $K_0=H(A)$ and $K_k=\ker\partial_{k-1}$, $1\le k\le n+1$, so that $\partial_k\colon P_k\to K_k$ is surjective with kernel $K_{k+1}$, and $K_{n+1}=0$. We follow the construction in the proof of Proposition \ref{prop:semifreeexist}, with the modules $P_k$ of the resolution, with zero differential, and with $f_k\colon(P_k,0)\to A_k$ a lift to $\ker\dff_{A_k}$ of $\kappa_k\partial_k\colon P_k\to H(A_k)$, which exists since $P_k$ is projective, where $\kappa_k\colon K_k\to H(A_k)$ is an isomorphism, and $\kappa_0=\id$. Then $H(f_k)=\kappa_k\partial_k$ is surjective with kernel $K_{k+1}$, and by \eqref{LESinHom} for $A_{k+1}=\cone(f_k)$, $H(i)=0$ and $H(j)\colon H(A_{k+1})\to P_k$ is injective with image $\ker H(f_k)=K_{k+1}$, where $j\colon A_{k+1}\to\sfS P_k$ is the projection; put $\kappa_{k+1}=H(j)^{-1}$. Since $K_{n+1}=0$, $A_{n+1}$ is acyclic, and the proof of Proposition \ref{prop:semifreeexist}, with $A_{n+1}$ in place of $A_\infty$, yields the finite projective flag $P=P_0\oplus\cdots\oplus P_n$ of finitely generated projective modules and the quasi-isomorphism $\pi\colon P\to A$. By Proposition \ref{Kprojprop}, $\pi$ is a $K$-projective resolution.
    \qed
\end{proof}

\end{document}